%% file: main.tex
\documentclass[11pt]{article}

\usepackage{graphicx}
\usepackage{latexsym,amsmath,amsfonts,amscd, amsthm, dsfont}
\usepackage{bm,color}
\usepackage{epsfig,verbatim,epstopdf,graphics}
\usepackage{subfigure}
\usepackage{changebar}
\usepackage{multirow}

\usepackage[ruled,vlined]{algorithm2e}

\usepackage{yhmath}%�ӻ���ͷ��������  \wideparen
 \usepackage{booktabs} %\cmidrule(lr){2-5} \cmidrule(lr){6-9}
 \usepackage{tikz}
\usepackage{verbatim}
\usetikzlibrary{arrows,backgrounds,snakes,shapes}

 \numberwithin{equation}{section}

\graphicspath{{./}{./figure/}}
\allowdisplaybreaks

\newtheoremstyle{plainNoItalics}{}{}{\normalfont}{}{\bfseries}{.}{ }{}

\theoremstyle{plain}
\newtheorem{thm}{Theorem}[section]

\theoremstyle{plainNoItalics}

\newtheorem{exa}[thm]{Example}

\newcommand{\bx}{{\bf x}}

\newcommand{\bv}{{\bf v}}
\newcommand{\bB}{{\bf B}}
\newcommand{\bE}{{\bf E}}
\newcommand{\bJ}{{\bf J}}

\newcommand{\ba}{{\bf a}}
\newcommand{\bU}{{\bf U}}

\newcommand{\br}{{\bf r}}

\newcommand{\beq}{\begin{equation}}
\newcommand{\eeq}{\end{equation}}
\newcommand{\bit}{\begin{itemize}}
\newcommand{\eit}{\end{itemize}}
\newcommand{\be}{\begin{eqnarray}}
\newcommand{\ee}{\end{eqnarray}}
\newcommand{\beno}{\begin{eqnarray*}}
\newcommand{\eeno}{\end{eqnarray*}}

\makeatletter

\newcommand{\Rmnum}[1]{\expandafter\@slowromancap\romannumeral #1@}
\makeatother

\begin{document}

\baselineskip=1.8pc

%\vspace*{.10in}

%=============  title  =========================
\input{title}
\input{intro}

\input{algorithm}
\input{numerical}

\input{conclusion}

%\appendix
%\input{appendix}

\bibliographystyle{abbrv}
\bibliography{refer}

\end{document}

%% file: title.tex
%!TEX root = main.tex
\begin{center}
{\bf
A Low Rank Tensor Representation of Linear Transport and Nonlinear Vlasov Solutions and Their Associated Flow Maps
}
\end{center}

\vspace{.2in}
\centerline{
 Wei Guo\footnote{
Department of Mathematics and Statistics, Texas Tech University, Lubbock, TX, 70409. E-mail:
weimath.guo@ttu.edu. Research
is supported by NSF grant NSF-DMS-1620047 and NSF-DMS-1830838.
},
Jing-Mei Qiu\footnote{Department of Mathematical Sciences, University of Delaware, Newark, DE, 19716. E-mail: jingqiu@udel.edu. Research supported by NSF grant NSF-DMS-1818924, Air Force Office of Scientific Research FA9550-18-1-0257.}
}

\bigskip
\noindent
{\bf Abstract.} We propose a low-rank tensor approach to approximate linear transport and nonlinear Vlasov solutions and their associated flow maps. The approach takes advantage of the fact that the differential operators in the Vlasov equation is tensor friendly,  based on which we propose a novel way to dynamically and adaptively build up low-rank solution basis by adding new basis functions from discretization of the PDE, and removing basis from an SVD-type truncation procedure. For the discretization, we adopt a high order finite difference spatial discretization and a second order strong stability preserving multi-step time discretization. We apply the same procedure to evolve the dynamics of the flow map in a low-rank fashion, which proves to be advantageous when the flow map enjoys the low rank structure, while the solution suffers from high rank or displays filamentation structures. Hierarchical Tucker decomposition is adopted for high dimensional problems. An extensive set of linear and nonlinear Vlasov test examples are performed to show the high order spatial and temporal convergence of the algorithm with mesh refinement up to SVD-type truncation, the significant computational savings of the proposed low-rank approach especially for high dimensional problems, the improved performance of the flow map approach for solutions with filamentations.

\vfill

\newpage

%% file: intro.tex
%!TEX root = main.tex
\section{Introduction}

In this work, we propose a novel approach aiming for {resolution of} the challenge of curse of dimensionality in grid-based methods when solving high dimensional nonlinear Vlasov equation as a kinetic description of collisionless plasma. 
The nonlinear Vlasov equation models a collisionless plasma, via distribution functions $f_s(\bx,\bv,t)$ for particle species $s$ with charge $q_s$ and mass $m_s$. The complete Vlasov-Maxwell model reads as follows:
\beq
\frac{\partial f_s}{\partial t}  + {\bf v} \cdot \nabla_{\bf x} f_s + \frac{q_s}{m_s}({\bf E} + {\bf v} \times {\bf B}) \cdot \nabla_{\bf v} f_s = 0,
\label{eq: vlasov}
\eeq
\begin{align}
&-\frac{1}{c^2} \frac{\partial  \bE}{\partial t} +\nabla \times \bB = \mu_0  \bJ,\quad  \frac{\partial  \bB}{\partial t} +\nabla \times \bE =  0,\label{eq:maxwell}\\
&\nabla \cdot \bE =\frac{\rho}{\varepsilon_0},  \quad\nabla \cdot \bB = 0,\notag
\end{align}
where $c$ is the speed of light, 
$\varepsilon_0$ and $\mu_0$ are the vacuum permittivity and
permeability, respectively. $\mathbf{E}$ and $\mathbf{B}$ are the electric and magnetic fields. The sources for Maxwell's equations \eqref{eq:maxwell}, i.e., the macroscopic charge density $\rho$ and the current density $\bJ$, are obtained from the distribution functions $f_s$:
\begin{equation}\label{eq:macro}
\rho=\sum_{s} q_{s} \int f_{s} \mathrm{~d} \bv, \quad \bJ=\sum_{s} q_{s} \int f_{s} \bv \mathrm{d} \bv.
\end{equation}
Observe that the Vlasov equation \eqref{eq: vlasov} is a six-dimensional nonlinear transport equation in phase space.
Among many existing challenges for deterministic Vlasov simulations (e.g. multiscale features, nonlinearity, formation of filamentation structures), the curse of dimensionality and the associated huge computational cost have been a key obstacle for realistic high-dimensional simulations. 
Note that the celebrated particle-in-cell method can generate qualitative results with reasonable computational cost in high dimensions, while the inherent statistical noise of such a method prevents accurate capture of physics of interest  \cite{kraus2017gempic}.
Hence, most existing deterministic schemes are only applicable to reduced lower-dimensional models ($d\le3$) in the literature. The sparse grid approach is considered as a viable framework for dimension reduction in the Vlasov simulations, yet the computation can still be very expensive for large $d$ (e.g., $d=6$) as the curse of dimensionality is not be fully removed \cite{bungartz2004sparse,kormann2016sparse,tao2018sparseguo}.
Another related approach is the reduced order modeling (ROM) \cite{benner2015survey}: typically a low-dimensional reduced subspace is constructed in an offline training phase for approximating the solution manifold. Then the surrogate solution for any desired parameter can be computed very efficiently from the reduced model in the online phase. The nonlinear Vlasov models are hyperbolic in nature, thus may not have low-rank/low-dimensional structures, if snapshot of solutions are taken at different instances of time, as opposed to parabolic problems in a reduced order modeling framework. It is related to the slow decay of the Kolmogorov $N$-width of the solution manifold for transport-dominated problems \cite{greif2019decay}. 

%Motivated by the emerging low-rank tensor decomposition techniques, there have been a few pioneering works in exploring the low-rank  solution structure of the Vlasov equation. These include the low-rank semi-Lagrangian (SL) method in the TT format developed in \cite{kormann2015semi};  a low-rank method based on the canonical polyadic (CP) format developed in \cite{ehrlacher2017dynamical}; the method proposed in \cite{einkemmer2018low,einkemmer2020low} in seeking a set of dynamic low-rank bases by a tangent space projection. Moreover, in \cite{hatch2012analysis} the HOSVD is applied to analyze and compress high-dimensional gyrokinetic datasets  generated by a full-rank spectral method, leading to efficient data compression especially in velocity domain. In \cite{dektor2020dynamically, dektor2021dynamic}, dynamic tensor approximations for high dimensional linear and nonlinear PDEs are proposed based on functional tensor decomposition and dynamic tensor approximation. There are recent work of on low-rank methods with asymptotic preserving property for multi-scale models \cite{einkemmer2021asymptotic, einkemmer2021efficient, chen2020random}.  

Inspired by the existing understanding of the low-rank solution structure for the Vlasov dynamics, as well as the observation that the differential operator in the Vlasov equation \eqref{eq: vlasov} can be represented in the tensorized form, in this paper we consider a novel way to (a) dynamically and adaptively build up low-rank solution basis, and (b) determine the low-rank solutions in a tensor format with well-established high order finite difference upwind weighted essentially non-oscillatory (WENO) method coupled with the SSP multi-step time discretizations \cite{gottlieb2011strong}, which offers more computational savings compared with the SSP multi-stage RK method in the low-rank tensor framework. We will first demonstrate our proposed idea for a reduced dimensionless 1D1V Vlasov-Poisson system. Compared with the recent work on dynamic tensor approximations with constant rank such as \cite{einkemmer2018low,einkemmer2020low, dektor2020dynamically, dektor2021dynamic}, our proposed approach is based on a procedure of adding from RHS of PDEs and removing basis by SVD truncation; hence not only the basis, but also the rank of the solution are dynamically evolving. 
Motivated from the filamentation phenomenon of the Vlasov solution, we propose a low-rank approach to evolve the flow map of solution, followed by fetching solution values at the feet of characteristics. Such an approach displays advantages for problems whose flow maps are of low rank, yet their solutions are not necessarily of low rank. However, further development is needed in a more general setting, e.g. in handling general boundary conditions. 
Then we discuss the extension to general high-dimensional cases in light of the Hierarchical Tucker (HT) decomposition. 
The HT format \cite{hackbusch2009new,grasedyck2010hierarchical} is motivated by the classical Tucker  format (also known as the tensor subspace format) \cite{tucker1966some,de2000multilinear}. It is developed by considering a dimension tree and taking advantage of the hierarchy of the nested subspaces and associated nested basis. A \emph{quasi-optimal} low-rank approximation in the HT format can be computed stably via the hierarchical  high order singular value decomposition (HOSVD)\cite{hackbusch2009new,grasedyck2010hierarchical,hackbusch2012tensor}. The HT format attains a storage complexity that is linearly scaled with the dimension,  hence striking a perfect balance between data complexity and numerical feasibility. We note that an alternative way of representing low-rank tensor is via the tensor train (TT) format \cite{oseledets2011tensor}, which can be thought of as a special type of the HT format, which has a degenerate dimension tree and enjoys a simpler structure. In this paper, we focus on the HT format, with a balanced dimension tree that separate the physical space ${\bf x}$ and phase space ${\bf v}$ dimensions. We also solve Poisson's equation in the HT format by adopting a low-rank conjugate gradient method \cite{grasedyck2018distributed}.

There have been a few pioneering works in exploring the low-rank  solution structure of the Vlasov equation with tensor decompositions. These include the low-rank semi-Lagrangian (SL) method in the TT format developed in \cite{kormann2015semi};  a low-rank method based on the canonical polyadic (CP) format developed in \cite{ehrlacher2017dynamical}; the method proposed in \cite{einkemmer2018low,einkemmer2020low} in seeking a set of dynamic low-rank bases by a tangent space projection. Moreover, in \cite{hatch2012analysis} the HOSVD is applied to analyze and compress high-dimensional gyrokinetic datasets  generated by a full-rank spectral method, leading to efficient data compression especially in velocity domain. In \cite{dektor2020dynamically, dektor2021dynamic}, dynamic tensor approximations for high dimensional linear and nonlinear PDEs are proposed based on functional tensor decomposition and dynamic tensor approximation. There are recent work of on low-rank methods with asymptotic preserving property for multi-scale models \cite{einkemmer2021asymptotic, einkemmer2021efficient, chen2020random}.

This paper proposed a dynamic low-rank approach with HT decomposition for the high dimensional nonlinear Vlasov model, when applying the finite difference WENO method coupled with second order SSP multistep method as the high order discretization. The flow map approach is also developed and its effectiveness is demonstrated via several linear and nonlinear examples. 
The organization of the paper is the following. Section 2 illustrate the main spirit of the low-rank approach, as well as the low-rank flow map approach, via a simple 1D1V Vlasov model.  Section 3 is on the extension to high dimensional problem using the HT decomposition of tensors. Section 4 presents extensive numerical results for linear hyperbolic equations and nonlinear kinetic models. Finally, the conclusion is given in Section 5.

%% file: algorithm.tex
\section{Low rank representation of Vlasov solution}

Inspired by existing understanding of the low-rank solution structure for Vlasov dynamics (e.g. Landau damping and two-stream instabilities) \cite{kormann2015semi, einkemmer2018low, ehrlacher2017dynamical}, as well as the observation that the differential operator in the Vlasov equation \eqref{eq: vlasov} can be represented in the tensorized form, we propose a novel way to dynamically and adaptively build up low-rank solution basis, and determine the low rank solutions in a tensor format. We will first demonstrate our proposed idea in a simplified 1D1V setting ($d$=2) using a high order spatial differentiation operator, coupled with a first order forward Euler time discretization. Here the high order spatial differential operator could come from the spectral collocation method \cite{hesthaven2007spectral} or flux-based finite difference approximation \cite{tadmor1998approximate}. We will discuss the extension of the algorithm to high order temporal discretization, followed by extension to general 3D3V ($d$=6) problems by using the hierarchical Tucker decomposition of tensors \cite{grasedyck2010hierarchical} or tensor train decomposition \cite{oseledets2011tensor}. Here and below, we denote $d$ as the dimension of the problem. 

\subsection{A low rank Vlasov solver in a simplified 1D1V setting}

We consider a simplified 1D1V VP system 
\beq
\frac{\partial f}{\partial t}  
+  {\bf{v}} \cdot \nabla_{\bf{x}}  f 
+ {\bf{E}} (t, {\bf{x}}) \cdot \nabla_{\bf{v}}  f = 0,
\label{vlasov1}
\eeq
\beq
 {\bf E}(t, {\bf x}) = - \nabla_{\bf x} \phi(t, {\bf x}),  \quad -\triangle_{\bf x} \phi (t, {\bf x}) = {{\bf \rho} (t, {\bf x})} - \rho_0,
\label{poisson}
\eeq
which describes the probability distribution function $f(t, {\bf x}, {\bf v})$ of electrons in collisionless plasma. 
Here ${\bf E}$ is the electric field and $\phi$ is the self-consistent electrostatic potential. $f$ couples to the long range fields via the charge density, ${\bf \rho}(t, {\bf x}) = \int_{\mathbb{R}} f(t, {\bf x}, {\bf v}) d {\bf v} -1$, where we take the limit of uniformly distributed infinitely massive ions in the background.  The VP system describes the movement of electrons due to self-induced electric field ${\bf E}$ determined by the Poisson equation. 

\noindent
{\bf The proposed low rank Vlasov solver.} 
Our proposed low rank Vlasov solver is built base on the assumption that our solution at time $t^n$ has a low-rank representation in the form of 
\begin{equation}
\label{eq: fn1}
%f(x, v, t^n) = \sum_{j=1}^{r^n} \left(C_j^n \ \ U_j^{(1), n}(x) \otimes U_j^{(2), n}(v)\right),
f(x, v, t^n) = \sum_{j=1}^{r^n} \left(C_j^n \ \ U_j^{(1), n}(x) \cdot U_j^{(2), n}(v)\right),
\end{equation}
where $\left\{U_j^{(1),n}(x)\right\}_{j=1}^{r^n}$ and $\left\{U_j^{(2),n}(v,n)\right\}_{j=1}^{r^n}$ are a set of low rank {unit length orthogonal} basis in $x$ and $v$ directions respectively, $C_j^n$ is the coefficient for the basis $U_j^{(1),n} \cdot U_j^{(2),n}$, and $r^n$ is the rank of the  tensor representation. 
For the Vlasov dynamics, we propose to adaptively update our low-rank basis, and hence the coefficient $C$, the basis $U^{(1)}$, $U^{(2)}$ and the rank $r$ are all time dependent with superscript $n$. 

\begin{figure}[h!]
	\centering
	\includegraphics[height=40mm]{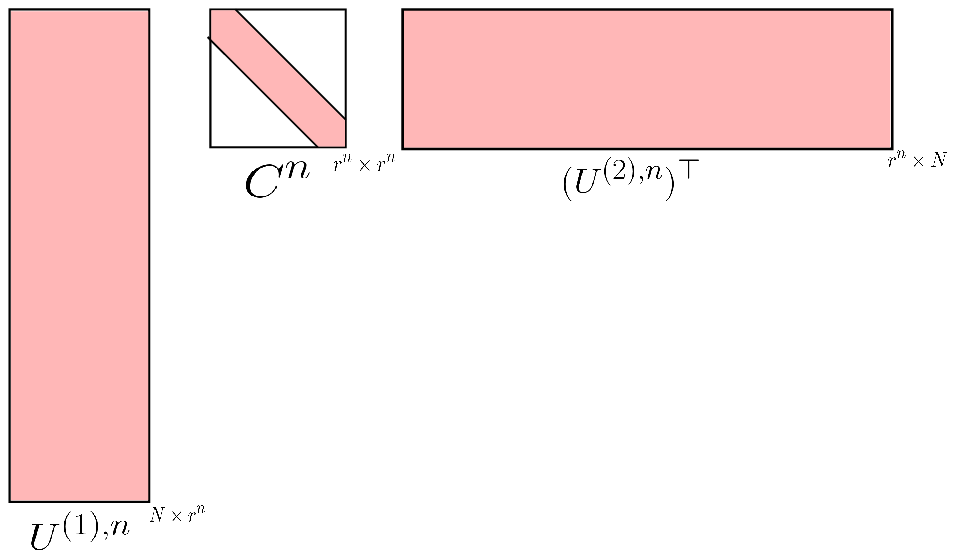}
	\label{fig:1}
\end{figure}

We choose to work with solutions on uniformly distributed $N$ grid points in each dimension for $U_j^{(1), n}(x)$ and $U_j^{(2), n}(v)$ respectively; thus $f^n$ in equation \eqref{eq: fn1} with grid point discretization can be written in the following tensor product form and in a matrix form
\begin{equation}
\label{eq: fn2}
f^n =\sum_{j=1}^{r^n} \left(C_j^n \ \ {\bf U}_j^{(1), n} \otimes {\bf U}_j^{(2), n}\right) = {\bf U}^{(1),n}\cdot C^n\cdot({\bf U}^{(2),n})^\top.
% \mbox{plot a rectangular matrix of size N by r, r by r, r by N}
\end{equation}
Here columns of ${\bf U}^{(1),n}$, i.e., ${\bf U}_j^{(1), n}$, $j=1 \cdots {r^n}$, are point values of $U_j^{(1), n}(x)$, $j=1 \cdots {r^n}$ at uniformly distributed $N$ grid points in $x$-direction; similarly rows of ${\bf U}^{(2),n}$ are point values of $U_j^{(2), n}(v)$ at uniform grids in $v$-direction. $C^n$ is a diagonal matrix of size $r^n$ representing the coefficient for the tensor product basis. See Figure~\ref{fig:1} for illustration. 
%In the following, depending on the context, we will use the low rank representation in a functional form \eqref{eq: fn1} and its discretized version in a matrix form \eqref{eq: fn2} interchangeably. 

We consider a simple first order forward Euler discretization of \eqref{vlasov1}, to illustrate the main idea in dynamically and adaptively updating basis and solutions. 
\begin{enumerate}
\item {\em Add basis.} 
Consider forward Euler discretization of the Vlasov equation \eqref{vlasov1} 
%with $D_x$ and $D_v$ being high order locally mass conservative discretization of differentiation terms, we have 
\begin{equation}
\label{eq: fn3}
f^{n+1} = f^n - \Delta t (v D_x (f^n) + E^n(x) D_v (f^n)).
%f^{n+1} = f^n - \Delta t ( D_x\otimes vI (f^n) + E^n(x)\otimes D_v (f^n)).
\end{equation}
Here $D_x$ and $D_v$ represent a high order locally mass conservative discretization of spatial differentiation terms. It could be the spectral collocation method \cite{hesthaven2007spectral} or the flux-based finite difference type numerical differentiation \cite{tadmor1998approximate}. We also assume that $E^n$ can be computed from \eqref{poisson} accurately and efficiently, see discussions on Poisson solver below. Thanks to {the tensor friendly form} of the Vlasov equation, $f^{n+1}$ can be evolved from $f^n$ \eqref{eq: fn2}, and be represented in the following low-rank format:
%With the tensor representation , we have
%\begin{eqnarray*}
\begin{equation}
\label{eq:lowrankmethod}
f^{n+1} = \sum_{j=1}^{r^n} C_j^n \left[  \left( {\bf U}_j^{(1), n} \otimes {\bf U}_j^{(2), n}\right)  - \Delta t 
\left( D_x {\bf U}_j^{(1), n} \otimes v \star {\bf U}_j^{(2), n} + E^n \star {\bf U}_j^{(1), n} \otimes D_v {\bf U}_j^{(2), n}
\right)\right],
%\end{eqnarray*}
\end{equation}
where $\star$ demotes an element-wise multiplication operation. 
Here we see that the number of basis has increased from $r^n$ (for $f^n$) to $3r^n$ (for $f^{n+1}$) in a single step update. In particular, 
a basis in $f^n$, e.g. ${\bf U}_j^{(1), n} \otimes {\bf U}_j^{(2), n}$ has evolved into three basis with 
\begin{equation}
\left\{ {\bf U}_j^{(1), n} \otimes {\bf U}_j^{(2), n}, \ \ D_x {\bf U}_j^{(1), n} \otimes v \star {\bf U}_j^{(2), n},\, \ \ E^n\star {\bf U}_j^{(1), n} \otimes D_v {\bf U}_j^{(2), n} \right\}.
\end{equation}
This step is illustrated in Figure~\ref{fig:remove} (a). 
The computational cost of the `adding basis' step is $\mathcal{O}(r N \log(N))$ if a global spectral differentiation is performed and is $\mathcal{O}(rN)$ if a local finite difference type numerical differentiation is performed.
%The computational cost of the `adding basis' step is $\mathcal{O}(r d N \log(N))$ if a global spectral differentiation is performed and is $\mathcal{O}(rdN)$ if a local finite difference type numerical differentiation is performed.
%Note that the tensor formulation \eqref{eq:lowrankmethod} allows a computational cost of $\mathcal{O}(r d N \log(N))$, assuming the matrix-vector multiplication $D_x U_j^{(1), n}$ and $D_v U_j^{(2), n}$ are computed via FFT. This is in contrast to $O(N^d \log(N))$ for the traditional full grid spectral collocation method.

\item {\em Remove basis.} If no basis is removed, then the rank of the tensor approximation would grow exponentially as time evolves. Hence, the removing basis procedure is crucial for the efficiency of the low-rank method. A SVD-type truncation procedure is proposed as following. 
We start with the pre-compressed solution $f^{n+1}$ from \eqref{eq:lowrankmethod}, see Fig.~\ref{fig:remove}(a) in which the red and blue parts refer to the old and newly added basis, respectively. These new set of basis is not necessarily orthogonal; so we perform the Gram-Schmidt process, e.g., QR decomposition to orthogonalize the basis, see the cyan matrices in Fig.~\ref{fig:remove}(b).  Then we apply a truncated SVD to the product of three $3r^n\times3r^n$ matrices based on a  prescribed threshold, see Fig.~\ref{fig:remove}(b-c). In this step, the rank of $f^{n+1}$ is being reduced from $3r^n$ to $r^{n+1}$. By combining the orthogonal matrices from QR decomposition and SVD, see the cyan and red matrices in Fig.~\ref{fig:remove}(c), we obtain the compressed solution $f^{n+1}$ with the updated basis $U^{(1),n+1}$ and $U^{(2),n+1}$, see Fig.~\ref{fig:remove}(d). The computational cost of the `removing basis' step is $\mathcal{O}(r^2 N + r^3)$.
%The computational cost of the `removing basis' step is $\mathcal{O}(r^2 d N + r^3)$.

\begin{figure}[h!]
	\centering
	\subfigure[]{\includegraphics[height=30mm]{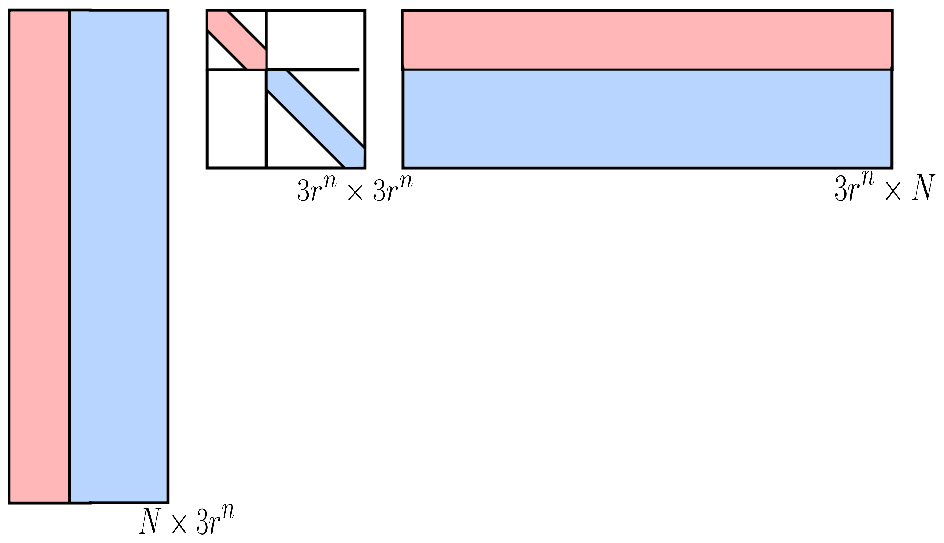}}\quad\quad\quad
	\subfigure[]{\includegraphics[height=30mm]{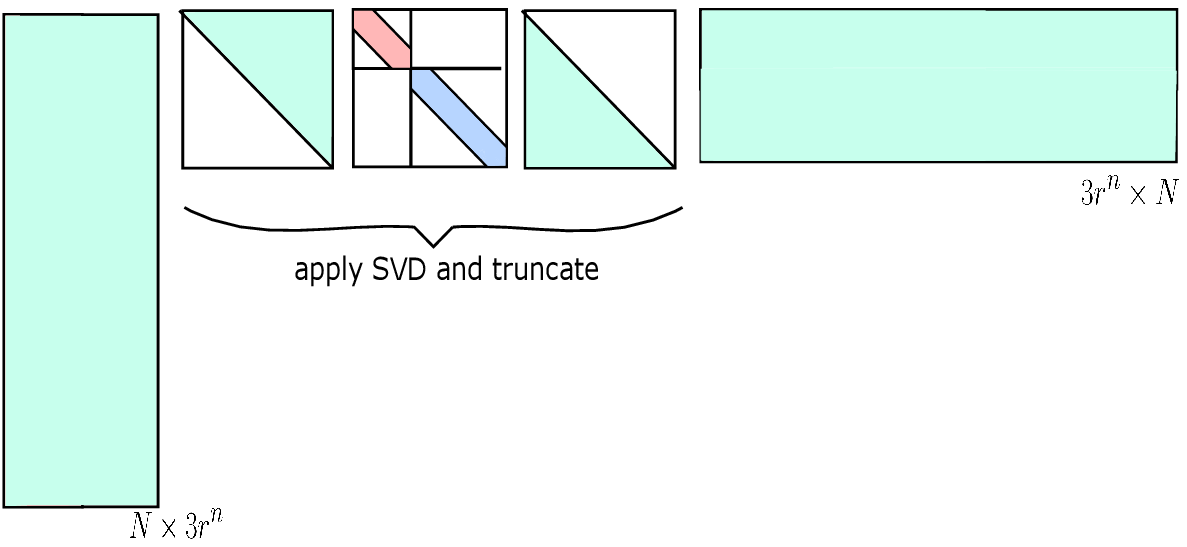}}\\
	\subfigure[]{\includegraphics[height=34mm]{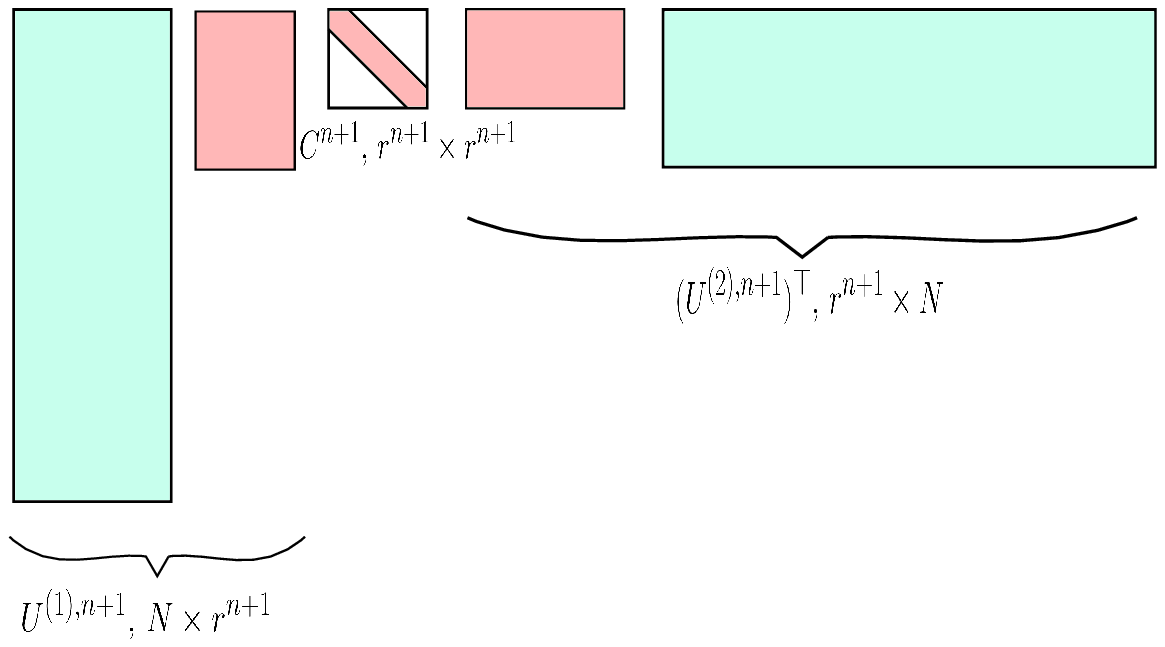}}\quad\quad\quad
	\subfigure[]{\includegraphics[height=34mm]{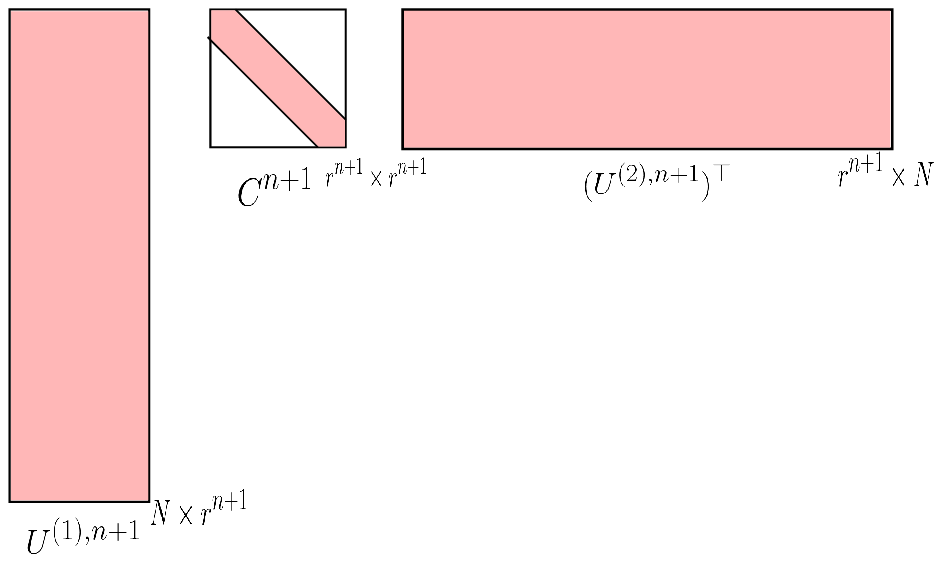}}
	\caption{Remove basis.\label{fig:remove}}
\end{figure}
\end{enumerate}

\noindent
{\bf High order spatial discretization of $D_x$ and $D_v$.} $D_x$ and $D_v$ can be viewed as the differentiation matrices for the corresponding variables. We perform fifth order finite difference method, derived from computing fluxes based on the upwind principle and taking flux differences to ensure local mass conservation \cite{johnson1997advanced}. In particular, with the upwind principle, one has to use different differentiation operator $D_x$ for positive and negative $v$'s. We let \[
v^+ =\max(v, 0), \quad v^- = \min(v, 0), \qquad E^{n,+} = \max(E^n, 0), \quad E^{n, -} = \min(E^n, 0), 
\] 
and let $D^{\pm}_x$ be the upwind differentiation operator in $x$-direction corresponding to $v^{\pm}$ respectively. Similarly, we let $D^{\pm}_v$ be the upwind differentiation operator in $v$-direction  corresponding to $E^{n, \pm} $ respectively. 
 Thus, similar to \eqref{eq:lowrankmethod} in the add basis step, we have, with the upwind differentiation operator,    
\be
\label{eq:lowrankmethod_upwind1}
f^{n+1} &&= \sum_{j=1}^{r^n} C_j^n \left[  \left( {\bf U}_j^{(1), n} \otimes {\bf U}_j^{(2), n}\right) \right. \\ 
&&- \Delta t 
\left( D^+_x {\bf U}_j^{(1), n} \otimes v^+ \star {\bf U}_j^{(2), n} + D^-_x {\bf U}_j^{(1), n} \otimes v^- \star {\bf U}_j^{(2), n} \right.\nonumber\\ 
&&\left.\left.+ E^{n, +} \star {\bf U}_j^{(1), n} \otimes D^+_v {\bf U}_j^{(2), n}+E^{n, -} \star {\bf U}_j^{(1), n} \otimes D^-_v {\bf U}_j^{(2), n}
\right)\right],\nonumber
\ee
Here we see that the number of basis has increased from $r^n$ (for $f^n$) to $5r^n$ in the add basis step. 
%\QQ{Move this to high D discussion. For a dimensional problem, the increase of rank is from $r^n$ (for $f^n$) to $ (1+ 2 m)r^n$ (for $f^{n+1}$), where $m$ is the number of RHS terms and is usually associated with the dimension of the problem $d$. } 
WENO type reconstructions can be applied for those flux functions to avoid numerical oscillations due to under-resolution of Vlasov solutions. We note that other type of spatial discretizations are possible, such as the spectral method, with a global differentiation matrix \cite{hesthaven2007spectral} can be applied. In this case, there is no upwind biased differentiation operator and the growth of basis is three folded, i.e. from $r^n$ (for $f^n$) to $3r^n$. We also remark that, despite the growth of basis, in the removing basis step, the solution rank will remain low, if the solution displays a low rank structure.   

%WENO type reconstructions can be applied at flux functions 

% For the spectral collocation method,   spectral differentiation matrices for the corresponding variables; while for the finite difference WENO method, extra consideration of upwind mechanism has to be taken into account. 

\noindent
{\bf High order temporal discretization.}
In the above procedure of `adding and removing basis' for the Vlasov equation, the set of basis as well as their coefficients are being updated in each time step. Such an idea of updating the solution can be extended to high order accuracy in time by a strong stability preserving (SSP) multi-step method or a SSP multi-stage RK method \cite{gottlieb2011strong}. 
For example, for a second order SSP multi-step method, 
\begin{eqnarray*}
f^{n+1} &=& \frac14 f^{n-2} + \frac34 f^n - \frac32\Delta t (v D_x (f^n) + E^n(x) D_v (f^n)),
\end{eqnarray*} 
the rank will increase from $r^{n}$ to $r^{n-1} + 5 r^n$ if an upwind differentiation such as \eqref{eq:lowrankmethod_upwind1} is used. For a third order SSP multi-step method 
\begin{eqnarray*}
f^{n+1} &=& \frac{11}{27} f^{n-3} + \frac{16}{27} f^{n} - \frac{16}{9} \Delta t (v D_x (f^n) + E^n(x) D_v (f^n))+ \frac{4}{9} \Delta t (v D_x (f^{n-1}) + E^{n-1}(x) D_v (f^{n-1}))
\end{eqnarray*} 
the rank will increase to $r^{n-3} +4r^{n-1}+5r^n$ in the ``add basis" step if an upwind differentiation described above is used. The SSP multi-step methods have advantages, compared with the multi-stage RK methods, in that the rank expand per time step is much smaller, if no rank-truncation is performed at intermediate RK stages. Aggressive rank-truncation (with relatively large threshold for truncating singular values) at intermediate RK stages may lead to temporal order reduction, as the temporal accuracy of RK methods rely on delicate combination of RK intermediate solutions; while the mild rank-truncation with small threshold will lead to a greater computational expense due to the faster growth of rank in the ``add basis" step. 
%will not reduce the rank of intermediate RK solutions much. 
% \QQ{Old text:
%As a reference for comparison, for the SSP RK2 method, the rank will expand 9 fold; while for the SSP RK3, the rank will expand $31$ times. Although the removing basis step will be able remove the redundant basis, the third order SSP multi-stage RK method is still expensive to implement, especially when consider the high dimensional setting.} 

\noindent
{\bf Poisson solver.} To solve the Poisson equation, we first compute the charge density $\rho(x, t) = \int f(x, v, t)dv - 1$. In the low rank format, we have from \eqref{eq: fn1}
\beq
\rho(x, t^n) = \int f(x, v, t^n)dv = \sum_{j=1}^{r^n} \left(C_j^n \ \ U_j^{(1), n}(x) \cdot \int U_j^{(2), n}(v)dv\right),
\eeq
with its discretized version
\beq
\label{eq: poisson_lr}
\rho^n = \sum_{j=1}^{r^n} C_j^n \ \ {\bf U}_j^{(1), n}  \ \Delta v \left ({\bf U}_j^{(2), n} \cdot \mathbf{1} \right).
\eeq
Here $ \Delta v \left ({\bf U}_j^{(2), n} \cdot \mathbf{1} \right)$ is the application of mid point rule for velocity integration, which is spectrally accurate for smooth solution and the zero boundary condition. A fast Fourier transform, or a high order finite difference Poisson solver can be applied to the Poisson solver. 
%\QQ{could move this part to high D.}

\noindent
{\bf Algorithm flow chart.} We organize the flow chart as Algorithm \ref{alg: low_rank} below for the low rank approach with fifth order finite difference for spatial discretization and second order SSP multi-step method for time discretization for the 1D1V VP system. 

\bigskip
\begin{algorithm}[H]
\label{alg: low_rank}
  \caption{Low rank approach with fifth order finite difference for spatial discretization and second order SSP multi-step method for time discretization for the 1D1V VP system.}
  \begin{enumerate}
\item Initialization:
  \begin{enumerate}
  \item 
   Initial distribution function $f(x, v, t=0)$ in a low rank format \eqref{eq: fn2}. 
 \end{enumerate}
\item For each time step evolution from $t^n$ to $t^{n+1}$. 
  \begin{enumerate}
  \item Compute $\rho(x, t^n) = \int f(x, v, t^n)dv -1 $ in the low rank format \eqref{eq: poisson_lr}, followed by computing $E(x, t^n)$ from the Poisson equation \eqref{poisson} by fast Fourier transform or a high order finite difference algorithm.
  \item Add basis: 
  \be
f^{n+1} &&=  \frac14 \sum_{j=1}^{r^{n-2}} \left(C_j^{n-1} \ \ {\bf U}_j^{(1), n-2} \otimes {\bf U}_j^{(2), n-2}\right)+ \frac34 \sum_{j=1}^{r^n} \left(C_j^n \ \ {\bf U}_j^{(1), n} \otimes {\bf U}_j^{(2), n}\right) \nonumber
\\\nonumber
&&- \frac32\Delta t \sum_{j=1}^{r^n} C_j^n \left[ 
\left( D^+_x {\bf U}_j^{(1), n} \otimes v^+ \star {\bf U}_j^{(2), n} + D^-_x {\bf U}_j^{(1), n} \otimes v^- \star {\bf U}_j^{(2), n} \right.\right.\\ \nonumber
&&\left.\left.+ E^{n, +} \star {\bf U}_j^{(1), n} \otimes D^+_v {\bf U}_j^{(2), n}+E^{n, -} \star {\bf U}_j^{(1), n} \otimes D^-_v {\bf U}_j^{(2), n}
\right)\right]. 
\ee
\item Remove basis as illustrated in Figure~\ref{fig:remove} and update the rank $r^{n+1}$, as well as ${\bf U}_j^{(1), n+1}$, ${\bf U}_j^{(2), n+1}$ and coefficients $C_j^{n+1}$, for $j=1\cdots r^{n+1}$.      
\end{enumerate}
\end{enumerate}
  \end{algorithm}
%For example, for a second order SSP RK method, 
%\begin{eqnarray*}
%f^{(1)} &=&  f^n - \Delta t (v D_x (f^n) + E^n(x) D_v (f^n))\\
%f^{n+1} &=& \frac12 f^{(1)} - \frac{\Delta t}{2} (v D_x (f^{(1)}) + E^n(x) D_v (f^{(1)})).
%\end{eqnarray*} 
%\begin{eqnarray*}
%f^{(1)} &=&  f^n - \Delta t (v D_x (f^n) + E^n(x) D_v (f^n))\\
%f^{(2)} &=& \frac34 f^n+ \frac14 f^{(1)} - \frac{\Delta t}{4} (v D_x (f^{(1)}) + E^n(x) D_v (f^{(1)}))\\
%f^{n+1} &=& \frac13 f^n+ \frac23 f^{(2)} - \frac{2\Delta t}{3} (v D_x (f^{(2)}) + E^n(x) D_v (f^{(2)})).
%\end{eqnarray*} 
%In a time step, for the SSP RK2 the rank would increase from $r^n$ to $9 r^n$; and for SSP RK3, the rank will increase from $r^n$ to $31r^n$. Perform SVD process will reduce the rank back.  

\bigskip
\noindent
{\bf Mass conservation and numerical stability of the algorithm.} The proposed low rank algorithm is locally mass conservative, due to the flux difference form of the differentiation operator that we employ. The global mass conservation error is up to the truncation threshold with accumulation in time in the remove basis step. The stability of the rank truncation algorithm is discussed in \cite{rodgers2019stability}. It is the consequence of the stability of the original time stepping algorithm and the stability of the rank-truncation algorithm in the $L^2$ norm. 

\subsection{A low rank representation of Vlasov flow map in a simplified 1D1V setting}

Let $(\mathcal{X}^*(x, v, t), \mathcal{V}^*(x, v, t))$ be the $x$-$v$ coordinate of feet of characteristics at $t=0$ of the VP system originated from $(x, v, t)$. $(\mathcal{X}^*, \mathcal{V}^*)$ satisfies the same characteristics evolution equation as the nonlinear VP system \eqref{vlasov1}-\eqref{poisson}, 
%\beq
%\frac{d\mathcal{X}}{dt} = \mathcal{V}, \quad \frac{d\mathcal{V}}{dt} = E,
%\eeq
with the initial conditions $\mathcal{X}^*(x, v, t=0) = x$ and $\mathcal{V}^*(x, v, t=0) = v$. 
Thus, $(\mathcal{X}^*, \mathcal{V}^*)$ satisfy the PDEs
\begin{eqnarray}
\label{eq: X_eq}
&\mathcal{X}^*_t + v\cdot \mathcal{X}^*_x + E \cdot\mathcal{X}^*_v =0, 
 \quad
 \mathcal{X}^*(x, v, t=0) = x,\\
% \mathcal{X}(t=0) = {\mathbf x}  \otimes {\mathbf 1},\\
\label{eq: V_eq}
&\mathcal{V}^*_t + v \cdot\mathcal{V}^*_x + E\cdot \mathcal{V}^*_v =0, \quad
 \mathcal{V}^*(x, v, t=0) = v.
% \mathcal{V}(t=0) = {\mathbf 1}  \otimes {\mathbf v}.
\end{eqnarray}
Here $E$ is the electrostatic field from the Poisson equation \eqref{poisson} in the VP system. 

For the computational discretization, we work with the same set of uniformly distributed computational grid. With such a mesh, the initial conditions in \eqref{eq: X_eq} and \eqref{eq: V_eq} are rank one tensors
\[
\mathcal{X}^*(t=0) = {\mathbf x}  \otimes {\mathbf 1}, \quad \mathcal{V}^*(t=0) = {\mathbf 1}  \otimes {\mathbf v},
\]
where $\mathbf x$ and $\mathbf v$ are coordinates of grid points for the corresponding direction. As the equations \eqref{eq: X_eq} and \eqref{eq: V_eq} enjoy the same tensor friendly structure as the original Vlasov equation \eqref{vlasov1}, the same proposed low rank approach can be applied. In particular, $\mathcal{X}^*$ and $\mathcal{V}^*$ at $t^n$ can be approximated in the form of 
\be
\label{eq: lr_X}
\mathcal{X}^{*, n} &=& \sum_{j=1}^{r_\mathcal{X}^n} \left(({C_{\mathcal{X}})_j^n} \ \ U_{\mathcal{X}, j}^{(1), n}(x) \otimes U_{\mathcal{X}, j}^{(2), n}(v)\right), \\
\label{eq: lr_V}
\mathcal{V}^{*, n} &=& \sum_{j=1}^{r_\mathcal{V}^n} \left(({C_{\mathcal{V}})_j^n} \ \ U_{\mathcal{V}, j}^{(1), n}(x) \otimes U_{\mathcal{V}, j}^{(2), n}(v)\right). 
\ee
From the fact that solution stays constant along characteristics $\frac{d}{dt}f(x(t), v(t), t) = 0$, we have
\beq
f(x, v, t) = f(\mathcal{X}^*, \mathcal{V}^*, t=0). 
\eeq

\begin{algorithm}[H]
  \caption{Low rank flow map approach for the 1D1V Vlasov system.}
  \label{alg: flow_map}
  \begin{enumerate}
\item Initialization:
  \begin{enumerate}
  \item 
   Initial distribution function $f(x, v, t=0)$.
   % that allows a pointwise evaluation. 
   %, either explicitly given in an analytical form, or stored in a low rank form.
  \item
 Rank one initial conditions for $\mathcal{X}^*$ and $\mathcal{V}^*$: ${\mathbf x}  \otimes {\mathbf 1}$, ${\mathbf 1}  \otimes {\mathbf v}.$
 \end{enumerate}
\item Evolution of $f$,  $\mathcal{X}^{*}$ and $\mathcal{V}^{*}$ from $t^n$ to $t^{n+1}$. 
  \begin{enumerate}
  \item Compute $\rho(x, t^n) = \int f(x, v, t^n)dv -1 $, and then $E(x, t^n)$ from the Poisson equation \eqref{poisson} by fast Fourier transform.
  \item
Use the low rank Algorithm \ref{alg: low_rank} for \eqref{eq: X_eq}-\eqref{eq: V_eq} to update solutions from $(\mathcal{X}^{*, n}, \mathcal{V}^{*, n})$ to $(\mathcal{X}^{*, n+1}, \mathcal{V}^{*, n+1})$ in a low-rank form \eqref{eq: lr_X}-\eqref{eq: lr_V}. 
\item 
\beq
f(x, v, t^{n+1}) = f(\mathcal{X}^{*, n+1}, \mathcal{V}^{*, n+1}, t=0).
\label{eq: flow_map_ini}
\eeq 
\end{enumerate}
\end{enumerate}
  \end{algorithm}

The flow map approach is advantageous for problems whose flow maps display low-rank structure, while their solutions are not necessarily of low rank. We demonstrate such advantages in several numerical examples in the following section. On the other hand, we note that the flow map approach has several computational issues to be addressed in a practical setting. Firstly, the evaluation of the initial condition at feet of characteristics \eqref{eq: flow_map_ini} requires analytic form of initial condition, or a low-rank form where an operator in the spirit of semi-Lagrangian interpolation is needed; secondly, the evaluation of the charge density in a low rank format with computational efficiency in a high dimensional setting needs to be developed; finally when source terms are involved, such direct evolution of the flow map needs to be adjusted. These issues impose limitation on the applicability of the flow map approach to high dimensional nonlinear problems. Addressing these computational challenges will be subject to our future research. 

%\subsection{The low rank algorithm to tensor friendly linear transport equations}

%forward characteristics flow map will lead to clustering of points as time evolves. 
%
%we adopt backward characteristics 
%
%Notations with time discretization
%\bit
%\item initial condition is $f_0(x, v)$
%\item
%${\bf fm}(x, v, t^{n}_m) = (fmx, fmy)(x, v, t^{n}_m)$ representing the $x$ and $y$ coordinates of the trace back points emanating from $(x, v, t^{n})$ back to $t^m$ assuming $m\le n$ since we consider backward characteristics. 
%\item 
%${\bf fm}(x, v, t^{n+1}_0)  =   {\bf fm}({\bf fm}(x, v, t^{n+1}_n), t^{n}_0)$
%\item 
%$f(x, v, t^n) = f_0({\bf fm}(x, v, t^{n}_0))
%$
%\eit
%
%Apply the spatial discretization $(x_i, v_j)$ of size $N \times N$ and tensor notations
%\bit
%\item
%$FMX^{n}_m$ is an $N \times N$ matrix, with  $(FMX^{n}_m)_{ij}  = fmx(x_i, v_j,  t^{n}_m)$. Similarly for $FMY^{n}_m$. 
%\item 
%To obtain $FMX^{n+1}_0$  
%\eit

\section{A low rank high dimensional Vlasov solver by Hierarchical Tucker Decomposition  of tensors}

The tensor networks have become an effective tool to obtain a low-rank approximation for high-dimensional problems. One such tensor format is the HT format, which enable us to extend the proposed methodology to high dimensions alleviating the curse of dimensionality. Below, as an example, we formulate a low-rank tensor algorithm for solving a 2D2V VP system \eqref{eq: vlasov}. 
%\begin{equation}\label{eq:vm3d}
%	f_t + v_1f_{x_1} + (E_1+v_2B_3)f_{v_1} + (E_2-v_1B_3)f_{v_2} = 0,
%\end{equation}
\begin{equation}\label{eq:vp4d}
	f_t + v_1f_{x_1} + v_2f_{x_2}+ E_1f_{v_1} + E_2f_{v_2} = 0,
\end{equation}
where the electric field $(E_1, E_2)$ is solved from the coupled Poisson's equation. 

\subsection{Hierarchical Tucker decompositions \cite{hackbusch2012tensor}}
 In this paper, we employ the HT format to compress the high order tensors, aiming to alleviate the curse of dimensionality. 
  The celebrated Tucker format (\cite{tucker1966some,de2000multilinear}) seeks to express an order $d$ tensor $\ba\in\mathbb{R}^{N_1\times\cdots\times N_d}$ as
 $$
 \ba = \sum_{j_1}^{r_1}\cdots\sum_{j_d}^{r_d}\mathbf{b}_{j_1,\ldots, j_d}\cdot\left(
 \bU_{j_1}^{(1)}\otimes\cdots\otimes\bU_{j_d}^{(d)}
 \right),
  $$ 
  where {$\bU^{(\mu)}:=\{\bU^{(\mu)}_{j_\mu}\}_{j_\mu=1}^{r_{\mu}}$ provide frames (or a basis if   $\{\bU^{(\mu)}_{j_\mu}\}_{j_\mu=1}^{r_{\mu}}$ are linearly independent) of linear space range$(\mathcal{M}^{(\mu)}(\ba))$, $\mu=1,\ldots,d$,  and $\mathbf{b}\in\mathbb{R}^{r_1\times\cdots r_d}$ is called the core tensor that glues all the $d$ frames. $\br=\{r_\mu\}_{\mu=1}^d$ is called the Tucker rank and we denote by $r=\max(\br)$. The low-rank tensor approximation of a tensor in the Tucker format can be computed through the HOSVD with quasi-optimal accuracy. Let $N=\max_\mu N_\mu$. The storage cost of the Tucker format scales $\mathcal{O}(r^d + drN)$ with  exponential dependence on dimension $d$. Hence, the Tucker format still suffers the curse of dimensionality and is only feasible in moderately high dimensions. The HT format is introduced to overcome the shortcoming. Denote the dimension index $D=\{1,2,\ldots,d\}$  and define a dimension tree $\mathcal{T}$ which is a binary tree containing a subset $\alpha\subset D$ at each node. Furthermore, $\mathcal{T}$ has $D$ as the root node and $\{1\},\  \{2\},\ldots,\ \{d\}$ as the leaf nodes. Each non-leaf node $\alpha$ has two children nodes denoted as  $\alpha_{l}$ and $\alpha_{r}$ with $\alpha = \alpha_{l}\bigcup\alpha_r$ and $\alpha_{l}\bigcap\alpha_r = \emptyset$. For example, the dimension tree given in Figure \ref{fig:dimtree} can be used to approximate $f(x_1,x_2,v_1,v_2)$ in \eqref{eq:vp4d} in the HT format. The efficiency of the HT format lies in the nestedness property \cite{hackbusch2009new}: for a non-leaf node $\alpha$ with two children nodes $\alpha_l,\,\alpha_r$, then 
 \begin{equation}
\text{range}(\mathcal{M}^{(\alpha)}(\ba))\subset  \text{range}(\mathcal{M}^{(\alpha_l)}(\ba)\otimes\mathcal{M}^{(\alpha_r)}(\ba)),
 \end{equation}
which implies that there exists a third order tensor $\bB^{(\alpha)}\in\mathbb{R}^{r_{\alpha_l}\times r_{\alpha_r}\times r_{\alpha}}$, known as the transfer tensor, such that
 \begin{equation} 
 \label{eq:htd_nested}
 \bU_{j_\alpha}^{(\alpha)} = \sum_{j_{\alpha_l}=1}^{r_{\alpha_l}}\sum_{j_{\alpha_r=1}}^{r_{\alpha_l}} \bB^{(\alpha)}_{j_{\alpha_l},j_{\alpha_r},j_{\alpha}}\bU_{j_{\alpha_l}}^{(\alpha_l)}\otimes \bU_{j_{\alpha_r}}^{(\alpha_r)},\quad j_\alpha =1,\ldots,r_\alpha.
 \end{equation}
 By recursively making use of \eqref{eq:htd_nested}, a tensor in the HT format stores a frame at each leaf node and a third order transfer tensor at each non-leaf node based on a dimension tree.  Denote $ \br = \{r_\alpha\}_{\alpha\in\mathcal{T}}$ as the hierarchical ranks. The storage of the HT format scales as $\mathcal{O}(dr^3+drN)$, where $r=\max \br$. If $r$ is reasonably low, then the HT format avoids the curse of dimensionality. In summary, the HT format is fully characterized by the three key components, including a dimension tree, frames at leaf nodes and transfer tensors at non-leaf nodes, see Figure \ref{fig:dimtree}. The HT format is ideal for simulating tensor friendly high-dimensional PDEs such as the Vlasov equation, as it allows efficient implementations of many operators by the low-rank method. These include addition of tensors, SVD-type truncation, element-wise multiplication (i.e., Hadamard product), and application of linear operators (e.g., spectral differentiation matrix) to tensors. 
 
 \begin{figure}
\centering
\subfigure[]{
 \begin{tikzpicture}[%sibling distance=10em,
     level/.style={sibling distance=40mm/#1},
  every node/.style = {shape=rectangle, rounded corners,
    draw, align=center,
    top color=white, bottom color=blue!20}
    ]
  \node {$\{1,\,2,\,3,\,4\}$}
    child { node {$\{1,\,2\}$} 
    	child{ node{$\{1\}$}}
	child{ node{$\{2\}$}}
    }
    child { node {$\{3,\,4\}$}
    	child{ node{$\{3\}$}}
	child{ node{$\{4\}$}} 
	};
\end{tikzpicture}}\quad
\subfigure[]{
\begin{tikzpicture}[%sibling distance=10em,
  every node/.style = {shape=rectangle, rounded corners,
    draw, %align=center,
    top color=white, bottom color=blue!20},
    level/.style={sibling distance=40mm/#1}
    ]
  \node {$\bB^{(1,2,3,4)}$}
      child { node {$\bB^{(1,2)}$} 
    	child{ node{$\bU^{(1)}$}}
	child{ node{$\bU^{(2)}$}}
    }
    child { node {$\bB^{(3,4)}$}
    	child{ node{$\bU^{(3)}$}}
	child{ node{$\bU^{(4)}$}} 
	};
\end{tikzpicture}}
\vspace{0.2in}
 \caption{Dimension tree $\mathcal{T}$ to express fourth-order tensors in the HT format. \label{fig:dimtree}}
\end{figure}
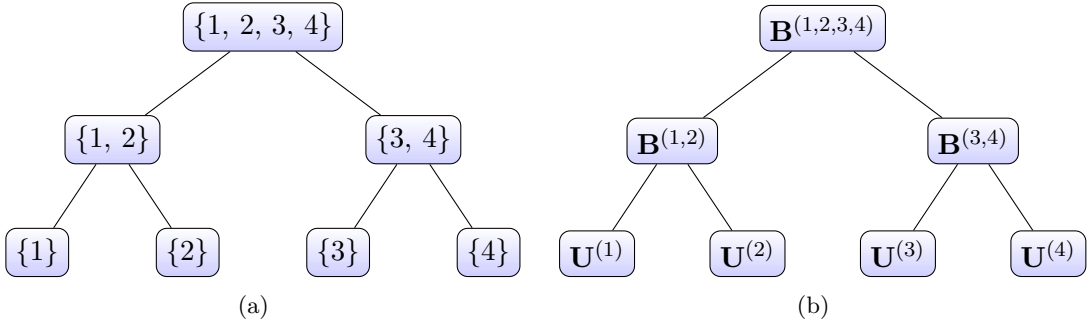

% \begin{figure}[h!]
%	\centering
%	\subfigure[]{\includegraphics[height=45mm]{dimtree3d.png}}
%	\caption{Dimension tree $\mathcal{T}$ to express $f(x_1,v_1,v_2,t)$ in the HTD format. \label{fig:dimtree}}
%\end{figure}

\subsection{A low-rank tensor method in HT for the 2D2V VP system}
Below, we formulate the low rank tensor method for solving the 2D2V VP system. We assume the solution $f^n$ at time step $t^n$ is expressed as the fourth-order tensor in the HT format with dimension tree $\mathcal{T}$ together with frames $\bU^{(1),n},\bU^{(2),n},\bU^{(3),n},\bU^{(4),n}$ at four leaf tensors, corresponding to directions $x_1,\,x_2,\,v_1,\,v_2$, respectively, and transfer tensors $\bB^{(1,2,3,4),n},\bB^{(1,2),n},\bB^{(3,4),n}$, see Figure \ref{fig:dimtree}. In particular, 
\beq
\label{eq:htd_f_nested}
f^n = \sum_{i_{12}=1}^{r_{12}}{\sum_{i_{34}=1}^{r_{34}}} \bB^{(1,2, 3, 4), n}_{i_{12},i_{34}, 1}\bU_{i_{12}}^{(1, 2), n}\otimes \bU_{i_{34}}^{(3, 4), n},
\eeq
with 
\beq
\label{eq: U12_htd}
\bU_{i_{12}}^{(1, 2), n} = \sum_{i_{1}=1}^{r_{1}}{\sum_{i_{2}=1}^{r_{2}}} \bB^{(1,2), n}_{i_{1},i_{2}, i_{12}}\bU_{i_{1}}^{(1), n}\otimes \bU_{i_{2}}^{(2), n},\quad i_{12} =1,\ldots,r_{12}, 
\eeq
and 
\beq
\label{eq: U34_htd}
\bU_{i_{34}}^{(3, 4), n} = \sum_{i_{3}=1}^{r_{3}}{\sum_{i_{4}=1}^{r_{4}}} \bB^{(3, 4), n}_{i_{3},i_{4}, i_{34}}\bU_{i_{3}}^{(3), n}\otimes \bU_{i_{4}}^{(4), n},\quad i_{34} =1,\ldots,r_{34}. 
\eeq
 Further, the electric field $E_1^n$ and $E_2^n$ are represented in the second order HT format. In the following, we discuss the details of low rank tensor method in the add basis, remove basis and Poisson solver steps of the algorithm. 
\begin{enumerate}
\item {\em Add basis.}
\begin{enumerate}
\item \underline{From the term ${\bf v} \cdot \nabla_{\bf x} f$.} 
We first focus on the Vlasov equation and consider the discretization of $v_1\cdot f_{x_1}^n$. The treatment of $v_2\cdot f_{x_2}^n$ is the same. 
Similar to the 1D1V case,  to account for upwinding, we split $v_1$ into $v_1^+ = \max(v_1,0)$ and $v_1^- = \min(v_1,0)$ and denote by $D_{x_1}^+$ and $D_{x_1}^-$ the biased upwind high order finite difference operators. Then $v_1\cdot f_{x_1}^n$ is approximated by the sum of two fourth-order tensors in the HT format denoted by $(v^+_1\otimes D_{x_1}^+) f^n$ and $(v^-_1\otimes D_{x_1}^-) f^n$. Both tensors have the same dimension tree, frames and transfer tensors as $f^n$, except that $(v^+_1\otimes D_{x_1}^+) f^n$ has frame $D^+_{x_1}\bU^{(1),n}$ in direction ${x_1}$ and frame $v^+_1\star\bU^{(3),n}$ in direction $v_1$, and $(v^-_1\otimes D_{x_1}^-) f^n$ has frame $D^-_{x_1}\bU^{(1),n}$ in direction ${x_1}$ and frame $v^-_1\star\bU^{(3),n}$ in direction $v_1$, see Figure~\ref{fig:dimtree_ex} (a).  
\item \underline{From the term ${\bf E} \cdot \nabla_{\bf v} f$.} 
We start from discussing the discretization of $E^n_1\cdot f^n_{v_1}$. Since $E_1^n$ is expressed in the HT format, we propose the following splitting strategy:
 \begin{equation}
 E^+_1 := \frac12\left(E_1^n + \alpha_1 \mathbf{1}_{x_1}\otimes \mathbf{1}_{x_2}\right),\quad E^-_1 := \frac12\left(E_1^n - \alpha_1 \mathbf{1}_{x_1}\otimes \mathbf{1}_{x_2}\right),
 \end{equation} 
where $\alpha_1= \max_{x_1,x_2} |E_1^n|$ is the maximum absolute value of $E_1^n$ over the physical domain.  The entries of $E^+_1$ and $E^-_1$ are nonnegative and nonpositive respectively.  Furthermore, $E^+_1$ and $E^-_1$  are still in the HT format. Then, $E^n_1\cdot f^n_{v_1}$ is approximated by the sum of two tensors $E_1^+ \otimes D_{v_1}^+\, f^n$ and $E_1^-\otimes D_{v_1}^-\, f^n$. % obtained by applying  $E_1^\pm$ and $D_{v_1}^\pm$ to $f^n$. 
 In particular, assume $E_1^+$ has {frames} $\bE^{+,(1)}_1$ and $\bE^{+,(2)}_1$ in direction $x_1$ and $x_2$ and transfer tensor $\bB^{(1,2)}$, 
\beq
\label{eq: E1_htd}
E^+_1 = \sum_{j_1=1}^{r_{1,E}} \sum_{j_2=1}^{r_{2,E}} \bB^{(1, 2)}_{E^+_1, j_1, j_2} \bE^{+,(1)}_{1, {j_1}} \otimes  \bE^{+,(2)}_{1, j_2}. 
\eeq
 Then, $E_1^+\otimes D_{v_1}^+\cdot f^n$ has a similar tree structure as $f^n$ with new frames $\bE^{+,(1)}_1\star \bU^{(1),n}$, $\bE^{+,(2)}_1\star \bU^{(2),n}$. From \eqref{eq: U12_htd} and \eqref{eq: E1_htd}, $\bE^{+,(1)}_1\star \bU^{(1),n}$ has $r_1 r_{1, E}$ frames with entry-wise multiplication
 \[
  \bE^{+,(1)}_{1, {j_1}} \star \bU_{i_{1}}^{(1), n},\quad j_1 = 1, \cdots r_{1, {E}},  \quad i_1 = 1, \cdots r_1. 
 \]
Correspondingly, the transfer tensor becomes $\bB^{(1,2)}_{E^+_1}\otimes^K\bB^{(1,2),n}$ with size  $r_1 r_{1, E} \times r_2 r_{2, E}$, where $\otimes^K$ denotes the generalized Kronecker product \cite{hackbusch2012tensor}.  
%For example, for two third order tensors $\ba\in \mathbb{R}^{N_1\times N_2\times N_3}$ and $\mathbf{b}\in \mathbb{R}^{M_1\times M_2\times M_3}$, we have $\ba\otimes^K\mathbf{b}\in \mathbb{R}^{N_1M_1\times N_2M_2\times N_3M_3}$ with entry  $(\ba\otimes^K\mathbf{b})_{ k_1,k_2,k_3}  = \ba_{i_1,i_2,i_3}\cdot \mathbf{b}_{j_1,j_2,j_3}$ where $k_\mu = (i_\mu-1)M_\mu + j_\mu$, $\mu=1,\,2,\,3$. 
 In addition, the frames in the $v_1$ direction are replaced with  $D^+_{v_1}\bU^{(1),n}$. See Figure~\ref{fig:dimtree_ex} (b) for illustration of the added tensor due to $E_1 \partial_{v_1} f$.  A similar extension can be done for the tensor from $E^n_2\cdot f^n_{v_2}$.
 \end{enumerate}
 When a second order SSP time integration method is used, then $f^{n+1}$ is approximated by $f^n$, together with about ten additional tensors in the adding basis step. We hence need to remove redundant basis to avoid exponential rank increase, meanwhile not compromising much accuracy.  
  \item {\em Remove basis.} 
 %As in the 1D1V case, we have to remove redundant basis to keep the hierarchical ranks of the numerical solution under control and meanwhile without compromising much accuracy. 
 We employ the standard hierarchical HOSVD for removing basis, which is implemented in the Matlab toolbox \texttt{htucker} \cite{kressner2014algorithm,tobler2012low}. The procedure is similar to the case of $d=2$, which consists of orthogonalizing  the frames and transfer tensors, computing the reduced Gramians and the associated eigen-decompostion with truncation (depending on the truncation threshold and the maximal rank allowed) at all nodes in the dimension tree, and then obtaining the truncated HT tensor $\tilde{f}^{n+1}$ with updated frames and transfer tensors. Such a truncation procedure costs $\mathcal{O}\left(d N r^{2}+(d-2) r^{4}\right)$, without exponential dependence on $d$, and ensures quasi-optimal truncation accuracy (\cite{grasedyck2010hierarchical})
  $$\|f^{n+1}-\tilde{f}^{n+1}\|_2\le\sqrt{2d-3}\|f^{n+1} -f_{best}\|_2,$$ 
  where $f_{best}$ is the best approximation to $f^{n+1}$ with the hierarchical ranks bounded by those  of $\tilde{f}^{n+1}$. We can further combine the addition and truncation procedures for improved efficiency and stability, which is proposed in  \cite{tobler2012low}.

\end{enumerate}

 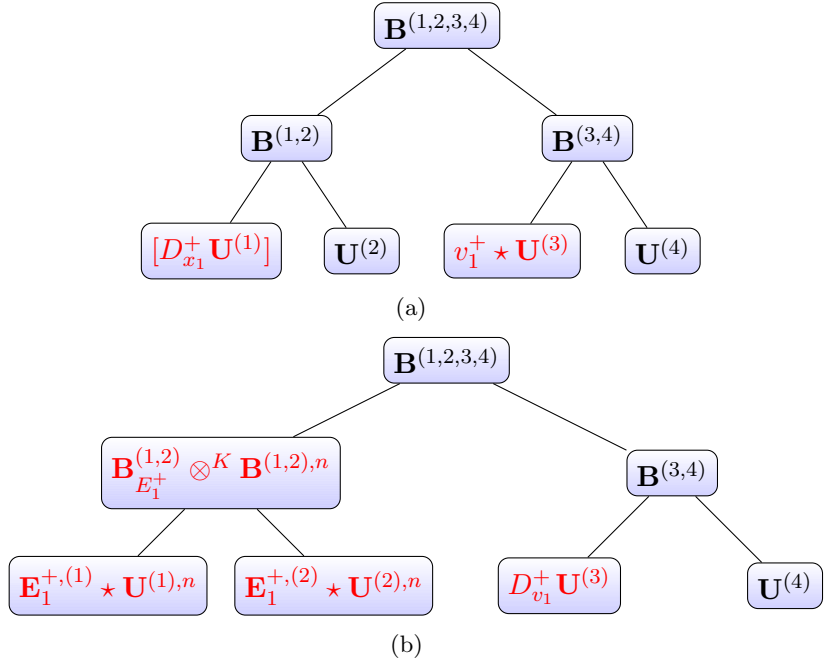
\begin{figure}[htp]
\centering
\subfigure[]{
\begin{tikzpicture}[%sibling distance=10em,
  every node/.style = {shape=rectangle, rounded corners,
    draw, %align=center,
    top color=white, bottom color=blue!20},
    level/.style={sibling distance=40mm/#1}
    ]
  \node {$\bB^{(1,2,3,4)}$}
      child { node {$\bB^{(1,2)}$} 
    	child{ node{\color{red}$[D_{x_1}^+\bU^{(1)}]$}}
	child{ node{$\bU^{(2)}$}}
    }
    child { node {$\bB^{(3,4)}$}
    	child{ node{\color{red}$v_1^+\star\bU^{(3)}$}}
	child{ node{$\bU^{(4)}$}} 
	};
\end{tikzpicture}}
\subfigure[]{
\begin{tikzpicture}[%sibling distance=10em,
  every node/.style = {shape=rectangle, rounded corners,
    draw, %align=center,
    top color=white, bottom color=blue!20},
    level/.style={sibling distance=60mm/#1}
    ]
  \node {$\bB^{(1,2,3,4)}$}
      child { node {\color{red}$\bB^{(1,2)}_{E^+_1}\otimes^K\bB^{(1,2),n}$} 
    	child{ node{\color{red}$\bE^{+,(1)}_1\star \bU^{(1),n}$}}
	child{ node{\color{red}$\bE^{+,(2)}_1\star \bU^{(2),n}$}}
    }
    child { node {$\bB^{(3,4)}$}
    	child{ node{\color{red}$D_{v_1}^+\bU^{(3)}$}}
	child{ node{$\bU^{(4)}$}} 
	};
\end{tikzpicture}}
\vspace{0.2in}
 \caption{(a) The added tensor from ${v_1} \partial_{x_1} f$; (b) the added tensor from the term ${E_1} \partial_{v_1} f$ in the HT format. \label{fig:dimtree_ex}}
\end{figure}

\subsection{Low rank method for solving Poisson's equation}
Now, we formulate a low-rank Poisson solver compatible with the Vlasov solver above. First, we need to compute the macroscopic charge density $\rho^n$ from $f^n$, that is
\begin{equation}
\label{eq:rho}
\rho^n_{i_1,i_2} =\Delta v_1\Delta v_2 \left(\sum_{i_3,i_4} f^n_{i_i,i_2,i_3,i_4}\right) - \rho_0, 
\end{equation}
with $\rho_0$ is the charge density of constant ion background. 
 Notice that the cost of direct summation scales as $\mathcal{O}(N^4)$, suffering the curse of dimensionality.  
As  $f^n$ is represented in the HT format, \eqref{eq:rho} can be computed via tensor contraction with cost $\mathcal{O}(rN+r^3)$. In particular, we first sum each frame vector in $\bU^{(3)}$ and $\bU^{(4)}$ and obtain $\mathbf{S}^{(3)}$ and $\mathbf{S}^{(4)}$, both of which are row vectors of size $r_3$ and $r_4$, respectively. See Figure \ref{fig: density_htd} (a). Then, we recursively merge the nodes enclosed by the dashed polygon in Figure~\ref{fig: density_htd} (a) and obtain a single node with transfer tensor $\bB^{(1,2)}_\rho$ in Figure~\ref{fig: density_htd} (b). 
%Hence, after the contraction, we obtain the second order tensor for approximating $\rho$.   
 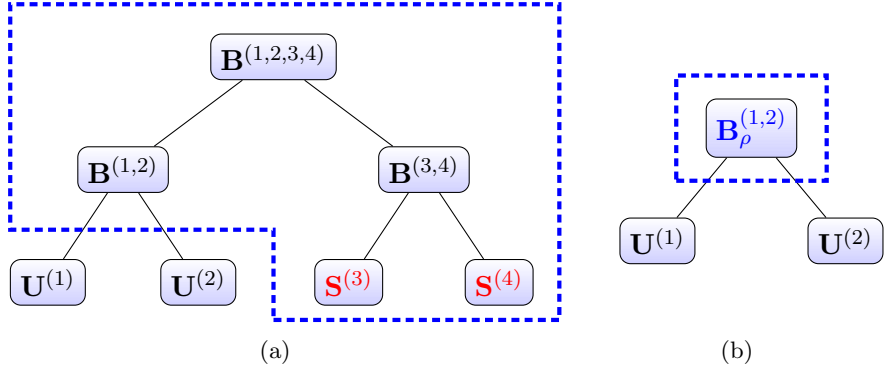
\begin{figure}
\centering
\subfigure[]{
\begin{tikzpicture}[%sibling distance=10em,
  every node/.style = {shape=rectangle, rounded corners,
    draw, %align=center,
    top color=white, bottom color=blue!20},
    level/.style={sibling distance=40mm/#1}
    ]
  \node {$\bB^{(1,2,3,4)}$}
      child { node {$\bB^{(1,2)}$} 
    	child{ node{$\bU^{(1)}$}}
	child{ node{$\bU^{(2)}$}}
    }
    child { node {$\bB^{(3,4)}$}
    	child{ node{\color{red}$\mathbf{S}^{(3)}$}}
	child{ node{\color{red}$\mathbf{S}^{(4)}$}}
	};
 %\draw[style=help lines] (-4,-4) grid (4,2);
  \draw[blue, ultra thick, densely dashed]     (0,-3.5) -- (3.8,-3.5) -- (3.8,0.7) -- (-3.5,0.7) -- (-3.5,-2.3) -- (0,-2.3) -- cycle;
\end{tikzpicture}}
\quad
\subfigure[]{
\raisebox{2em}{
\begin{tikzpicture}[%sibling distance=10em,
  every node/.style = {shape=rectangle, rounded corners,
    draw, %align=center,
    top color=white, bottom color=blue!20},
    level/.style={sibling distance=25mm/#1}
    ]
  \node {$\color{blue}\bB^{(1,2)}_\rho$}
    	child{ node{$\bU^{(1)}$}}
	child{ node{$\bU^{(2)}$}};
 %\draw[style=help lines] (-2,-2) grid (2,2);
  \draw[blue, ultra thick, densely dashed]     (-1,-0.7) -- (1,-0.7) -- (1,0.7) -- (-1,0.7) -- cycle;
\end{tikzpicture}}}
\caption{The contraction of the $f$ in the HT format to obtain $\rho$ via \eqref{eq:rho}.}
\label{fig: density_htd}
\end{figure}
Once the $\rho$ is computed in the HT format, we adopt a low-rank conjugate gradient (CG) method proposed in \cite{grasedyck2018distributed} for solving the Poisson equation \eqref{poisson}. Assuming the discretization of Laplacian with the spectral method, we have the linear system $\mathcal{A}\mathbf{\phi}=\mathbf{b}$, where $\mathcal{A}=D_{x_1}^2 \otimes I + I\otimes D_{x_2}^2$ denotes the Laplace operator with spectral differentiation matrices $D_{x_1}^2$, $D_{x_2}^2$ {in the tensor product form}. Note that the low-rank CG method requires truncation for each iteration to avoid exponential rank increase. 

%% file: numerical.tex
\section{Numerical results}
In this section we present a collection of numerical examples to demonstrate the efficiency and efficacy of the proposed low-rank tensor methods for simulating linear and nonlinear transport equations in high dimensions. In the simulations, fifth order upwind finite difference methods are employed for spatial discretization, together with a second order SSP multi-step method denoted by SSPML2 for temporal discretization. The numerical solutions of high dimensions are represented in the HT format \cite{hackbusch2012tensor}. Below, we denote the low-rank method for approximating the solution as approach I and the method for the flow maps as approach II. We compare the performance of approach I and II for linear and 1D1V VP system. For the 2D2V VP system, we only apply the approach I, due to the complication/high computational cost in evaluating the charge density for the approach II.   
%The methods are implemented in a serial manner, while parallel implementation is possible for improved efficiency for very large $d$, see \cite{grasedyck2018distributed}.

\subsection{Linear transport problems}
\begin{exa}\label{ex:adv} We first consider the advection equation with constant coefficients
	$$
	u_t + \sum_{m=1}^d u_{x_m}= 0,\quad \bx\in[-\pi,\pi]^d
	$$
with periodic conditions. We first consider a smooth initial condition with $d=4$
\begin{equation*}
u(\bx,t=0) = \exp\left(-2(x_1^2+x_2^2)\right)\sin(x_3+x_4), 
\end{equation*}
which can be expressed as a rank-two tensor in the CP format. Since the exact solution remains smooth and low-rank over time, we are able to test the accuracy and efficiency of the low-rank tensor method (approach I) with a properly chosen truncation threshold $\varepsilon$.  In the simulation, we set the same mesh size $N$ in each dimension and $\varepsilon = 10^{-6}$ and compute the solution up to $T=2\pi$. In Table \ref{tb:linear4d}, we report the $L^2$ error and associated order of convergence for the low-rank method. The fifth order linear FD method with upwinding is employed. Second order accuracy is observed due to SSPML2 used. In Figure \ref{fig:adv4d}, we report the time history of the hierarchical ranks of the solution in the HT format for $N=128$. It is observed that the hierarchical ranks of the solution stay very low and bounded by 3, making the method extremely efficient. Furthermore, CPU cost is 4.2s, 7.8s, 14.2, 30.1s, 65.3s for $N=16,\,32,\,64,\,128,\,256$, which grows linearly with the mesh size $N$. This contrasts strongly with the traditional full-rank counterpart which usually grow at the rate of $\mathcal{O}(2^{d+1})$. It is known that the sparse grid approach is incompetent in approximating Gaussian functions in high dimensions \cite{pfluger2010spatially}, while this example demonstrates the efficiency of the proposed low-rank tensor approach in this regard. Then, we consider a 2D discontinuous initial condition of a cross shape. The fifth order WENO method is employed for spatial discretization. The numerical solution is computed up to $t=2\pi$ and  plotted in Figure \ref{fig:adv2d}. Note that the solution remains a rank-two tensor over time despite being discontinuous, see the right plot in Figure \ref{fig:adv4d}. It is observed that the low-rank structure of the solution is well captured by the method. 
%In Figure \ref{fig:adv}, we plot the solution profiles of the problem with the discontinuous cross-shaped initial. WENO5 and WENO9 are used for spatial discretization, and we let $\varepsilon=10^{-5}$ for SVD truncation. The exact solution stays rank-two, while the ranks of the numerical solutions by WENO5 and WENO9 remain two and seven over time, respective. Due to the use of the WENO methodology, the low-rank solutions are observed to be essentially non-oscillatory.

\begin{table}[!hbp]
	\centering
	\caption{Example \ref{ex:adv}, $d=4$. $t=2\pi$. Convergence study.}
	\label{tb:linear4d}
	\begin{tabular}{|c|c|c|}
		\hline
		 $N$ & $L^2$ error & order \\\hline
16	&	2.56E-02	&		\\\hline
32	&	5.76E-03	&	2.15	\\\hline
64	&	1.41E-03	&	2.04	\\\hline
128	&	3.52E-04	&	2.00	\\\hline
256	&	8.09E-05	&	2.12	\\\hline

\hline

	\end{tabular}
\end{table}

\begin{figure}[h!]
	\centering
	\subfigure[]{\includegraphics[height=60mm]{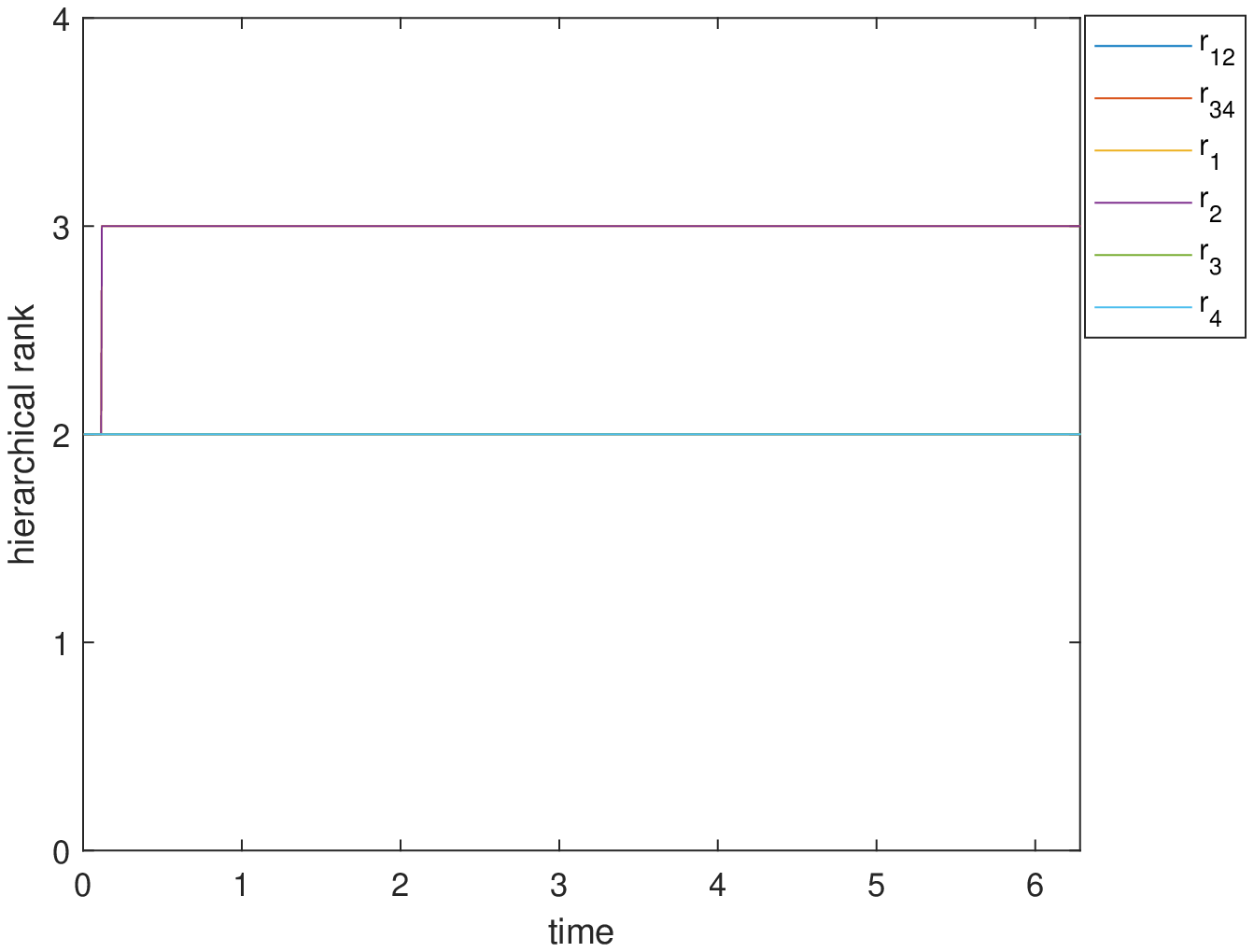}}
		\subfigure[]{\includegraphics[height=60mm]{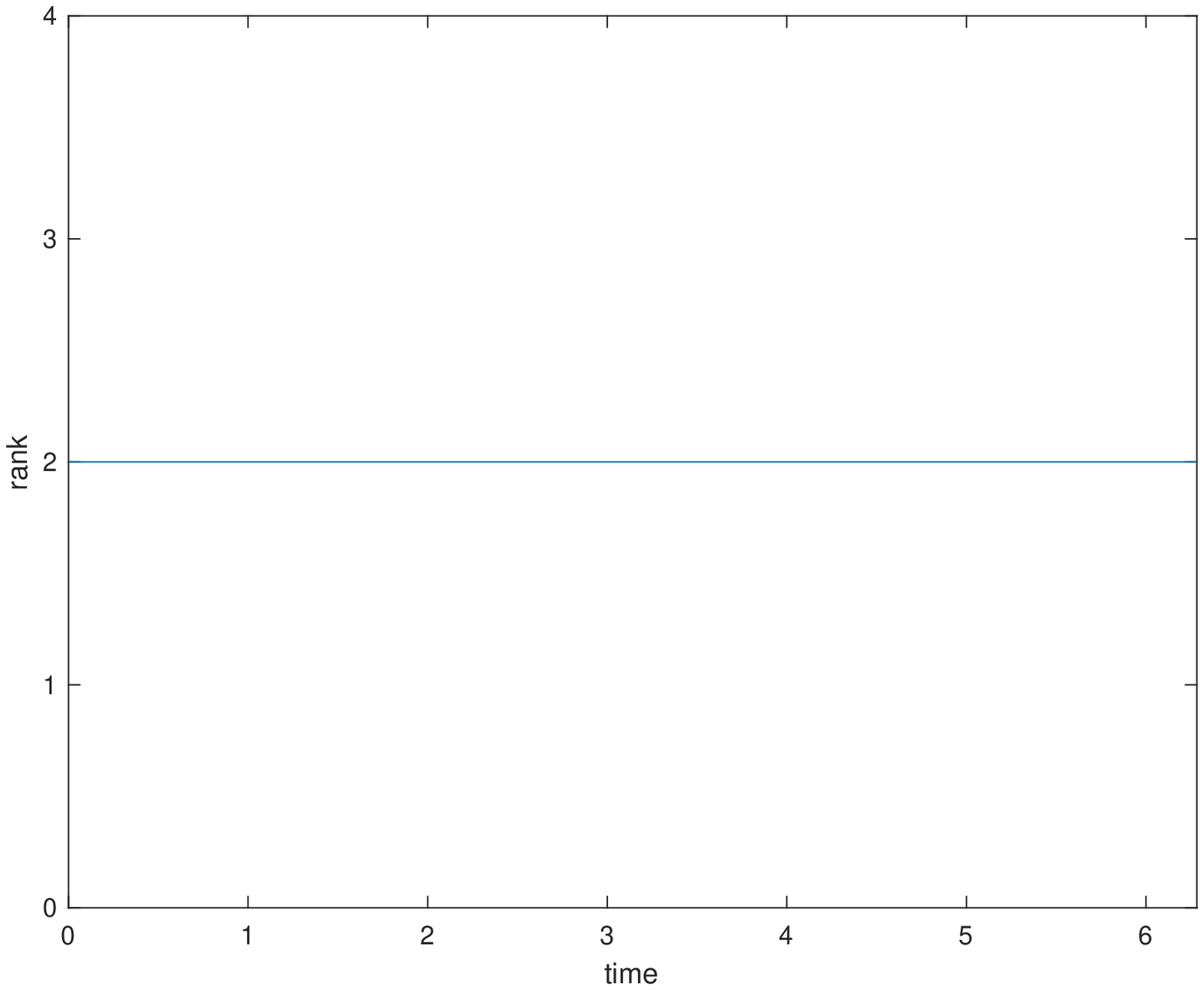}}
	\caption{Example \ref{ex:adv}. The time evolution of hierarchical ranks of the solution in the HT format. Left:	
	 $d=4$ for a smooth initial condition, $N=128$, $\varepsilon=10^{-6}$,  $t=2\pi$; right: $d=2$ for a cross-shaped initial condition, $N=128$, $\varepsilon=10^{-5}$, $t=2\pi$.}
	\label{fig:adv4d}
\end{figure}

\begin{figure}[h!]
	\centering
	\subfigure[]{\includegraphics[height=60mm]{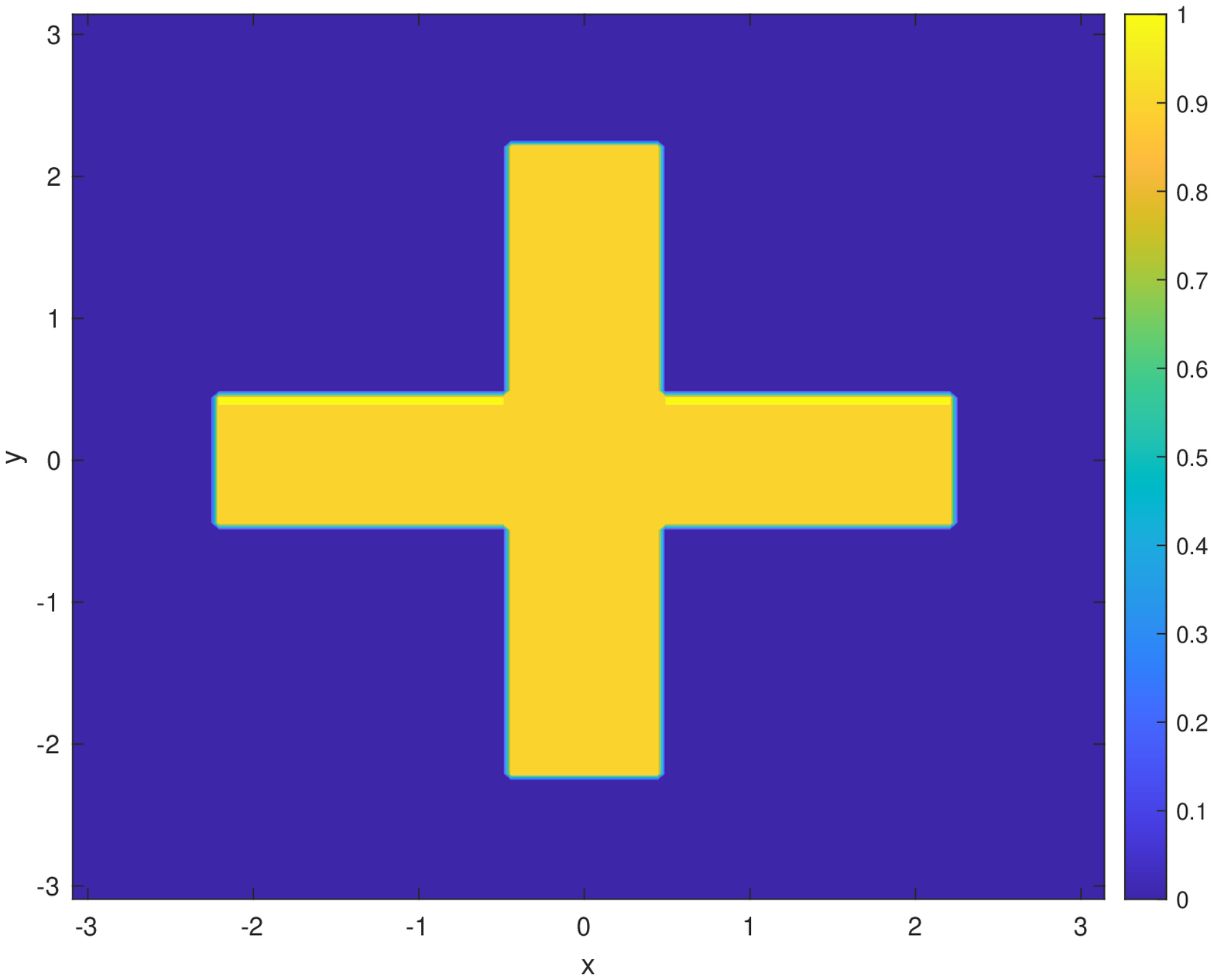}}
	\subfigure[]{\includegraphics[height=60mm]{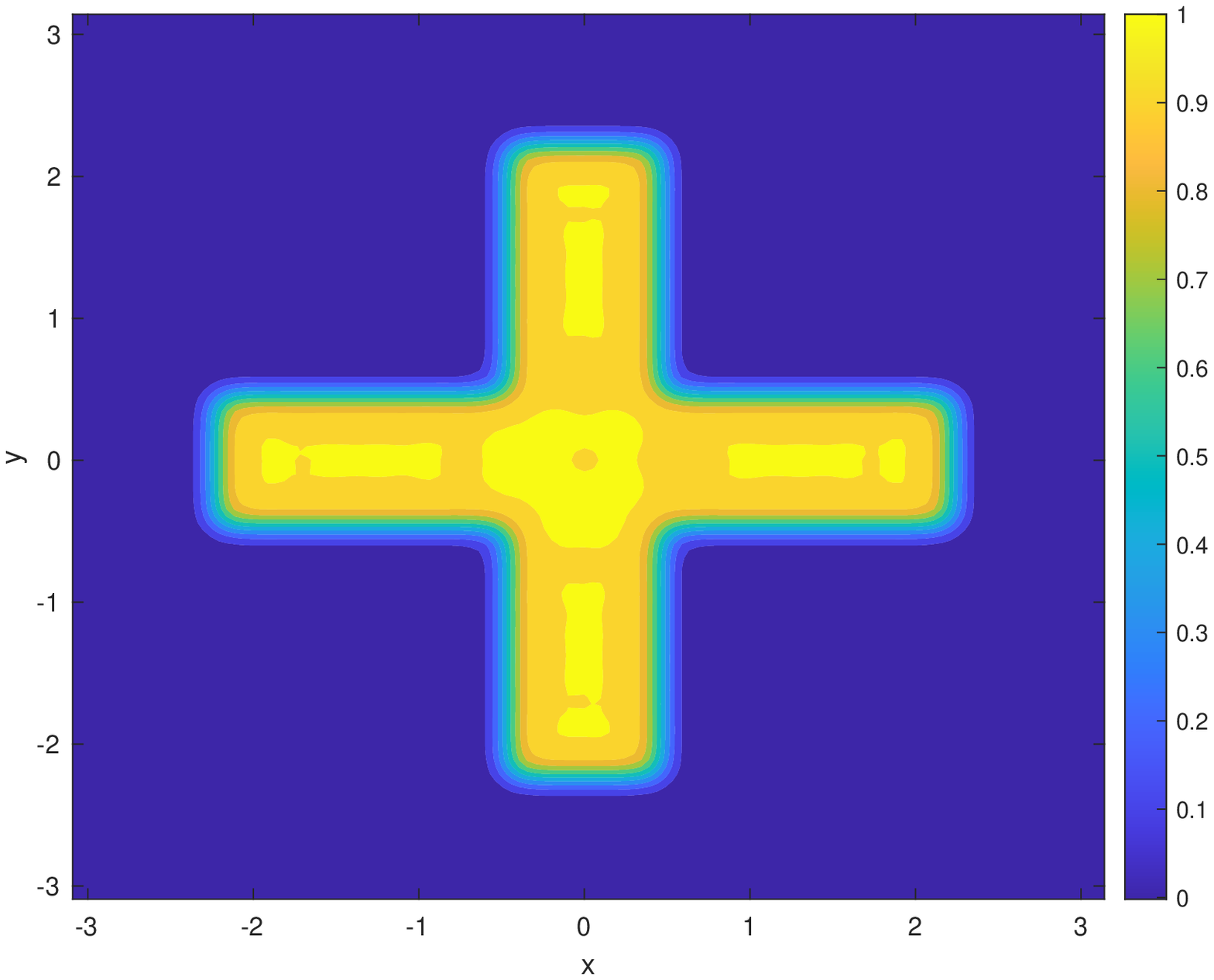}}
	\caption{Example \ref{ex:adv} with a cross-shaped initial condition. $d=2$.  $\varepsilon=10^{-5}$. $t=2\pi$. $N=128$. (a) Initial condition. (b) Numerical solution at $t=2\pi$.}
	\label{fig:adv2d}
\end{figure}
	\end{exa} 

\begin{exa}\label{ex:rot} For this example, we simulate the following 2D solid-body-rotation problem
	$$
	u_t - yu_x + xu_y = 0,\quad (x,y)\in[-\pi.\pi]^2
	$$
and compare the performance of two approaches. In particular, we consider two sets of initial conditions, including a smooth Gaussian hump
$$
u(x,v,t=0) = \exp(-x^2-5y^2),
$$
 and the non-smooth cross-shaped function considered in previous example. For the smooth initial condition, we report the convergence study of approach I for the $L^2$ error and orders of accuracy in Table \ref{ex:rot}. We employ the fifth order linear FD method with splitting for spatial discretization and let $\varepsilon=10^{-7}$ for truncation. We compute the solution up to $t=2\pi$, i.e. one full evolution. Second order accuracy is observed as expected. In addition, the numerical rank of the solutions remains relatively low and independent of mesh size $N$. Similar convergence behavior  is observed for approach II and hence omitted for brevity. We then consider the discontinuous initial condition and employ the fifth order WENO FD method with splitting for spatial discretization. In Figure \ref{fig:rot}, we report the numerical results and the time evolution of the numerical rank by approach I with mesh size $N=128$ up to $t=2\pi$. It is observed that the numerical solution suffers severe representation rank explosion, which is more pronounced with a finer mesh.  Furthermore, the numerical solution exhibits spurious oscillations due to the Gibbs phenomenon, and due to global nature of the basis  in the low-rank representation, the WENO methodology becomes less effective in controlling spurious oscillations.  In Figure \ref{fig:rot_II}, we report the contour plots of the solution together with the time history of the numerical rank of flow map $\mathcal{X}^*$. Note that the flow maps are smooth and remain very low-rank regardless of the solution profile. It is observed that the proposed method is able to capture the low-rank structure of the flow map and  high quality results are obtained. In Table \ref{tb:rot2d_II}, we report the convergence study of approach II for $\mathcal{X}^*$ and second order convergence is observed as expected.
 
\begin{table}[!hbp]
	\centering
	\caption{Example \ref{ex:rot}, $d=2$. $t=2\pi$. Approach I. }
	\label{tb:rot2d}
	\begin{tabular}{|c|c|c|c|}
		\hline
		$N$ & $L^2$-error & order  & rank\\\hline
16	&	3.29E-02	&		&	16	\\	\hline
32	&	6.16E-03	&	2.42	&	17	\\	\hline
64	&	3.85E-04	&	4.00	&	18	\\	\hline
128	&	2.88E-05	&	3.74	&	22	\\	\hline
256	&	9.07E-06	&	1.66	&	28	\\	\hline

	\end{tabular}
\end{table}

 \begin{figure}[h!]
 	\centering
 	\subfigure[]{\includegraphics[height=60mm]{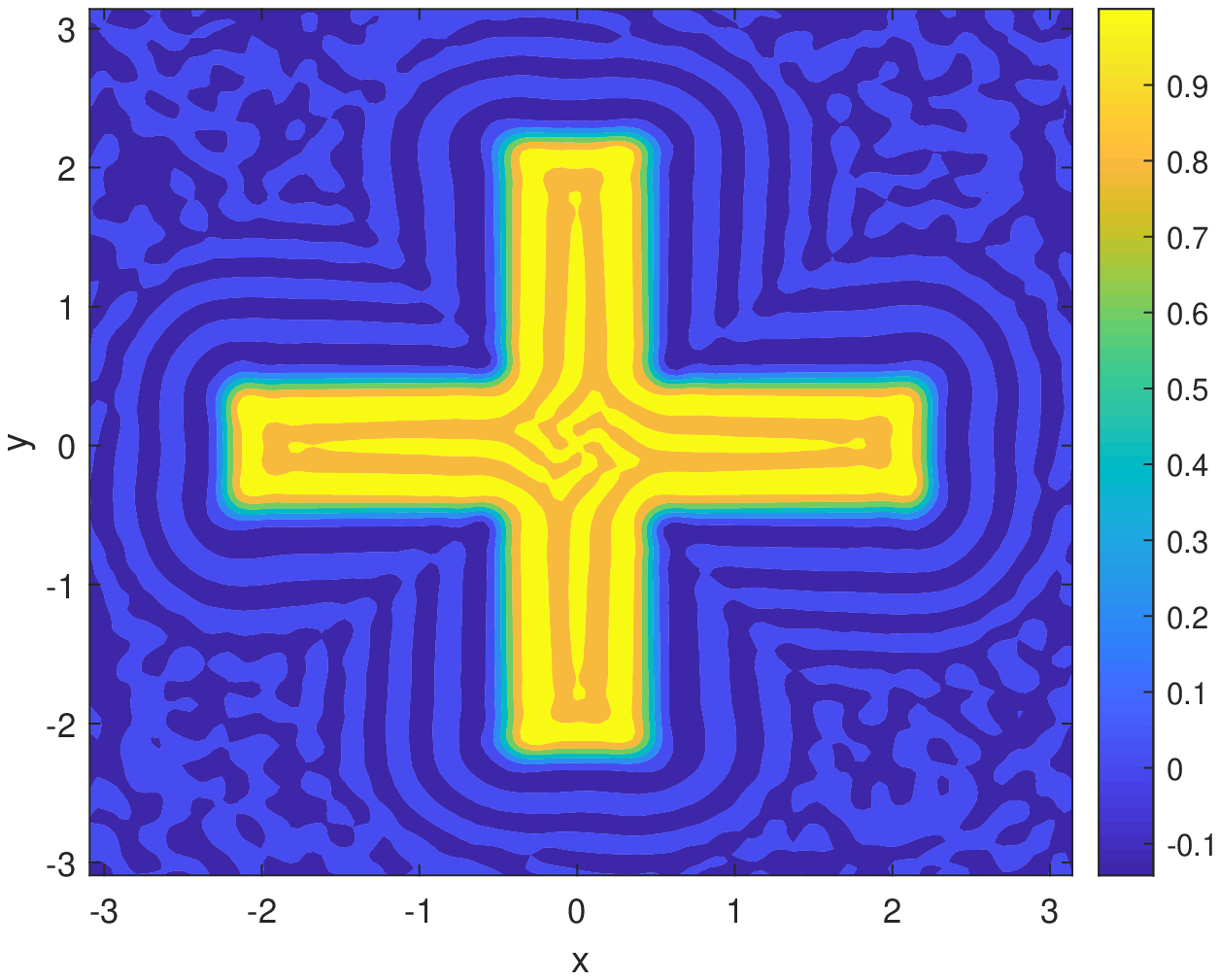}}
 	\subfigure[]{\includegraphics[height=60mm]{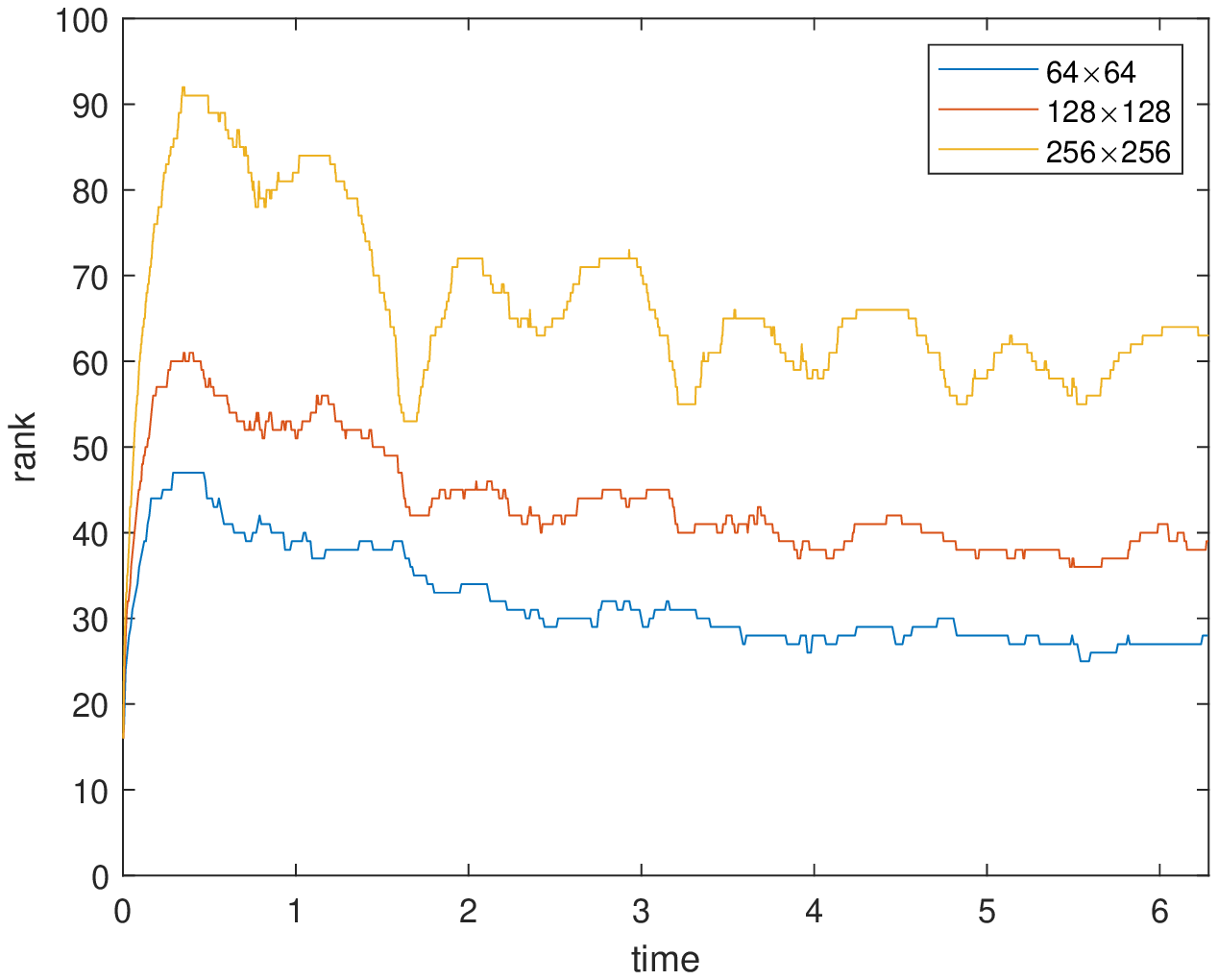}}
 	\caption{Example \ref{ex:rot}. $N = 128$, $\varepsilon=10^{-5}$. (a) Contour plot of the numerical solution by approach I, $t=2\pi$. (b) The time history of the numerical ranks. \label{fig:rot}}
 \end{figure}
 
  \begin{figure}[h!]
 	\centering
 	\subfigure[]{\includegraphics[height=60mm]{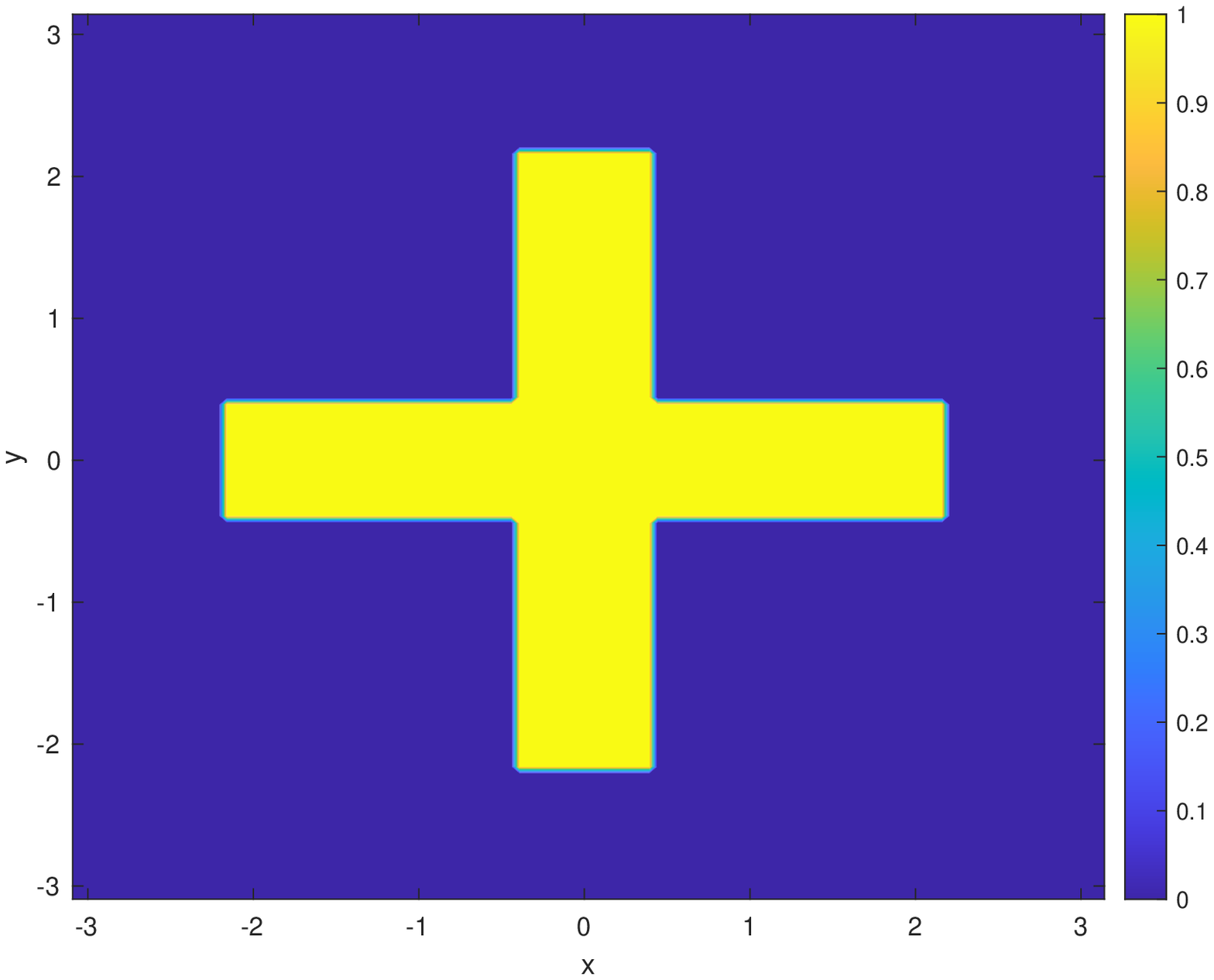}}
 	\subfigure[]{\includegraphics[height=60mm]{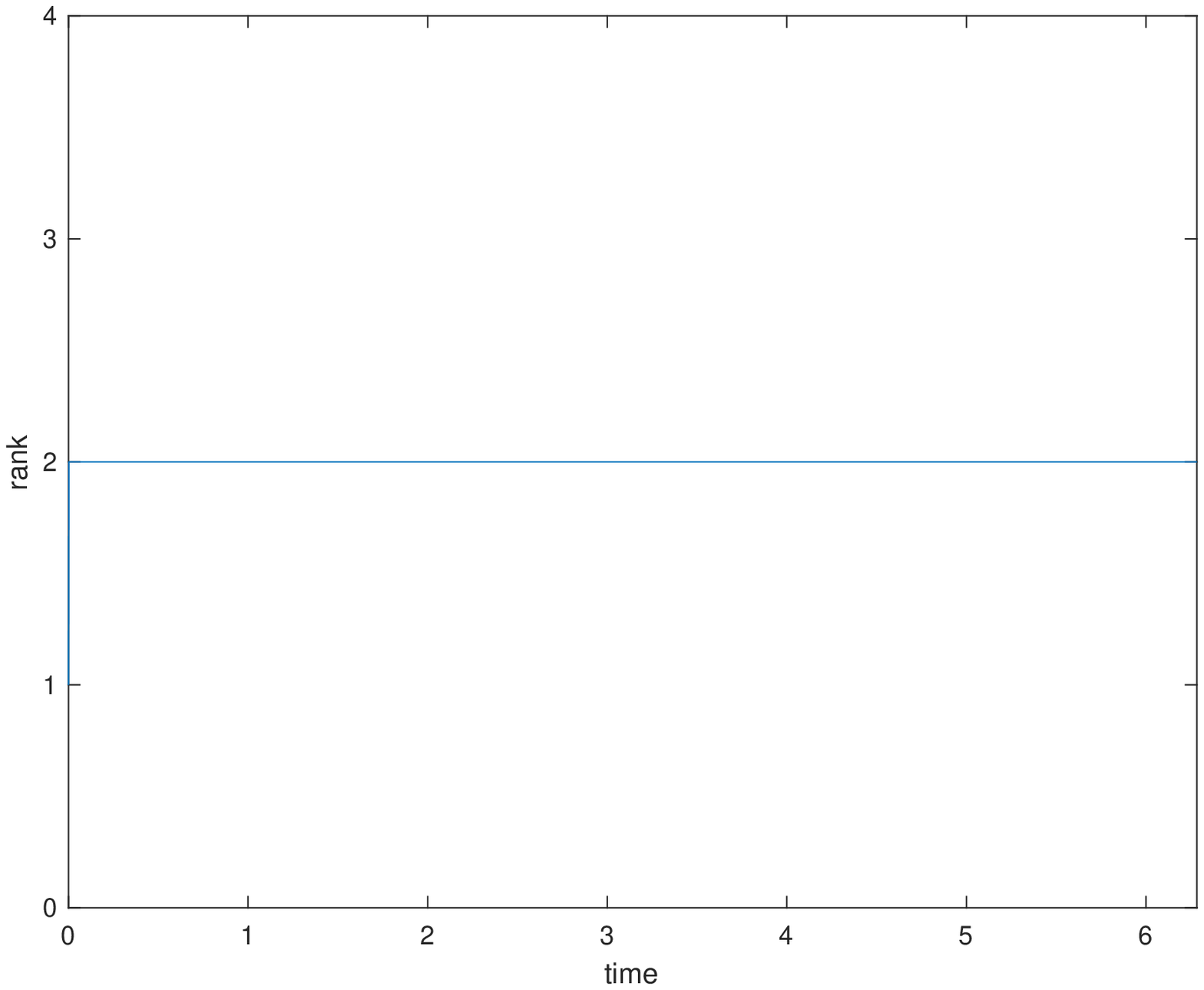}}
 	\caption{Example \ref{ex:rot}. $N=128$, $\varepsilon=10^{-5}$, $t=2\pi$. (a) Contour plot of the numerical solution by approach II. (b) The time history of the numerical rank of $\mathcal{X}^*$. \label{fig:rot_II}}
 \end{figure}

 \begin{table}[!hbp]
	\centering
	\caption{Example \ref{ex:rot}, $d=2$. $t=2\pi$. Approach II. Convergence for $\mathcal{X}^*$. }
	\label{tb:rot2d_II}
	\begin{tabular}{|c|c|c|c|}
		\hline
		$N$ & $L^2$-error & order  & rank\\\hline
16	&1.47E-03	&		& 2\\ \hline
32	&3.66E-04	&2.00&2\\ \hline
64	&9.15E-05	&2.00&2\\ \hline
128	&2.29E-05	&2.00&2\\ \hline

	\end{tabular}
\end{table}

\end{exa} 

\begin{exa}\label{ex:defor} We consider the swirling deformation flow, governed by the linear transport equation
$$
u_{t}-\left(\cos ^{2}\left(\frac{x}{2}\right) \sin (y) g(t) u\right)_{x}+\left(\sin (x) \cos ^{2}\left(\frac{y}{2}\right) g(t) u\right)_{y}=0, \quad(x, y) \in[-\pi, \pi]^{2}
$$
where $g(t) =\pi\cos\left(\frac{\pi t}{T}\right)$, and $T=1.5$, with periodic boundary conditions.
 We consider the same two initial conditions as in the previous example. Note that the solution profile would be deformed along the flow maps and return to its initial state at time $t=T$. Note that the underlying flow maps remain smooth and low-rank over time. We let $\varepsilon=10^{-5}$ for the truncation. The convergence study is summarized in Table \ref{tb:defor2d} for approach I.  Second order of convergence is observed. Then, we consider the discontinuous cross-shaped initial condition.  Similar to the solid body rotation example, the solution develops misaligned discontinuous  structures which are high-rank by nature. In Figure \ref{fig:defor_flow}, we report the contour plots of the solutions at $t=T/2$ and $t=T$ with mesh size $N_x\times N_y$. The time evolution of the numerical rank of the solution is reported in the first plot of Figure~\ref{fig:defor_rank}. It is observed that the method suffers severe rank explosion, and the numerical solution develops spurious oscillations due to the Gibbs phenomenon. Meanwhile,  since the flow maps are smooth and of low rank, approach II is able to generate a high quality result without oscillations, see Figure \ref{fig:defor_flow}. In addition, the numerical rank remains very low, leading to significant computation savings, see the second plot in Figure \ref{fig:defor_rank}. 
\begin{table}[!hbp]
	\centering
	\caption{Example \ref{ex:defor}, $d=2$. $t=T$. Approach I. }
	\label{tb:defor2d}
	\begin{tabular}{|c|c|c|c|}
		\hline
		$N$ & $L^2$-error & order  & rank\\\hline
16	&	1.90E-02	&		&	16	\\\hline
32	&	3.00E-03	&	2.66	&	17	\\\hline
64	&	1.92E-04	&	3.96	&	18	\\\hline
128	&	7.30E-06	&	4.72	&	22	\\\hline
256	&	7.32E-07	&	3.32	&	28	\\\hline	
	\end{tabular}
\end{table}

 \begin{figure}[h!]
	\centering
	\subfigure[]{\includegraphics[height=60mm]{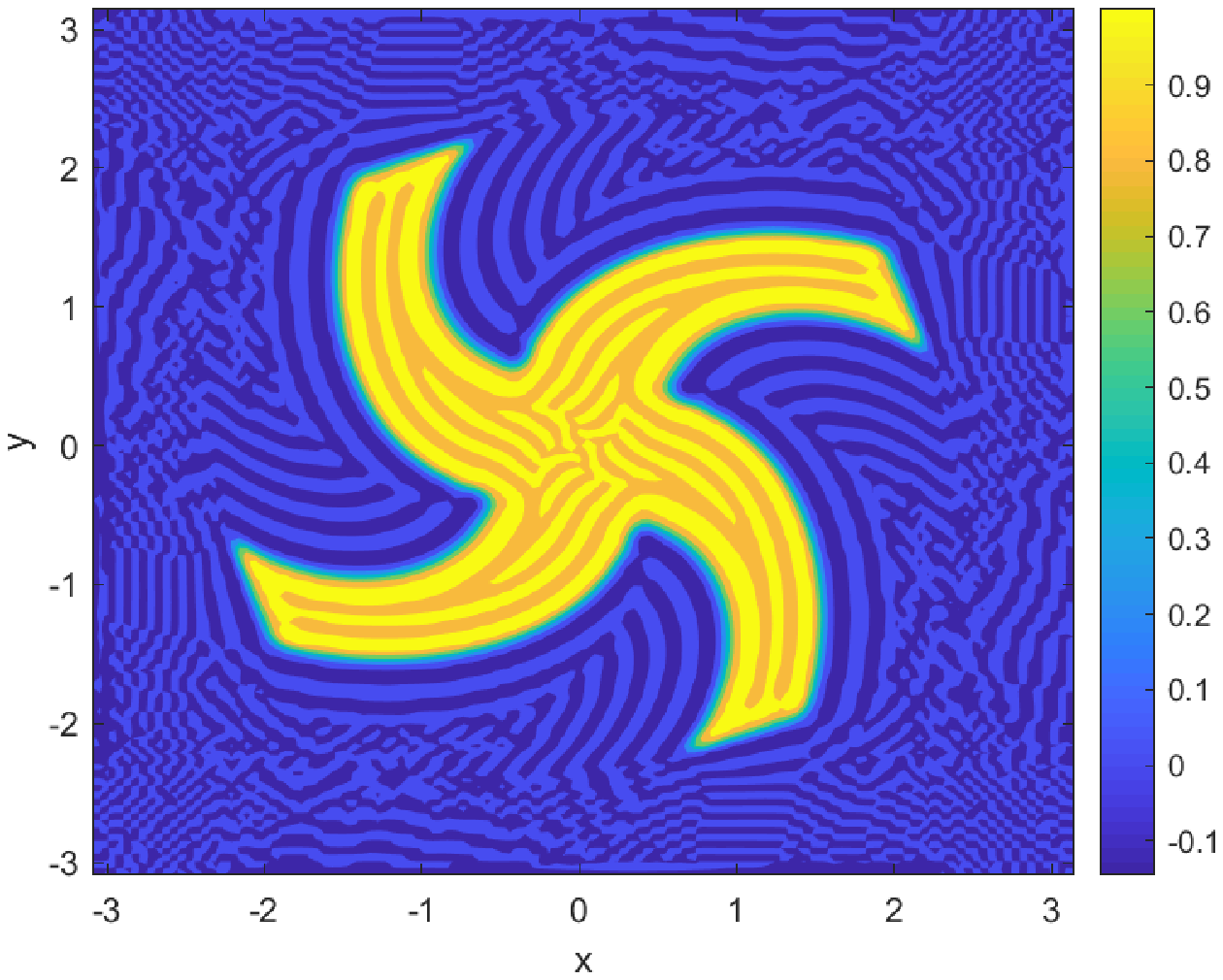}}
	\subfigure[]{\includegraphics[height=60mm]{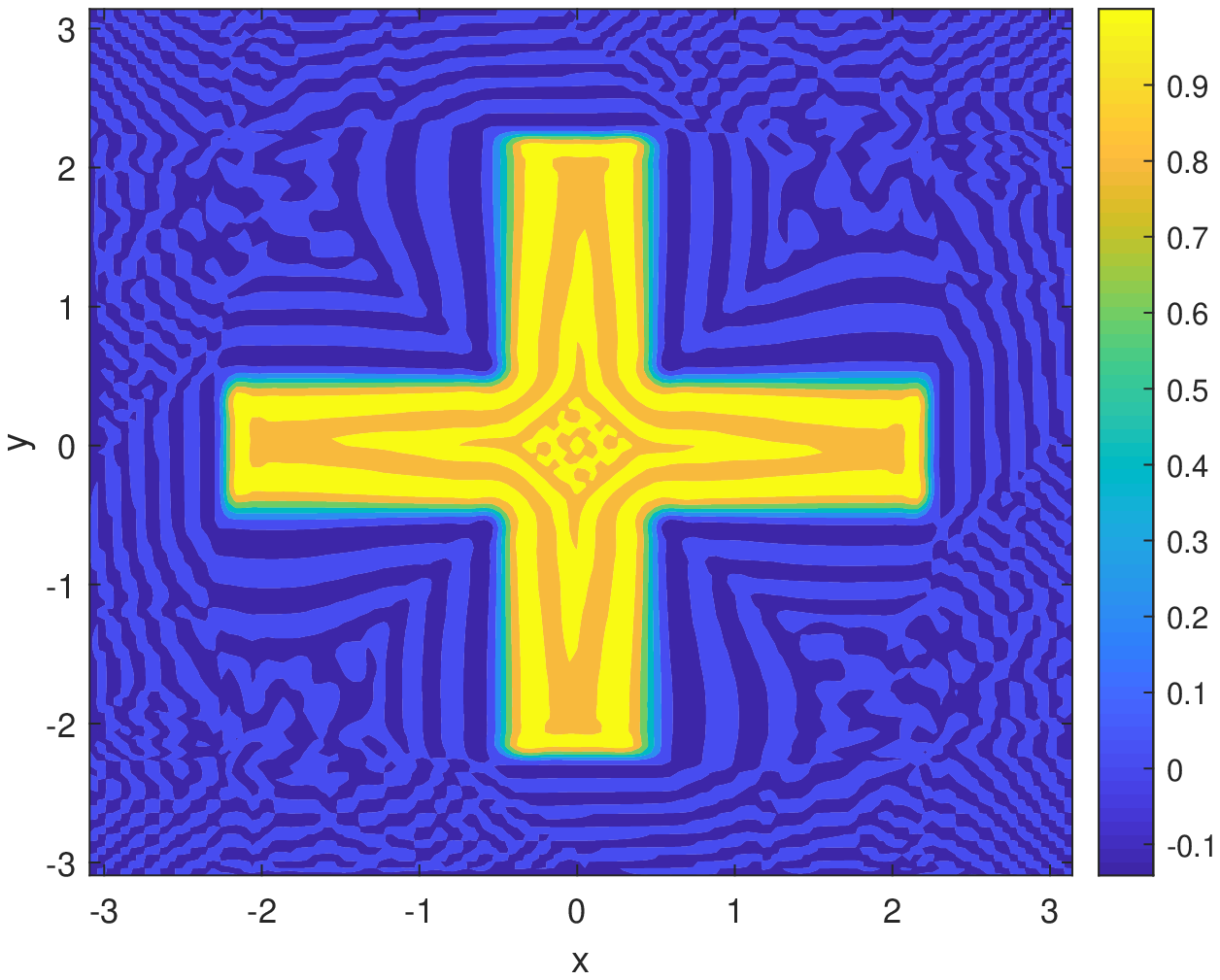}}
	\caption{Example \ref{ex:defor}. $N_x\times N_y=128\times128$, $\varepsilon=10^{-5}$. (a) Contour plot of the numerical solution by approach I, $t=T/2$. (b) Contour plot of the numerical solution by approach I, $t=T$. \label{fig:defor}}
\end{figure}	
 \begin{figure}[h!]
	\centering
	\subfigure[]{\includegraphics[height=60mm]{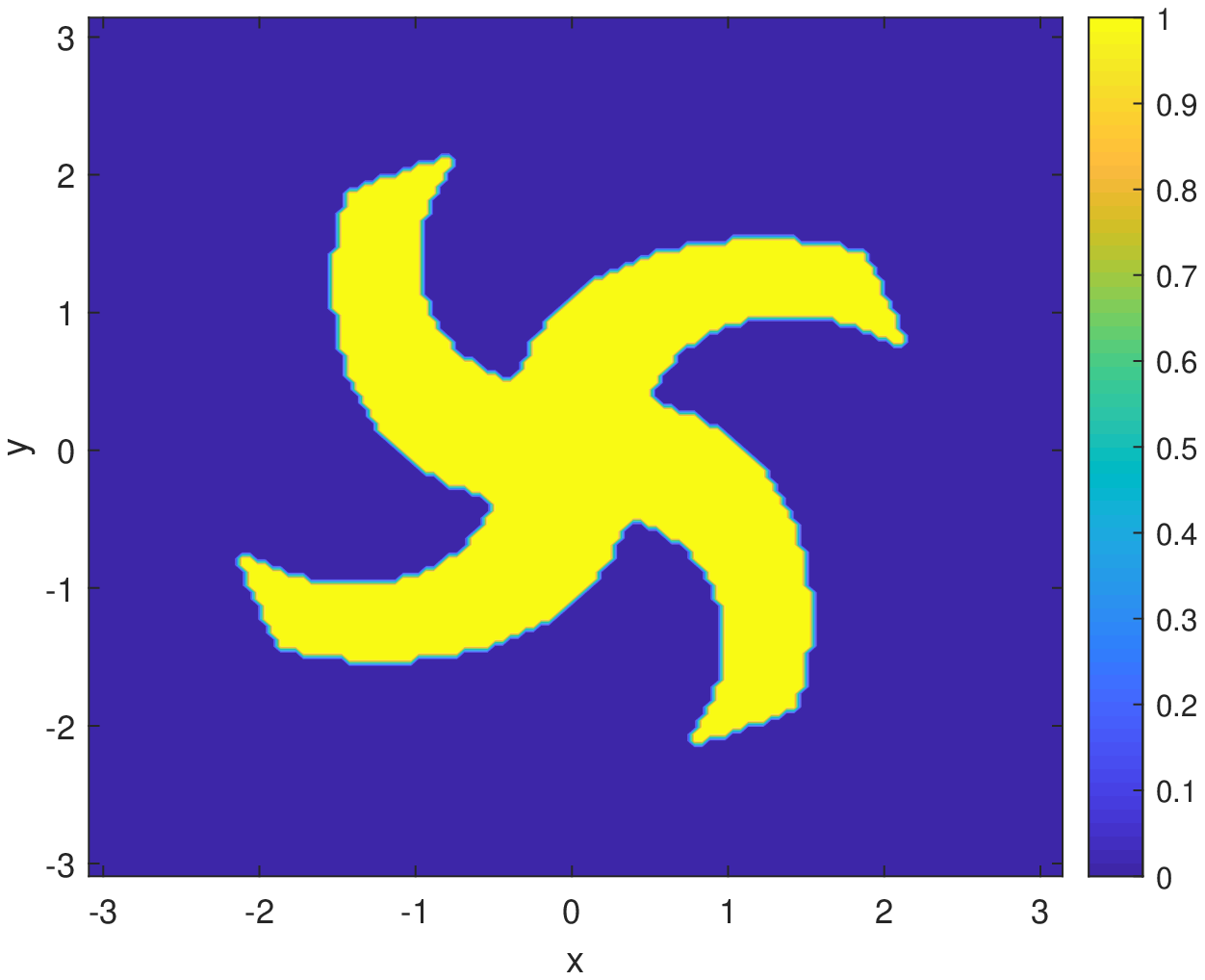}}
	\subfigure[]{\includegraphics[height=60mm]{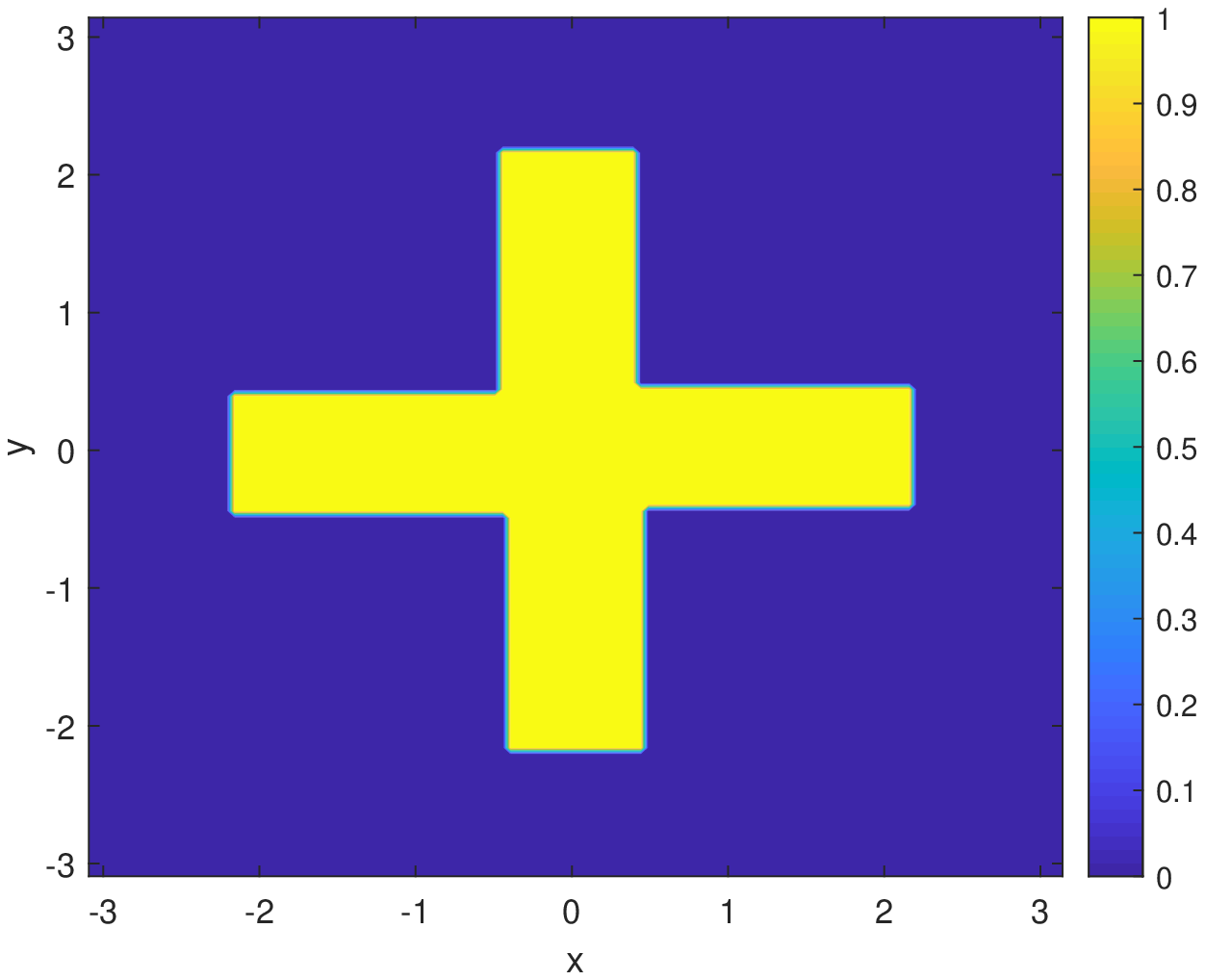}}
	\caption{Example \ref{ex:defor}. $N_x\times N_y=128\times128$, $\varepsilon=10^{-5}$. (a) Contour plot of the numerical solution by flow map approach, $t=T/2$. (b) Contour plot of the numerical solution by  flow map approach, $t=T$. (c) The time history of the numerical ranks of the solution. \label{fig:defor_flow}}
\end{figure}

 \begin{figure}[h!]
	\centering
\subfigure[]{\includegraphics[height=60mm]{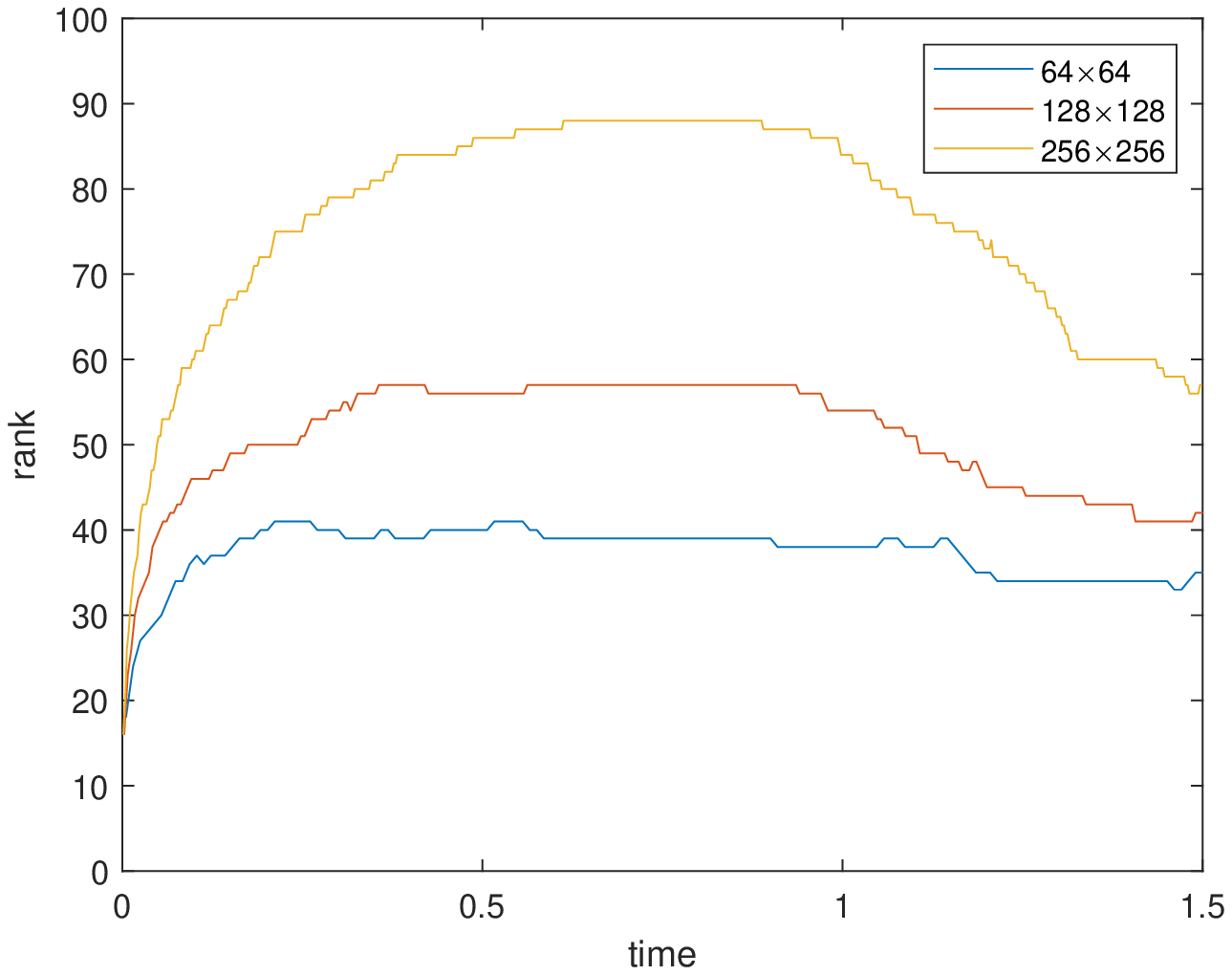}}
\subfigure[]{\includegraphics[height=60mm]{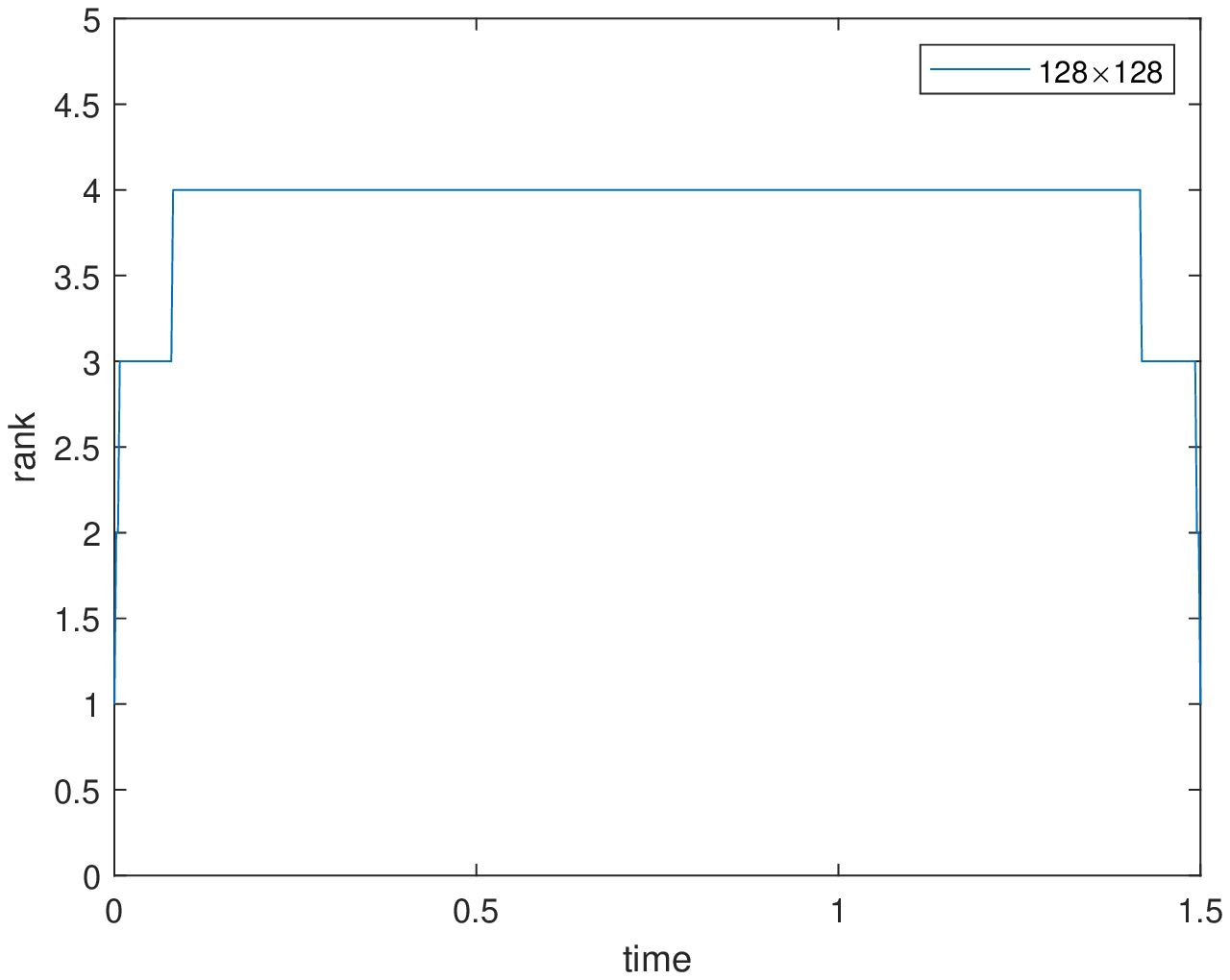}}
		\caption{Example \ref{ex:defor}.  The time history of the numerical ranks of the solution for Approach I (left) and II (right). \label{fig:defor_rank}}
\end{figure}

\end{exa}

\subsection{The 1D1V Vlasov-Poisson system}
\begin{exa}\label{ex:linearvp} Consider a linear 1D1V VP system with a given electric field $E(x)=\frac12\sin(\pi x)$. The initial condition is a smooth Gaussian hump
$$
u(x,v,t=0) = \exp(-20(x^2+v^2)).
$$ 
We test the accuracy for both approaches for $t=0.75$ and $t=1.5$. The reference solution is computed by solving the characteristic equation with very small time step. Note that the solution profile will deform into a thinner and thinner structure over time.  In Figures \ref{tb:linearVPaccuI}-\ref{tb:linearVPaccuII}, we report the convergence study by setting $\varepsilon=10^{-5}$ and $10^{-7}$. It is observed that, for approach I, once the solution is well resolved, high order accuracy can be observed. In the meantime, the numerical rank of the solution become larger and larger due to the deformation of the solution, and the error magnitude increases by comparing the results at $t=0.75$ and $t=1.5$. 
On the other hand,  the underlying flow maps are smoother and lower rank than the solution. The approach II is able to capture the very low-rank structure of the flow map and generate more accurate approximation. In particular, by choosing a large truncation parameter, i.e. $\varepsilon=10^{-5}$, the error from the discretization is dominated by the truncation error in the tensor decomposition, as opposed to approach I. Furthermore,
the truncation error accumulates over the time integration, observing that error increases the mesh is refined. If a smaller truncation parameter, i.e. $\varepsilon=10^{-7}$, then the discretization error dominates the truncation error, and high order convergence for the approximations of the flow map as well as the solution is observed. 
\end{exa}
\begin{table}[!hbp]
	\centering
	\caption{Example \ref{ex:linearvp}, $d=2$. Approach I.}
	\label{tb:linearVPaccuI}
	\begin{tabular}{|c|c|c|c|c|c|c|}
		\hline
 &	\multicolumn{3}{c|}{$\varepsilon=1e-5$}	& 	\multicolumn{3}{c|}{$\varepsilon=1e-7$}\\\hline
  &	\multicolumn{6}{c|}{$t=0.75$}\\\hline
		$N$ & L$^\infty$-error & order & rank & L$^\infty$-error & order & rank\\\hline
16	&	9.01E-02	&		&	13        & 9.01E-02	&		&	16    \\\hline
32	&	1.36E-02	&	2.73	&	19    &1.36E-02	&	2.73	&	25    \\\hline
64	&	8.65E-04	&	3.98	&	21   & 8.17E-04	&	4.06	&	31     \\ \hline
128	&	3.16E-04	&	1.45	&	20   &   4.67E-05	&	4.13	&	33     \\\hline
  &	\multicolumn{6}{c|}{$t=1.5$}\\\hline
  		$N$ & L$^\infty$-error & order & rank & L$^\infty$-error & order & rank\\\hline
16	&	2.94E-01	&		&	16		&		2.94E-01	&		&	16\\\hline
32	&	2.46E-01	&	0.26	&	25	& 		2.46E-01	&	0.26	&	31\\\hline
64	&	6.78E-02	&	1.86	&	41	& 		6.78E-02	&	1.86	&	48\\\hline
128	&	6.80E-03	&	3.32	&	45	& 		6.70E-03	&	3.34	&	67\\\hline	
	\end{tabular}
\end{table}

\begin{table}[!hbp]
	\centering
	\caption{Example \ref{ex:linearvp}, $d=2$. Approach II.}
	\label{tb:linearVPaccuII}
	\begin{tabular}{|c|c|c|c|c|c|c|c|c|}
		\hline
		&	\multicolumn{3}{c|}{$\varepsilon=1e-5$}	& 	\multicolumn{5}{c|}{$\varepsilon=1e-7$}\\\hline
		&	\multicolumn{8}{c|}{$t=0.75$}\\\hline
		$N$ & error of $X$ & error of $f$  & rank & error of $X$ & order & error of $f$ & order & rank\\\hline
16	&	8.12E-05		&	2.13E-04	&		8	&	5.72E-05	&		&	2.11E-04	&		&	8	\\\hline
32	&	1.43E-04	&			2.76E-04	&		7	&	1.43E-05	&	2.00	&	1.75E-05	&	3.59	&	10	\\\hline
64	&	2.23E-04	&			5.04E-04	&		7	&	3.98E-06	&	1.84	&	4.75E-06	&	1.88	&	10	\\\hline
128	&	3.90E-04	&			3.90E-04	&		7	&	1.29E-06	&	1.62	&	2.11E-06	&	1.17	&	10	\\\hline

&	\multicolumn{8}{c|}{$t=1.5$}\\\hline
		$N$ & error of $X$ & error of $f$  & rank & error of $X$ & order & error of $f$ & order & rank\\\hline
16	&	1.60E-03	&			1.32E-02	&		12	&	1.70E-03	&		&	1.31E-02	&		&	16	\\\hline
32	&	4.12E-04	&			1.20E-03	&		12	&	1.73E-04	&	3.29	&	1.00E-03	&	3.71	&	18	\\\hline
64	&	8.09E-04	&			1.60E-03	&		12	&	2.54E-05	&	2.77	&	3.56E-05	&	4.81	&	18	\\\hline
128	&	1.50E-03	&			2.60E-03	&		11	&	3.92E-06	&	2.70	&	5.21E-06	&	2.77	&	18	\\\hline

	\end{tabular}
\end{table}

\begin{figure}[h!]
	\centering
	\subfigure[]{\includegraphics[height=60mm]{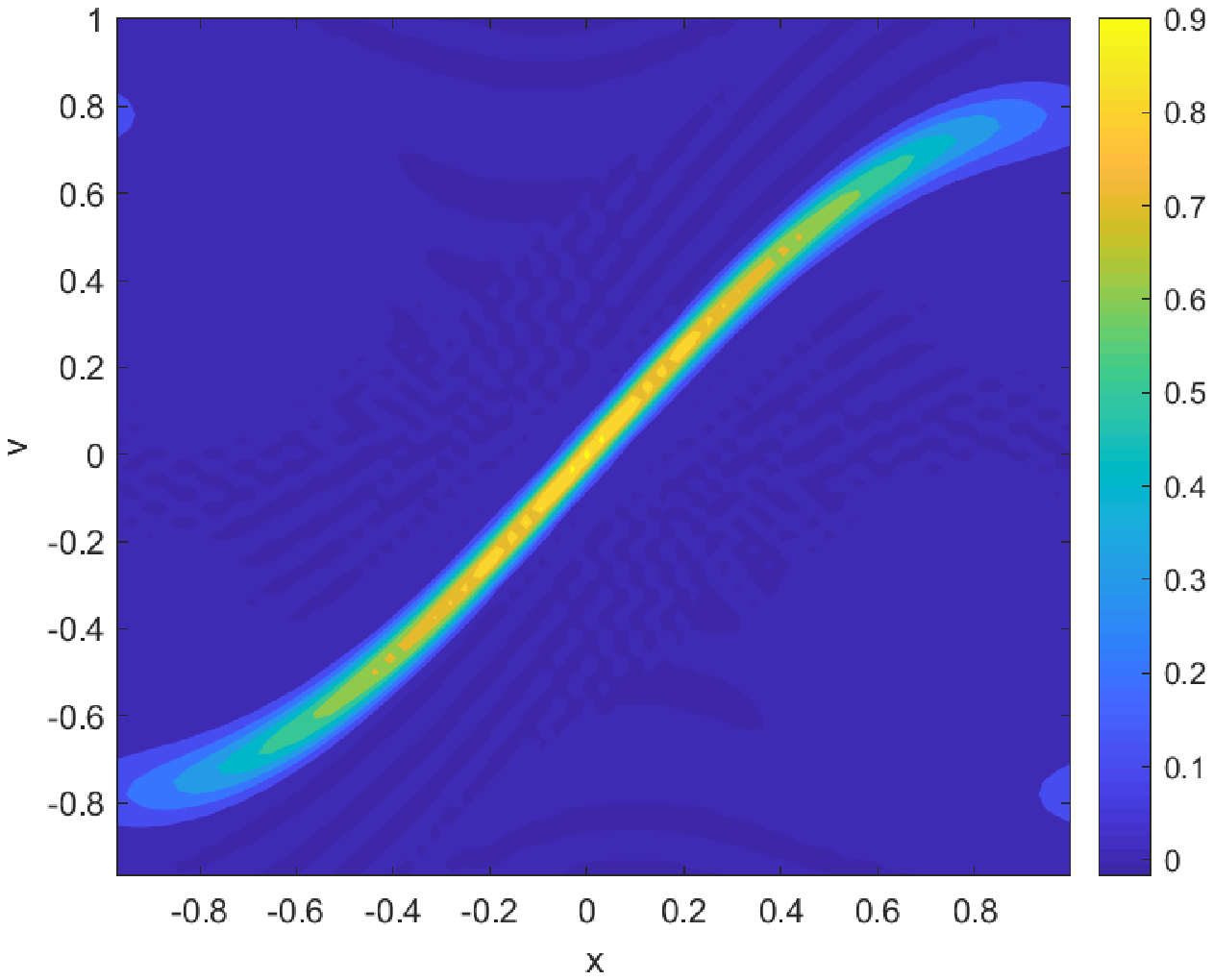}}
	\subfigure[]{\includegraphics[height=60mm]{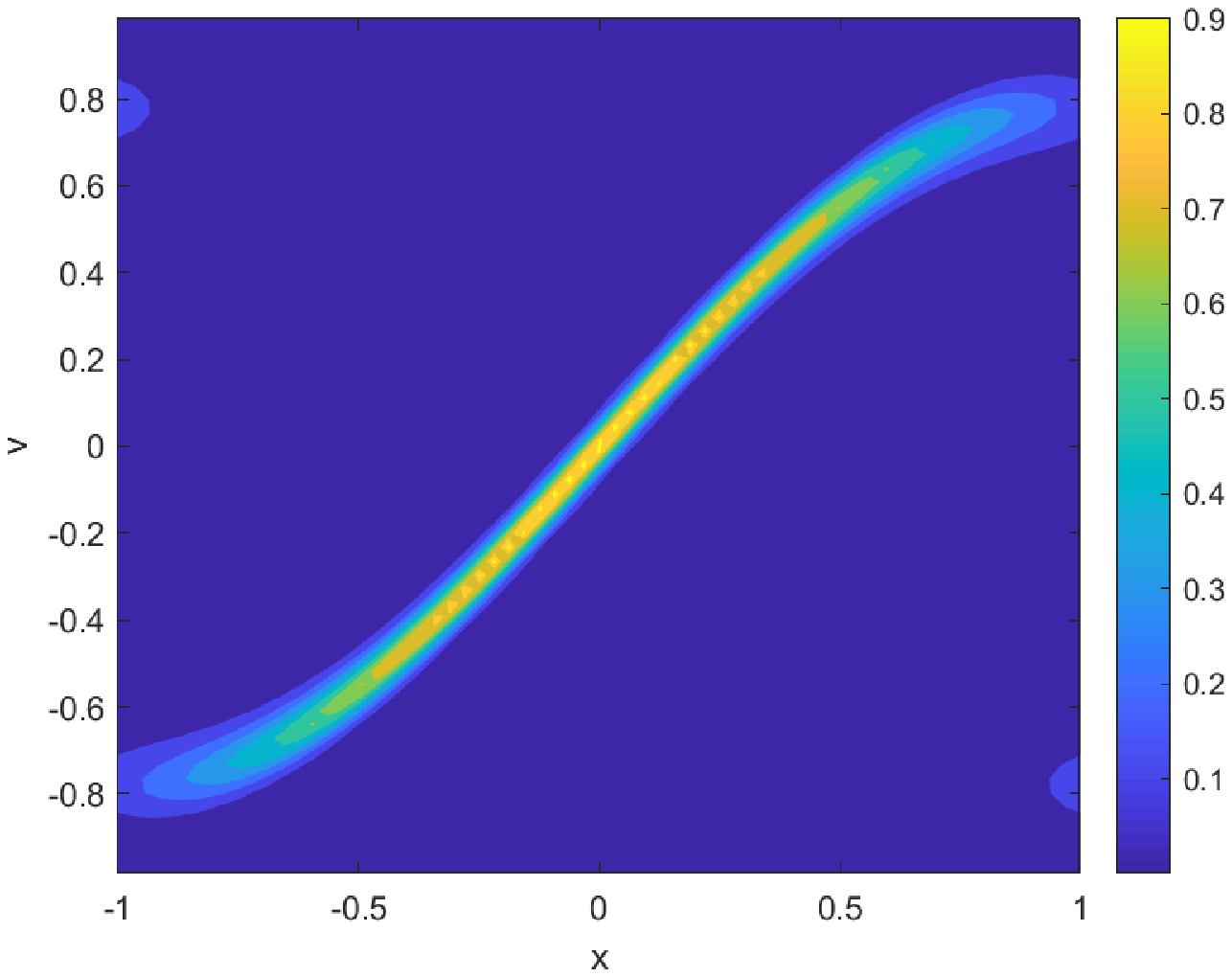}}
		\subfigure[]{\includegraphics[height=60mm]{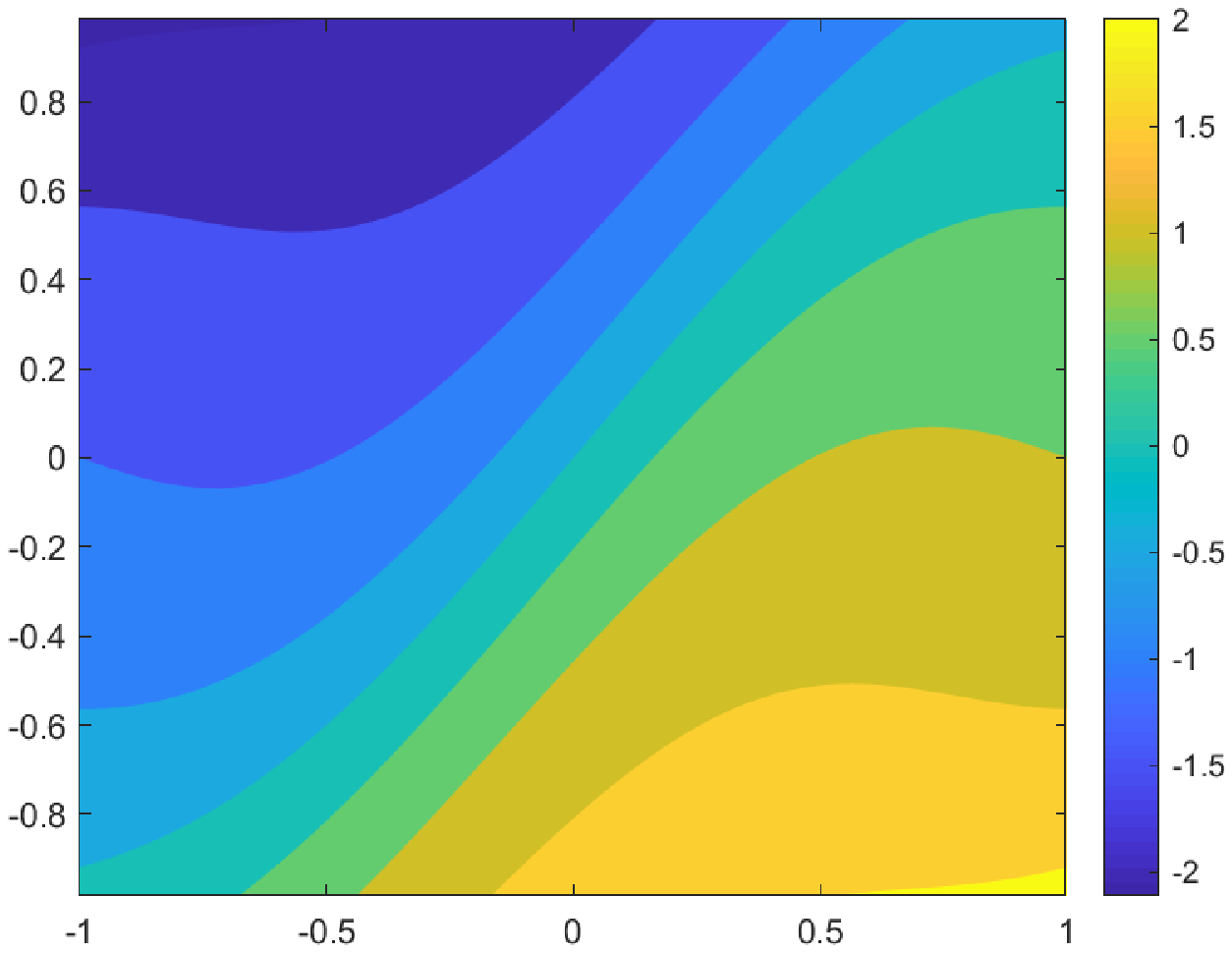}}
				\subfigure[]{\includegraphics[height=60mm]{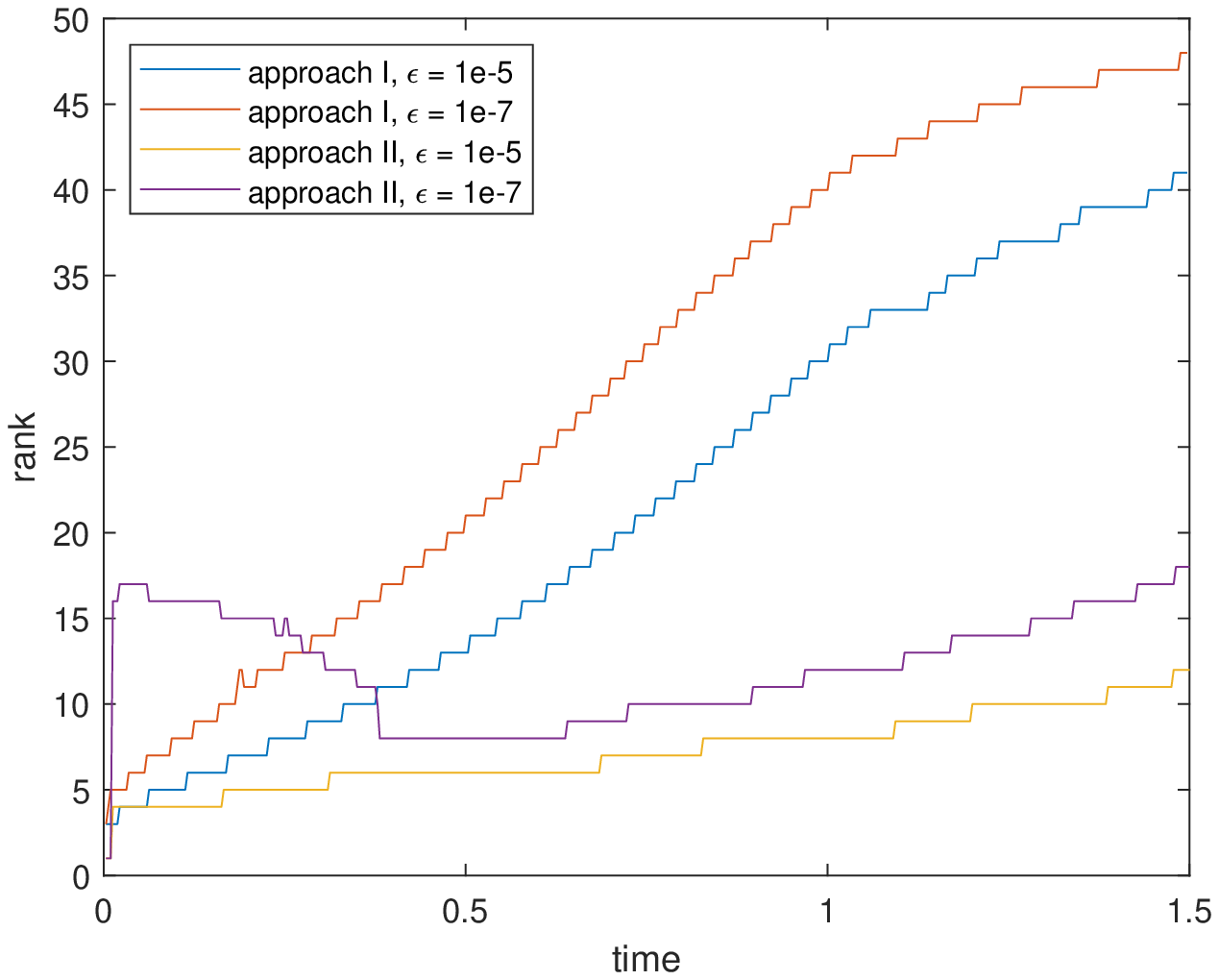}}
	\caption{Example \ref{ex:linearvp}. $N_x\times N_v=64\times64$. (a) Contour plot of the numerical solution by approach I, $t=1.5$, $\varepsilon=10^{-7}$. (b) Contour plot of the numerical solution by approach II, $t=1.5$, $\varepsilon=10^{-7}$. (c) Contour plot of flow map $\mathcal{X}^*$ by approach II, $t=1.5$, $\varepsilon=10^{-7}$. (d) The time evolution of the numerical ranks. \label{fig:linealVP}}
\end{figure}

\begin{exa}\label{ex:two} We consider the 1D1V two-stream instabilities with initial condition 
$$
f(x,v,t=0) = \frac{2}{7\sqrt{2\pi}}\left(1+5v^2\right)\left( 1+ \alpha\left(\left(\cos(2kx) + \cos(3kx) \right)/1.2 + \cos(kx) \right)\right)\exp\left(-\frac{v^2}{2}\right),
$$
where $\alpha = 0.01$ and $k=0.5$. We compare the performance of both approaches for the long term simulations. We set $\varepsilon=10^{-6}$. In Figure \ref{fig:two_his}, we plot the time evolution of the numerical rank of the solutions. It is observed that the rank of the flow map is much smaller than that of the solution when $t \le 20$. After that, the flow map rank starts growing to capture the underlying nonlinear dynamics, leading to improved efficiency.  In Figure \ref{fig:two_contour}, we report the contour plots of the solutions by both approaches at $t=30$ and $t=40$ with mesh size $N_x\times N_v=128\times 256$. Both methods can generate high quality results and filamentation structures are well captured. Furthermore, the flow map approach is able to resolve finer solution structures, which is ascribed to the fact that the flow maps are smoother than the solution itself for this example.    

\begin{figure}[h!]
	\centering
	\includegraphics[height=60mm]{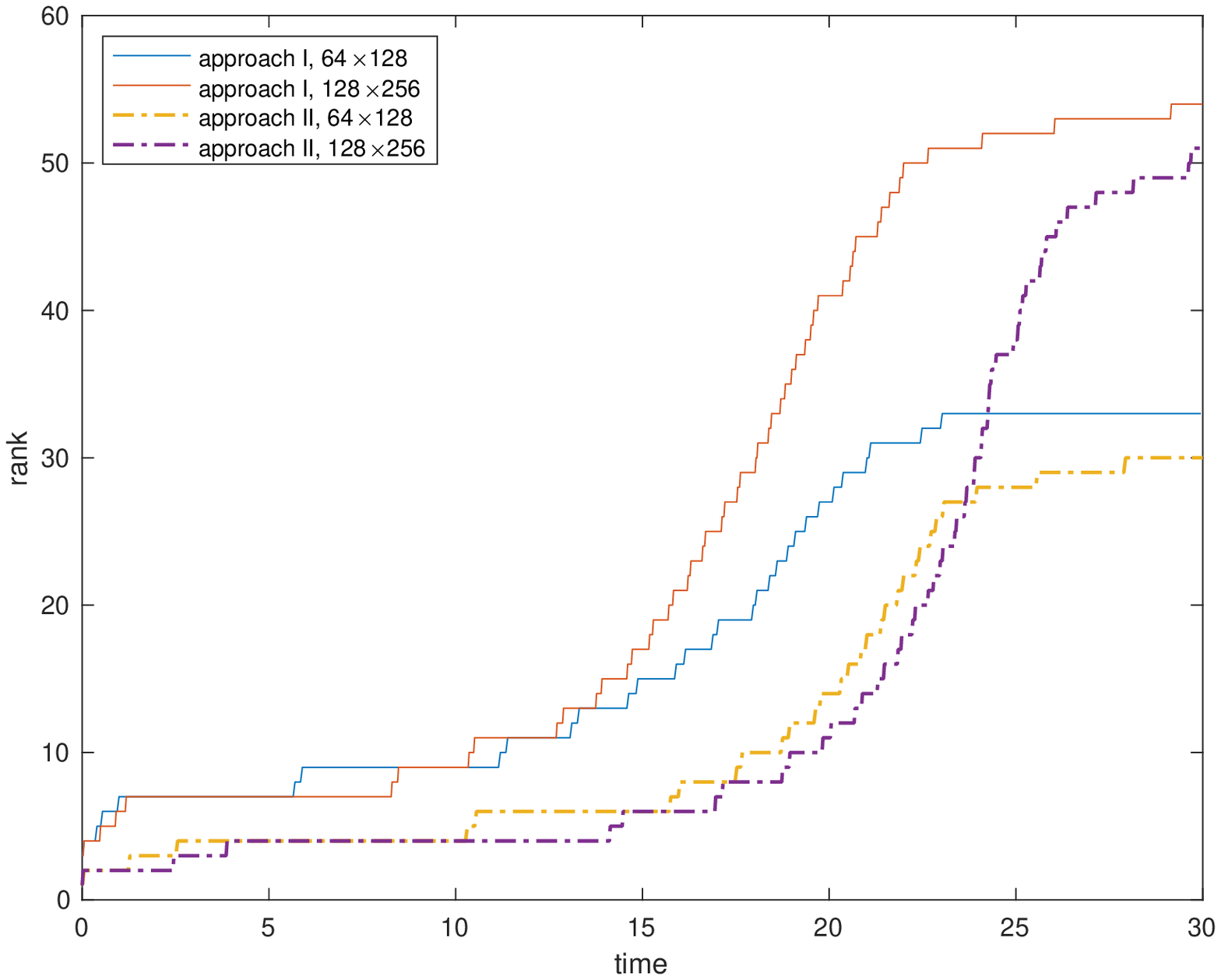}
	%\subfigure[]{\includegraphics[height=40mm]{weak2d_total_energy.eps}}
	\caption{Example \ref{ex:two}. Two-stream instabilities  1D1V. $d=1$. $\varepsilon=10^{-6}$. The time evolution of the numerical rank of the solutions.\label{fig:two_his}}
\end{figure}

\begin{figure}[h!]
	\centering
%	\subfigure[]{\includegraphics[height=60mm]{twostream1dii_WENO5_t30_64_128_contour_64.eps}}
%	\subfigure[]{\includegraphics[height=60mm]{twostream1dii_new_WENO_t30_64_128_contour_65.eps}}
	\subfigure[]{\includegraphics[height=60mm]{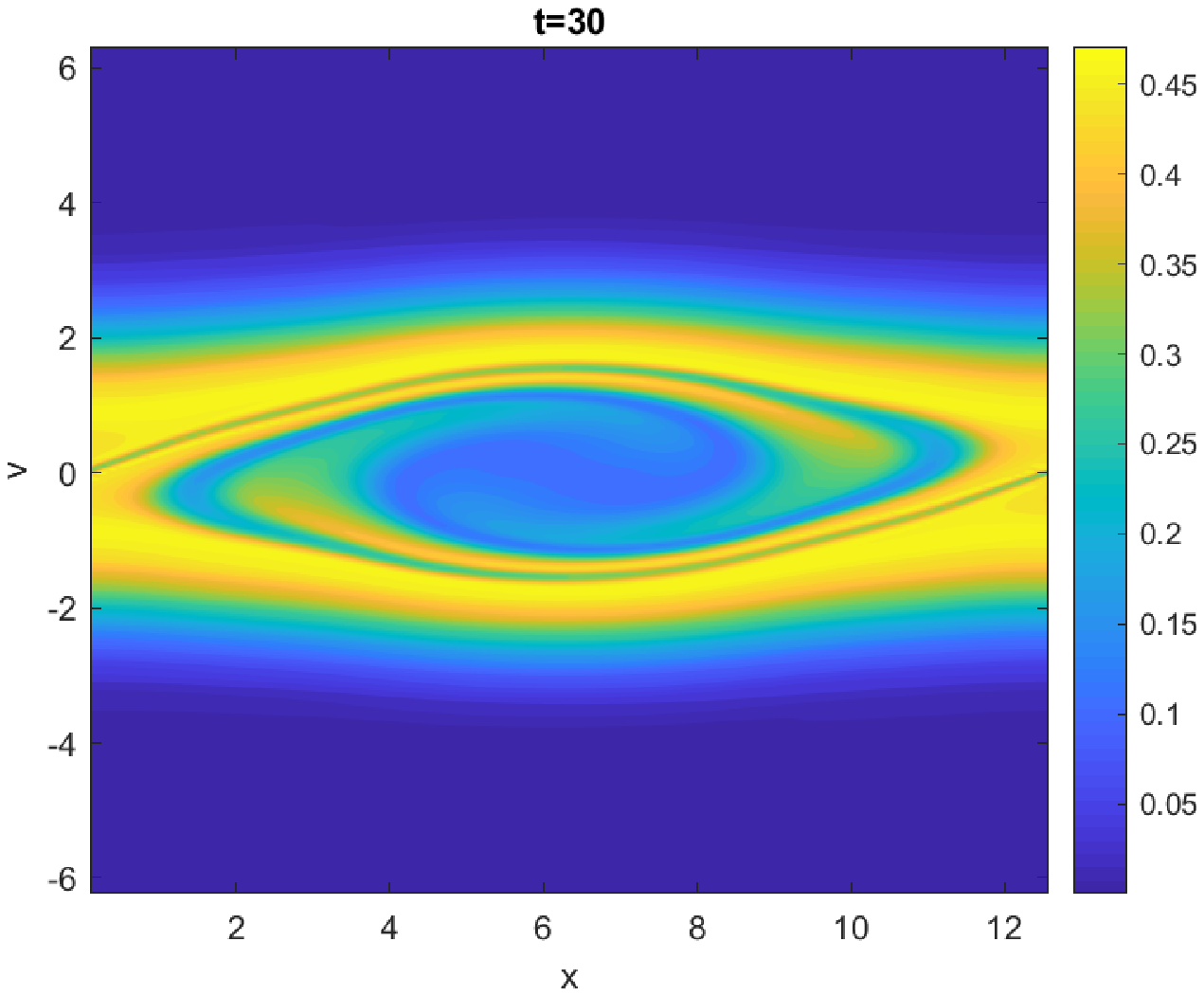}}
	\subfigure[]{\includegraphics[height=60mm]{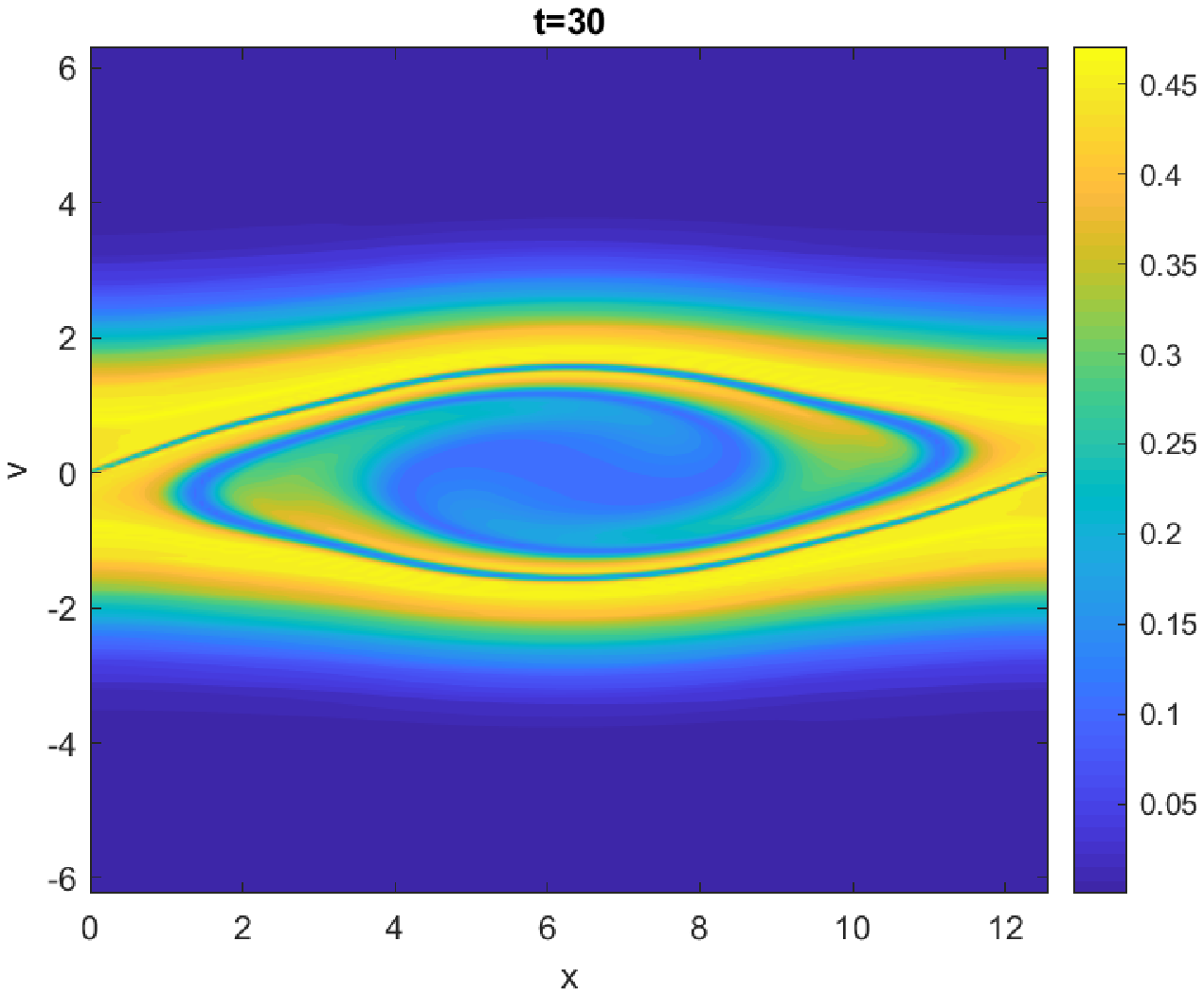}}
%	\subfigure[]{\includegraphics[height=60mm]{twostream1dii_WENO5_t40_64_128_contour_64.eps}}
%	\subfigure[]{\includegraphics[height=60mm]{twostream1dii_new_WENO_t40_64_128_contour_65.eps}}
	\subfigure[]{\includegraphics[height=60mm]{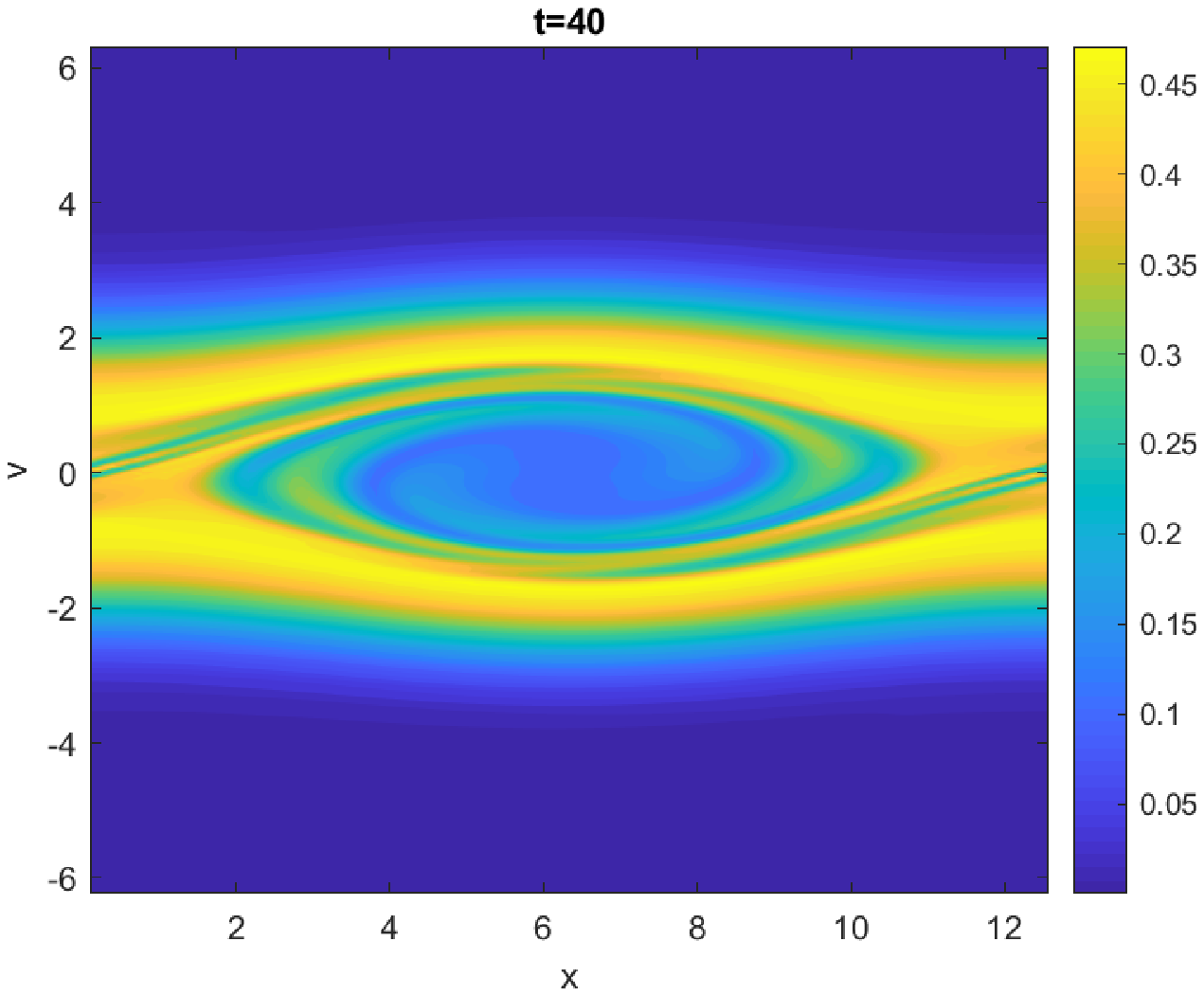}}
	\subfigure[]{\includegraphics[height=60mm]{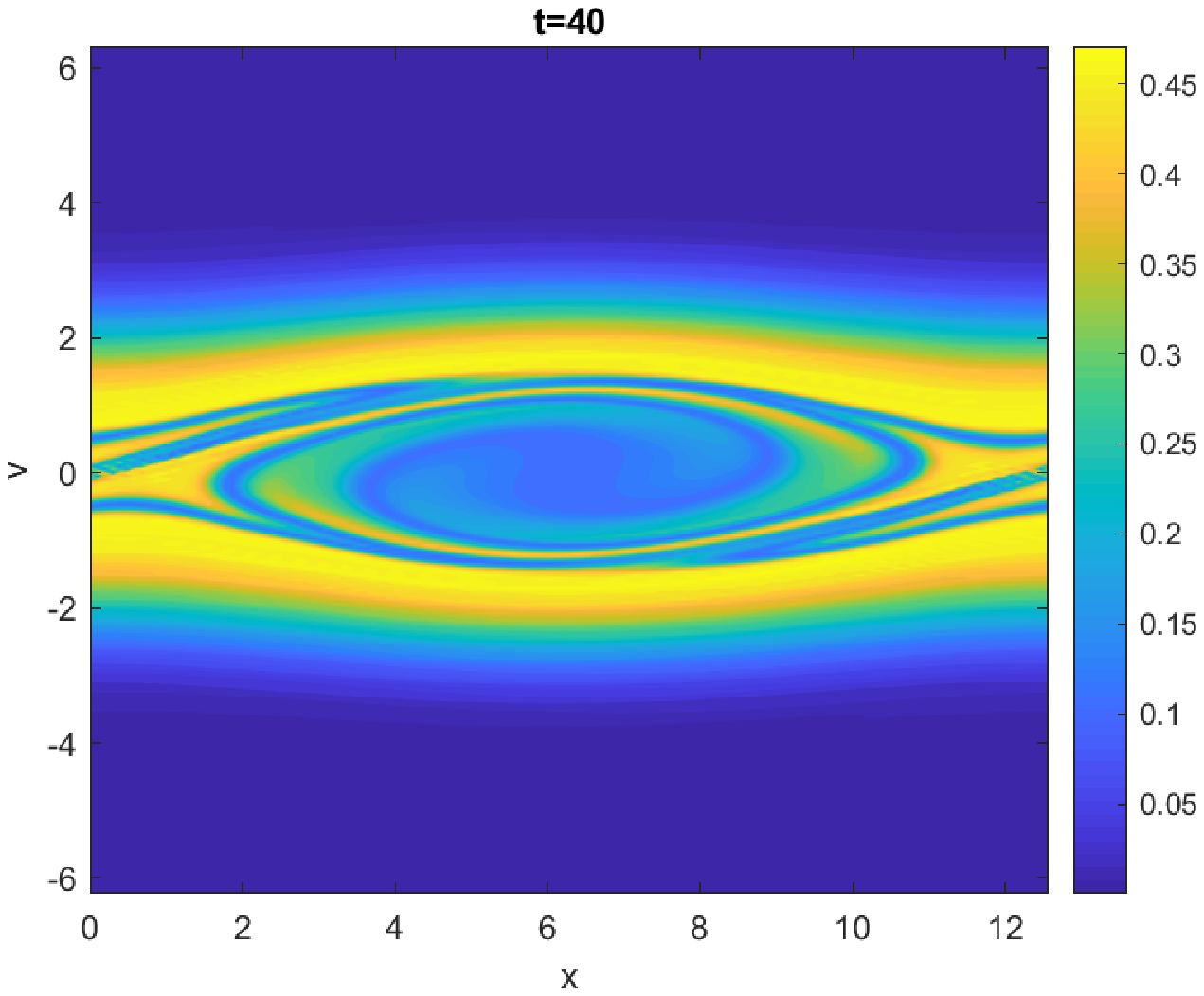}}

	%\subfigure[]{\includegraphics[height=40mm]{weak2d_total_energy.eps}}
	\caption{Example \ref{ex:two}. Two-stream instabilities  1D1V.  Contour plot of the solution. $\varepsilon=10^{-6}$. $N_x\times N_v=128\times256$.
 (a) approach I,  $t=30$. (b) approach II, $t=30$.
  (c) approach I, $t=40$. (d) approach II, $t=40$.
	\label{fig:two_contour}}
\end{figure}
\end{exa}

\subsection{The nonlinear 2D2V Vlasov-Poisson system}
\begin{exa}\label{ex:weak} 
We first consider the 2D2V weak Landau damping with initial condition 
\begin{equation}
\label{eq:weak}
f(\bx,\bv,t=0) =\frac{1}{(2 \pi)^{d / 2}} \left(1+\alpha \sum_{m=1}^{d} \cos \left(k x_{m}\right)\right)\exp\left(-\frac{|\bv|^2}{2}\right),
\end{equation}
where $d=2$, $\alpha=0.01$, and $k=0.5$. We set the computation domain as $[0,L_x]^2\times[-L_v,L_v]^2$, where $L_x=\frac{2\pi}{k}$ and $L_v=6$.  We simulate the problem using approach I and report the numerical results in Figure \ref{fig:weak}. It is observed that the low-rank method is able to  predict the correct damping rate of the electric energy and capture the low-rank structure of the solution. The hierarchical ranks remain low, leading to significant efficiency of the low-rank method. CPU time is $76.8s$, $117.5s$, and $265.6s$ for mesh size $N_x^2\times N_v^2=16^2\times32^2,\,32^2\times64^2,\,64^2\times128^2$, respectively. Similar to Example \ref{ex:adv}, the CPU time scales linearly with respect to $N$, as the proposed method is able to efficiently capture the low-rank solution structures.  The time histories of the relative error in the total partial number and the total energy are plotted in Figure \ref{fig:weak_his}. It is observed that the proposed method is able to conserve the physical invariants up to the scale of the truncation threshed  $\varepsilon$.

%In particular, we consider the following benchmark test including Landau damping and two stream instability. The performance of the proposed methods will be assessed in the aspects of the ability to resolve complex solution structure and to recover physical properties of interest, and the efficiency of computation, especially in high dimensions.
%For the nonlinear VP system, we do not apply the boundary fix for the approach II, since the performance only improves slightly while the cost increases a lot due to the costly interpolation of the electric field at the boundary points.

\begin{figure}[h!]
	\centering
	\subfigure[]{\includegraphics[height=60mm]{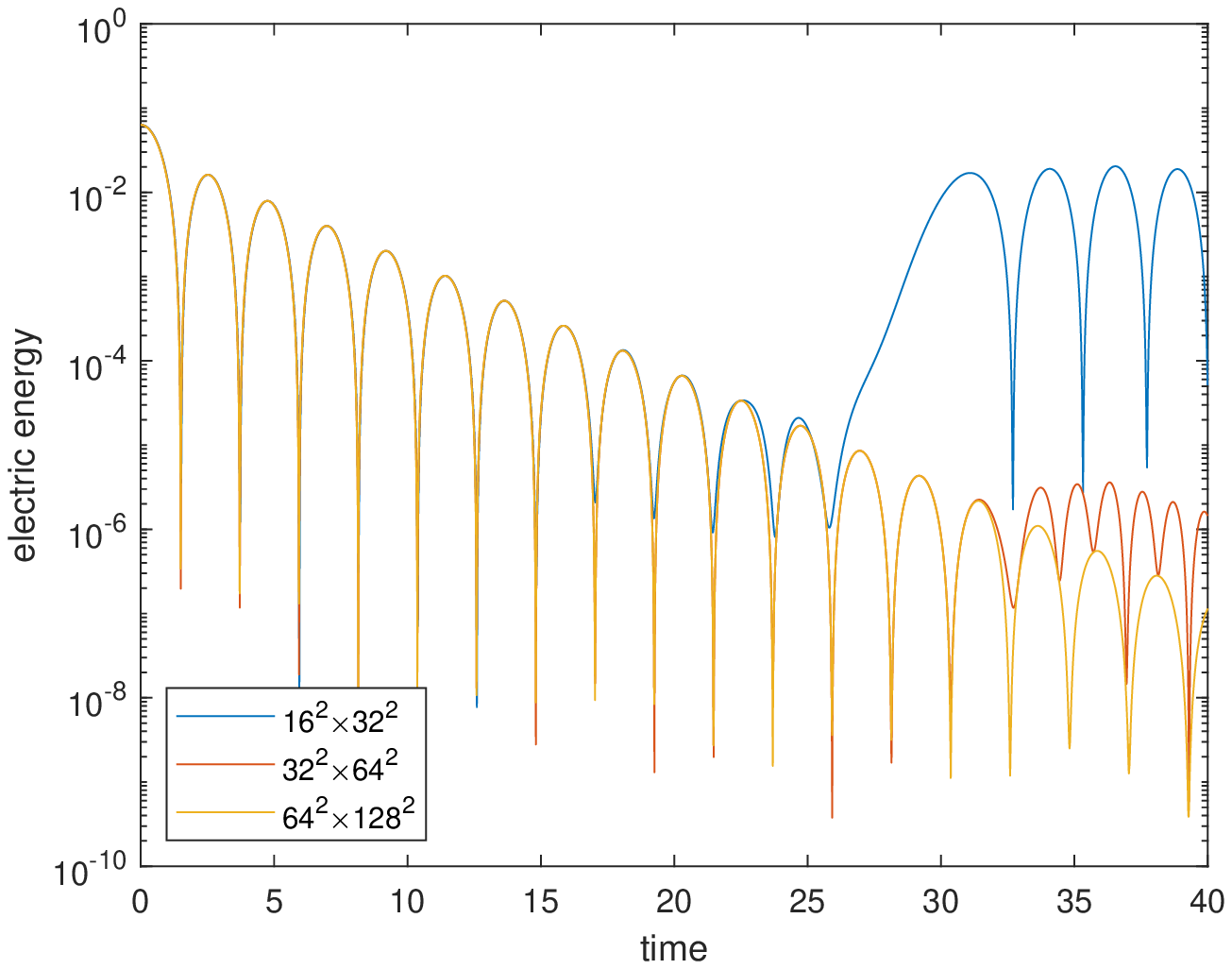}}
	\subfigure[]{\includegraphics[height=60mm]{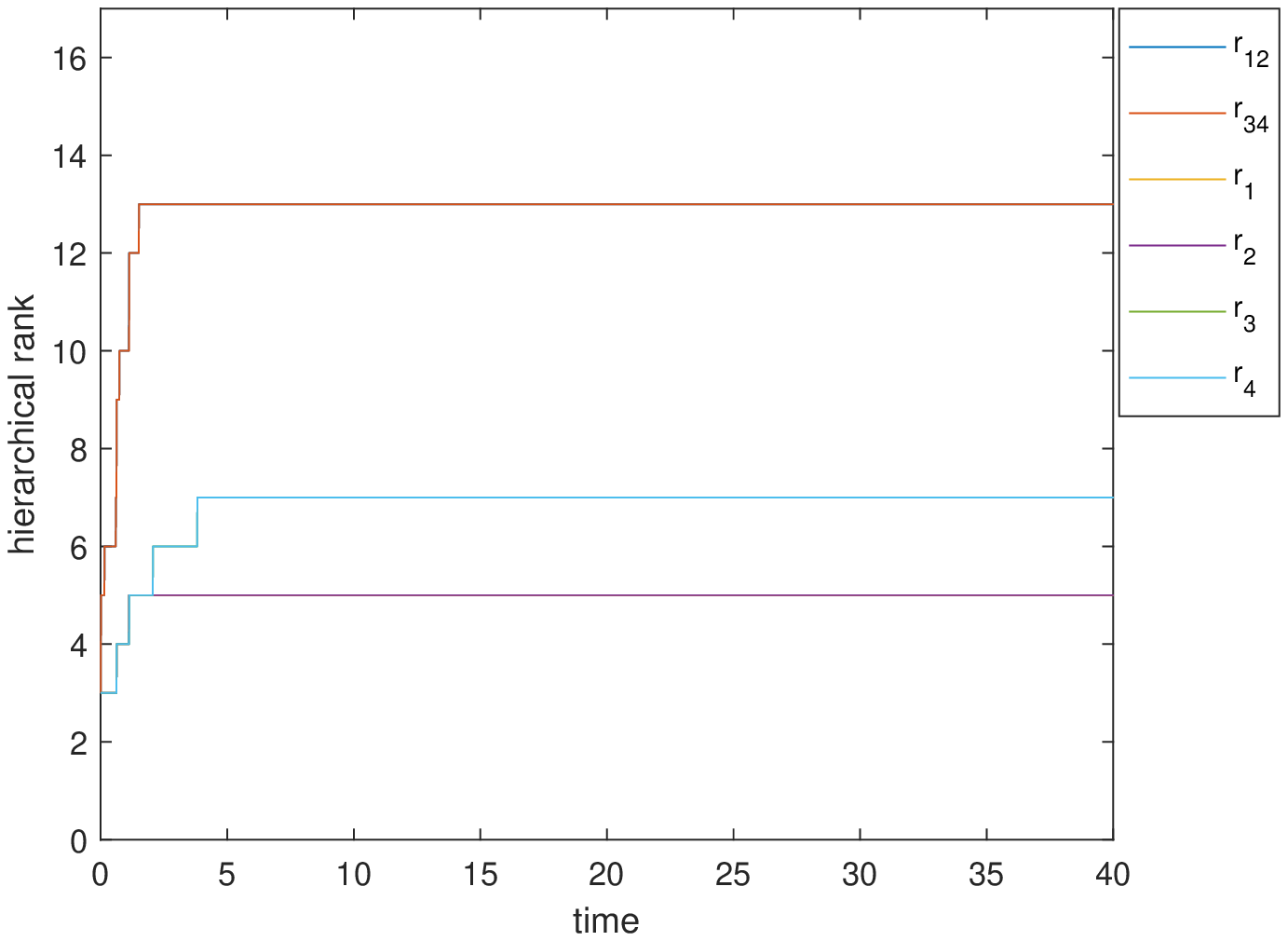}}
	%\subfigure[]{\includegraphics[height=40mm]{weak2d_total_energy.eps}}
	\caption{Example \ref{ex:weak}. Weak Landau damping 2D2V. $d=2$. $\varepsilon=10^{-6}$. Approach I. (a) The time evolution of  the electric energy. (b) Hierarchical ranks for mesh $N_x^2\times N_v^2=64^2\times 128^2$. For uniform meshes $N_x^2\times N_v^2=16^2\times32^2,\,32^2\times64^2,\,64^2\times128^2$. \label{fig:weak}}
\end{figure}

\begin{figure}[h!]
	\centering
	\subfigure[]{\includegraphics[height=60mm]{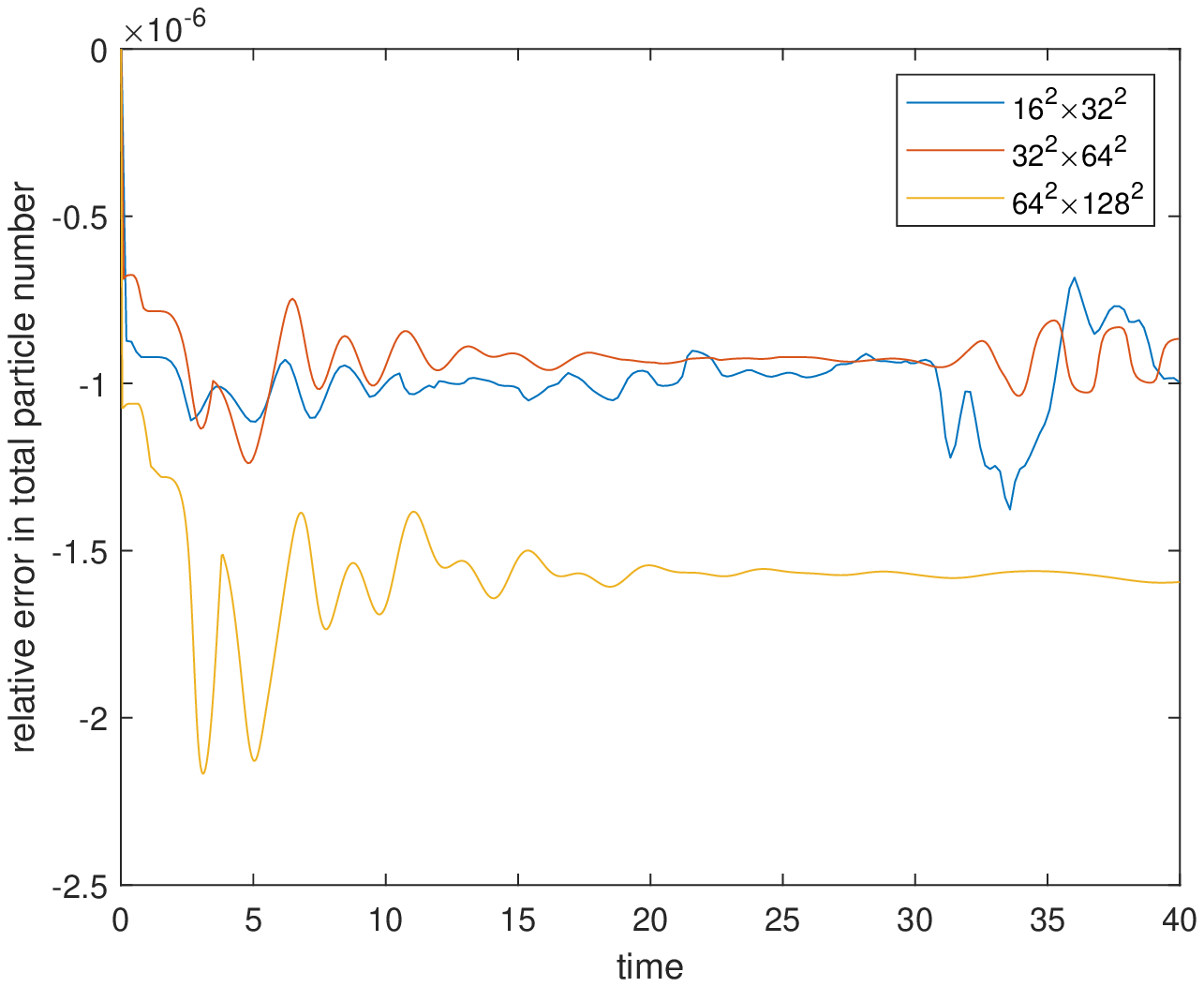}}
	\subfigure[]{\includegraphics[height=60mm]{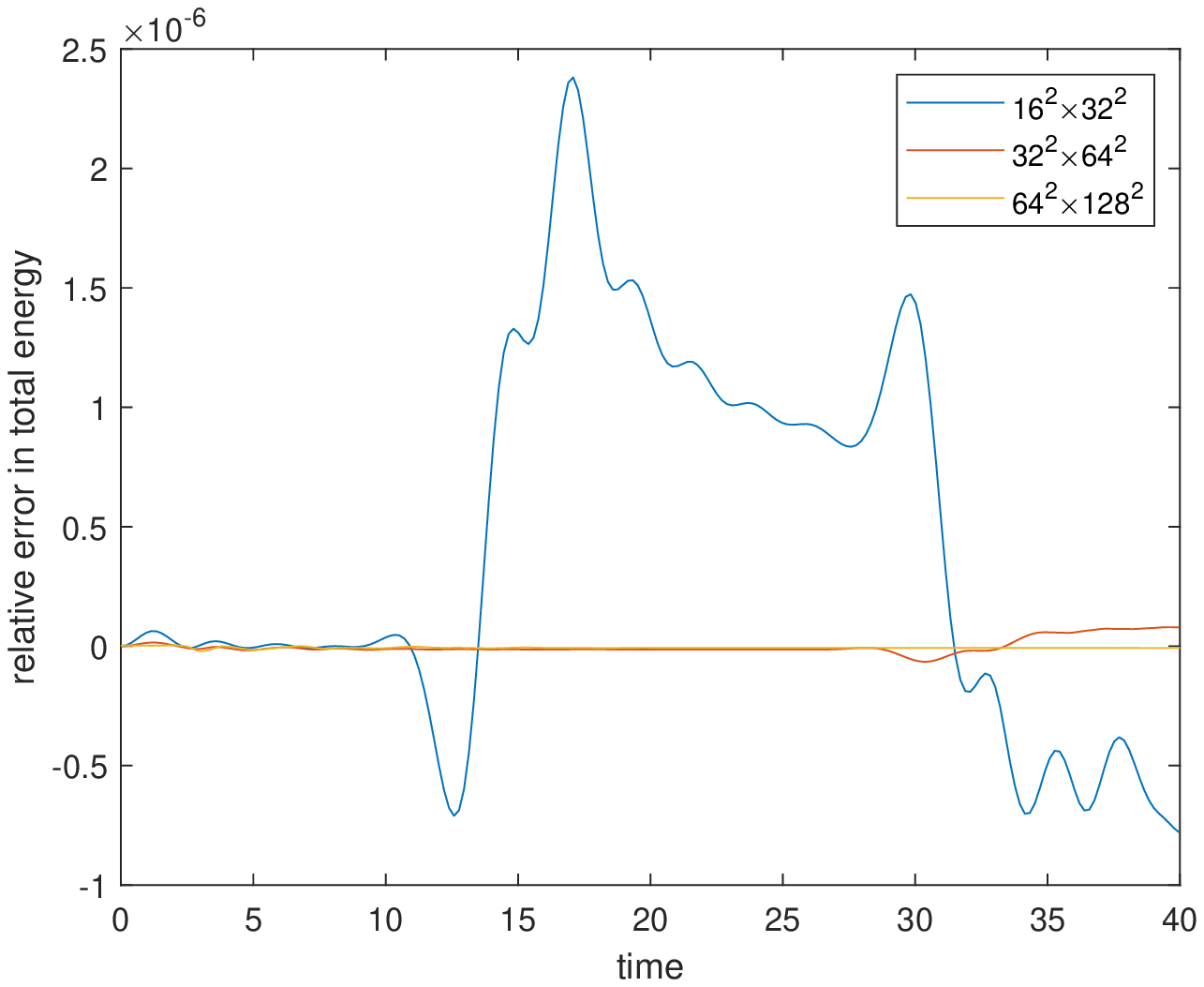}}
	%\subfigure[]{\includegraphics[height=40mm]{weak2d_total_energy.eps}}
	\caption{Example \ref{ex:weak}. Weak Landau damping 2D2V. $d=2$. $\varepsilon=10^{-6}$. Approach I.  The time evolution of  the relative errors. (a) total particle number, (b) total energy. \label{fig:weak_his}}
\end{figure}

\end{exa}

\begin{exa}\label{ex:strong} In this example, we consider the strong Landau damping with the initial condition \eqref{eq:weak} and $d=2$. The parameters are the same as in previous example except $\alpha=0.5$. Unlike the weak case, the nonlinear effect would play a dominant role and the solution would gradually develop filamentation structures. Note that the VP system enjoys a well-known time reversibility property; that is we evolve the VP system from initial condition $f(\bx,\bv,t=0)$ to $t=T$ and flip the velocity variable, i.e., $\tilde{f}(\bx,\bv,T) = f(\bx,-\bv,T)$, and then continue the evolution of the VP system for $\tilde{f}$ up to $t=2T$ and obtain $\tilde{f}(\bx,\bv,2T)$. 
Then,  $\tilde{f}(\bx,\bv,2T)=f(\bx,\bv,t=0)$. We can make use of the time reversibility property for the accuracy test. In particular, we let $T=1$ and compare the approximation of  $\tilde{f}(\bx,\bv,2T)$ with the initial condition. The truncation threshold is set to be $\varepsilon=10^{-6}$. The convergence study is summarized in Table \ref{tb:strong2d}. Second order of convergence is observed. Then we consider the long term simulation of the strong Landau damping.  We set $\varepsilon=10^{-3}$ and $r_{\max}=32$, and compute the solution up to $T=30$. In Figure \ref{fig:strong_2d_ranks}, we plot the time histories of the hierarchical ranks with three mesh sizes  $N_x^2\times N_v^2 = 32^2\times64^2$, $64^2\times128^2$, and $128^2\times256^2$. It is observed that the hierarchical ranks increase over time to capture the underlying filamentation structures before reaching $r_{max}$. The CPU time is $438.9s$, 
$1278.7s$, and $2606.4s$ which scales linearly with $N$. Note that for this example the total computational cost is dominated by the hierarchical HOSVD truncation. We set the same $r_{max}$ and hence the cost for truncation is similar for the three meshes.
In Figure \ref{fig:strong_2dcuts}, we plot the 2D cuts of the solutions at $(x_2,v_2)=(2\pi,0)$ and $(x_1,x_2)=(2\pi,2\pi)$ at $t=5,\ 15,\ 30$ with mesh size $N_x^2\times N_v^2 = 128^2\times256^2$.   It is observed  that the proposed method is able to capture the main structure of the solution. In Figure \ref{fig:strong2d_his}, we report the time histories of the electric energy and the relative errors in total mass and energy. 

\begin{table}[!hbp]
	\centering
	\caption{Example \ref{ex:strong}. $d=2$. $T=1$. $\varepsilon=10^{-6}$. Approach I. }
	\label{tb:strong2d}
	\begin{tabular}{|c|c|c|}
		\hline
		$N_x^2\times N_v^2$ & L$^2$-error & order  \\\hline
$16^2\times32^2$&	5.19E-02		&			\\	\hline
$32^2\times64^2$	&	3.62E-03	&	3.84		\\	\hline
$64^2\times128^2$	&	5.11E-04	&	2.82		\\	\hline
$128^2\times256^2$	&	1.57E-04	&	1.71		\\	\hline
	\end{tabular}
\end{table}

\begin{figure}[h!]
	\centering
	\subfigure[]{\includegraphics[height=40mm]{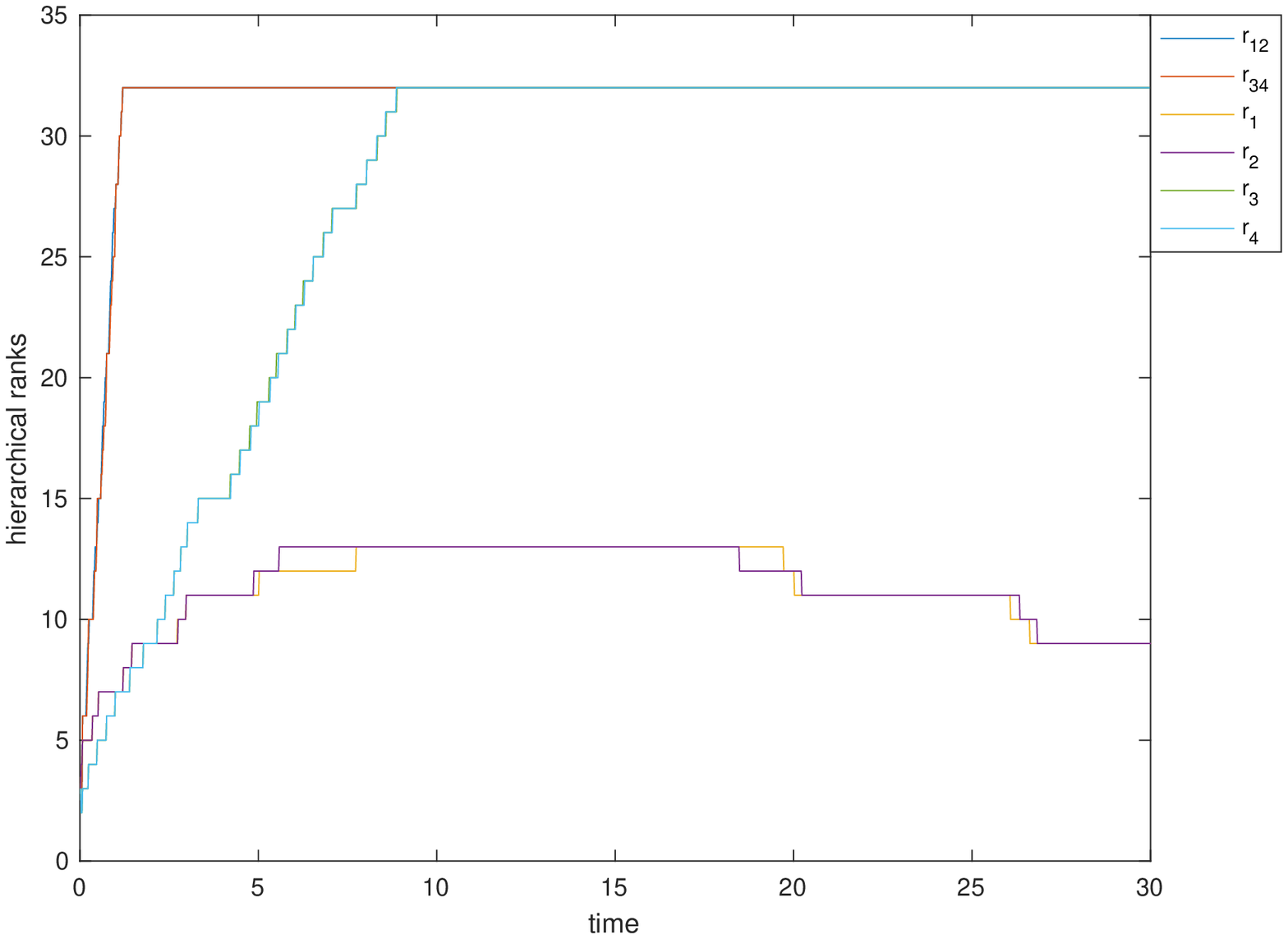}}
	\subfigure[]{\includegraphics[height=40mm]{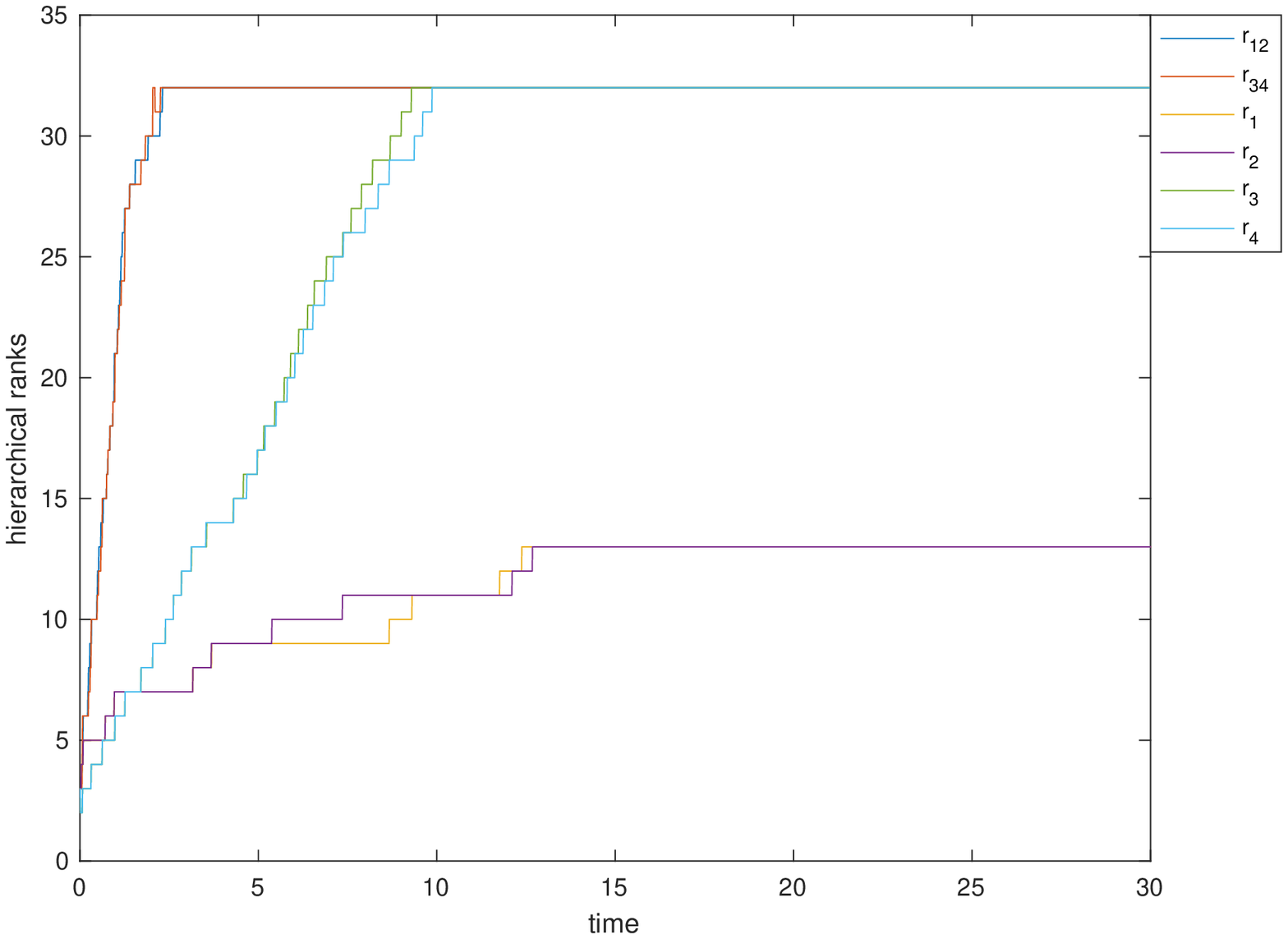}}
	\subfigure[]{\includegraphics[height=40mm]{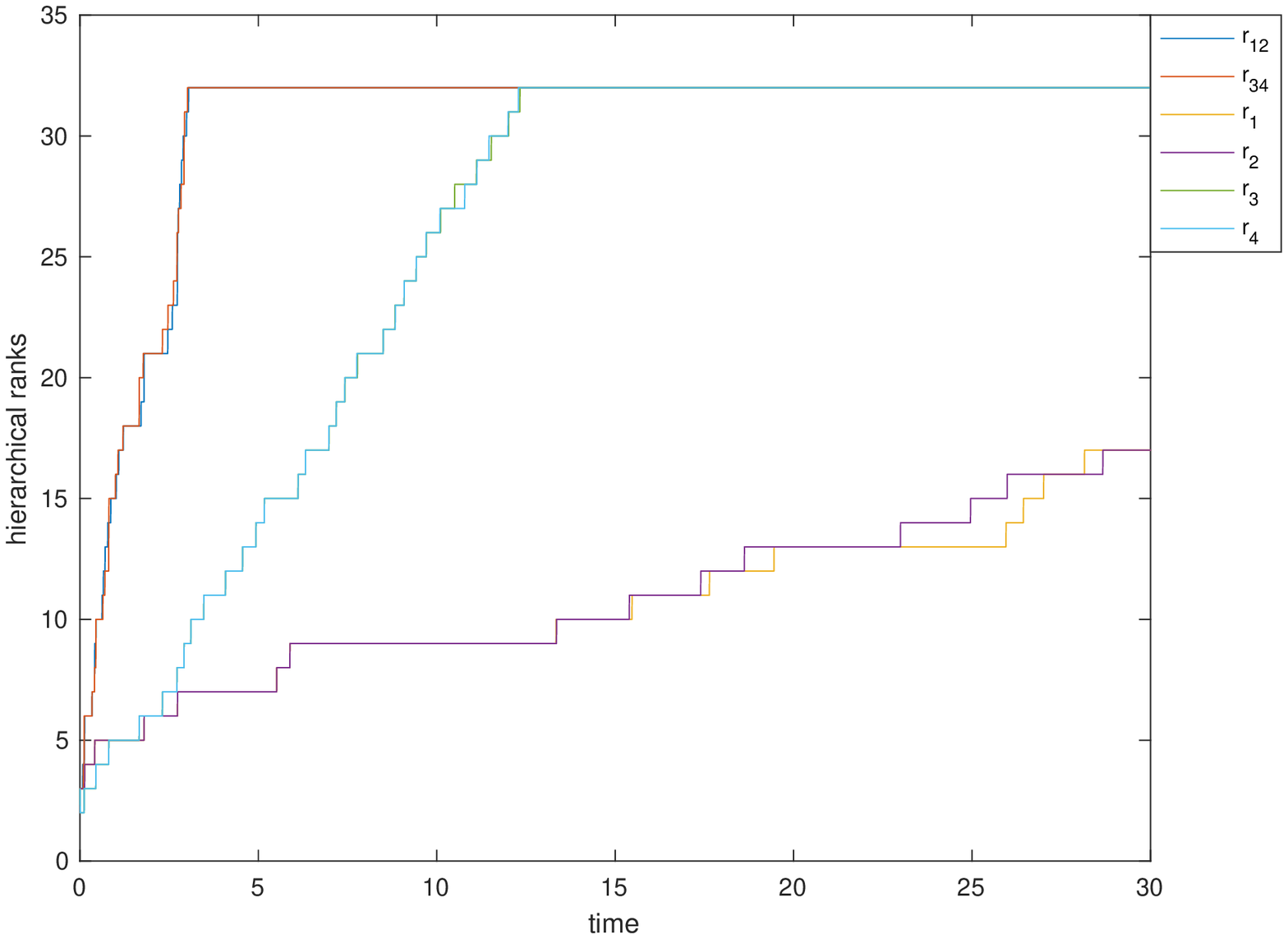}}
	%\subfigure[]{\includegraphics[height=40mm]{weak2d_total_energy.eps}}
	\caption{Example \ref{ex:strong}. Strong Landau damping 2D2V. Hierarchical ranks. $d=2$. $\varepsilon=10^{-3}$. $r_{max}=32$. Approach I. (a) $N_x^2\times N_v^2=32^2\times64^2$. (b) $N_x^2\times N_v^2=64^2\times128^2$. (c) $N_x^2\times N_v^2=128^2\times256^2$.
	\label{fig:strong_2d_ranks}}
\end{figure}

\begin{figure}[h!]
	\centering
	\subfigure[]{\includegraphics[height=60mm]{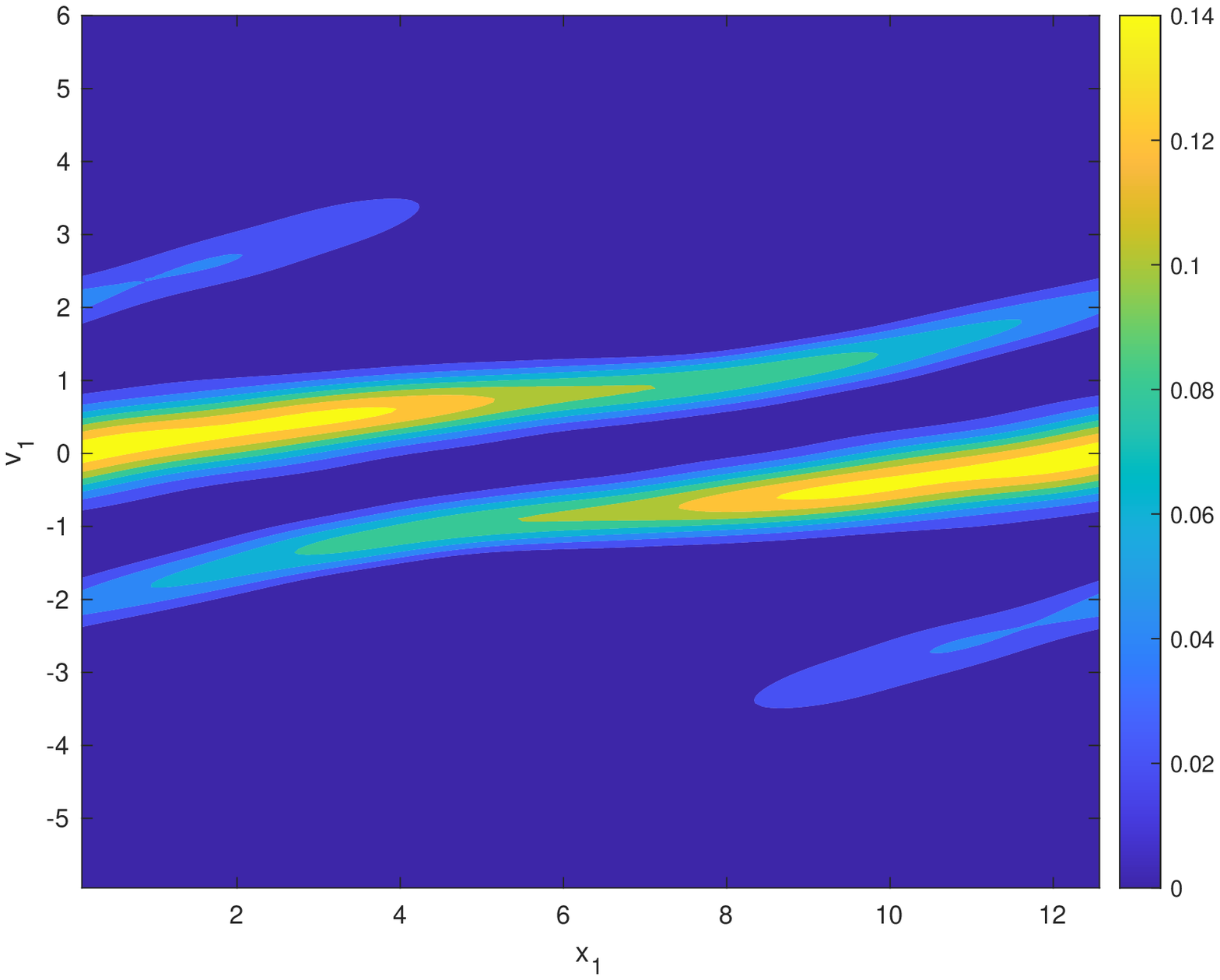}}
	\subfigure[]{\includegraphics[height=60mm]{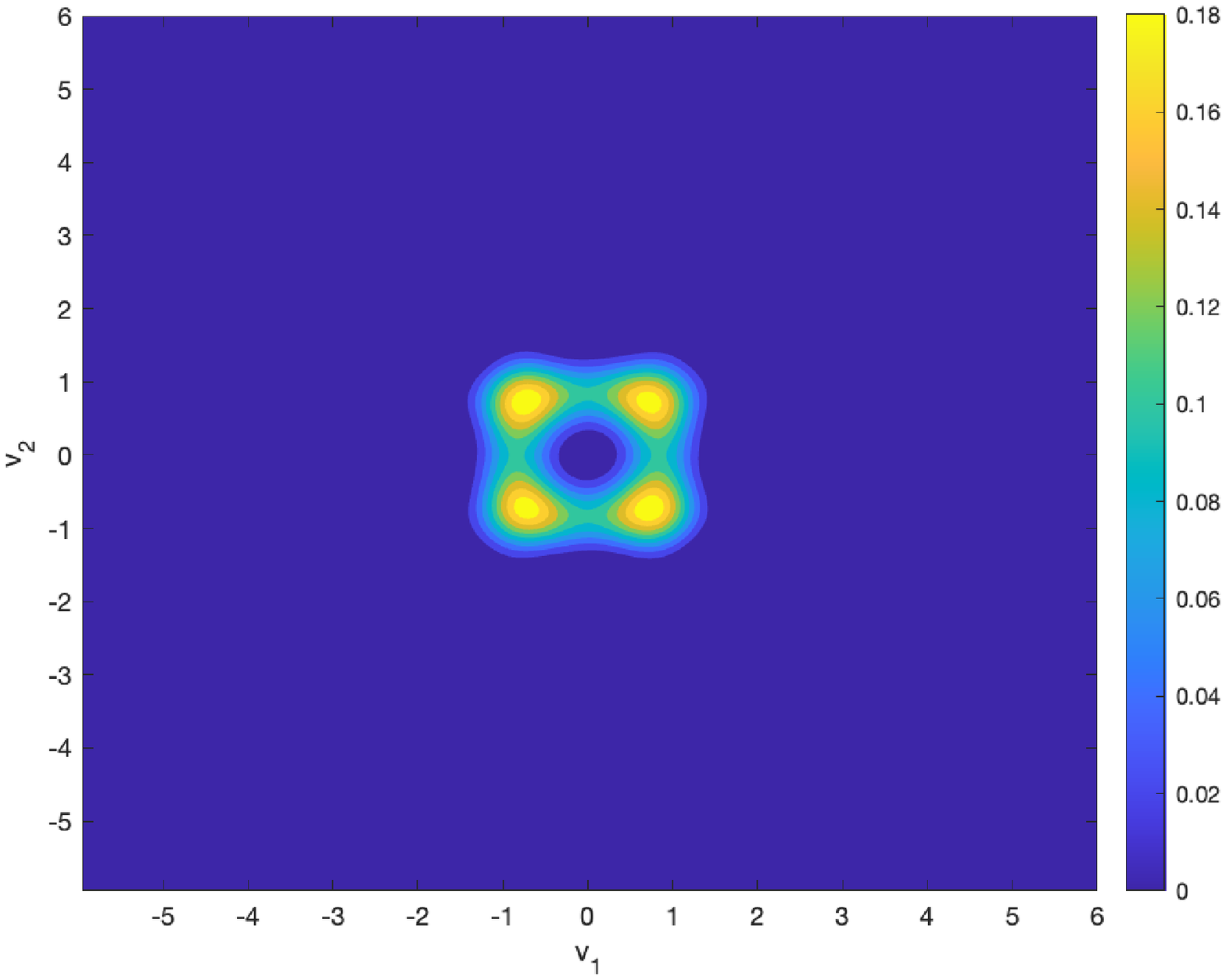}}
	\subfigure[]{\includegraphics[height=60mm]{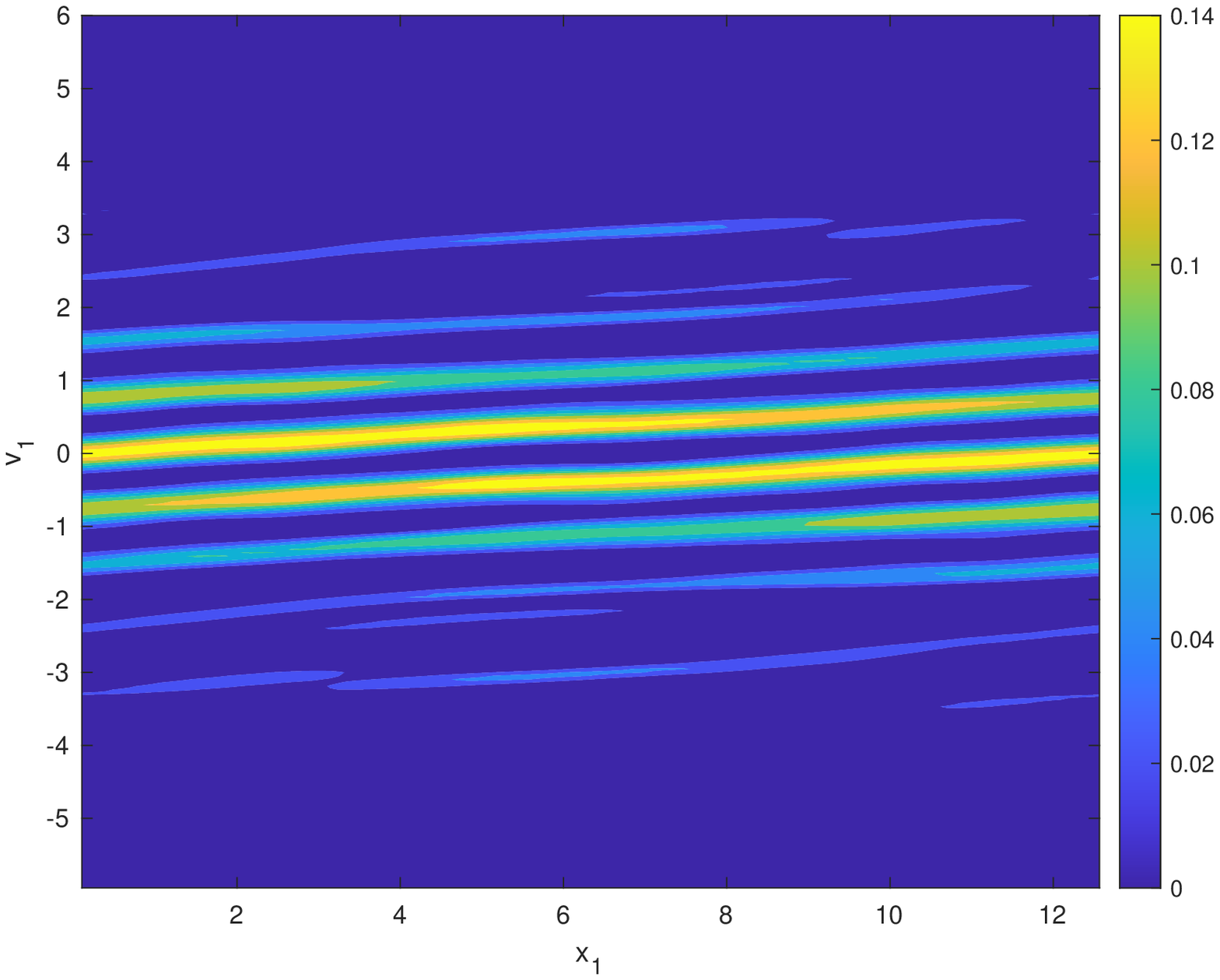}}
	\subfigure[]{\includegraphics[height=60mm]{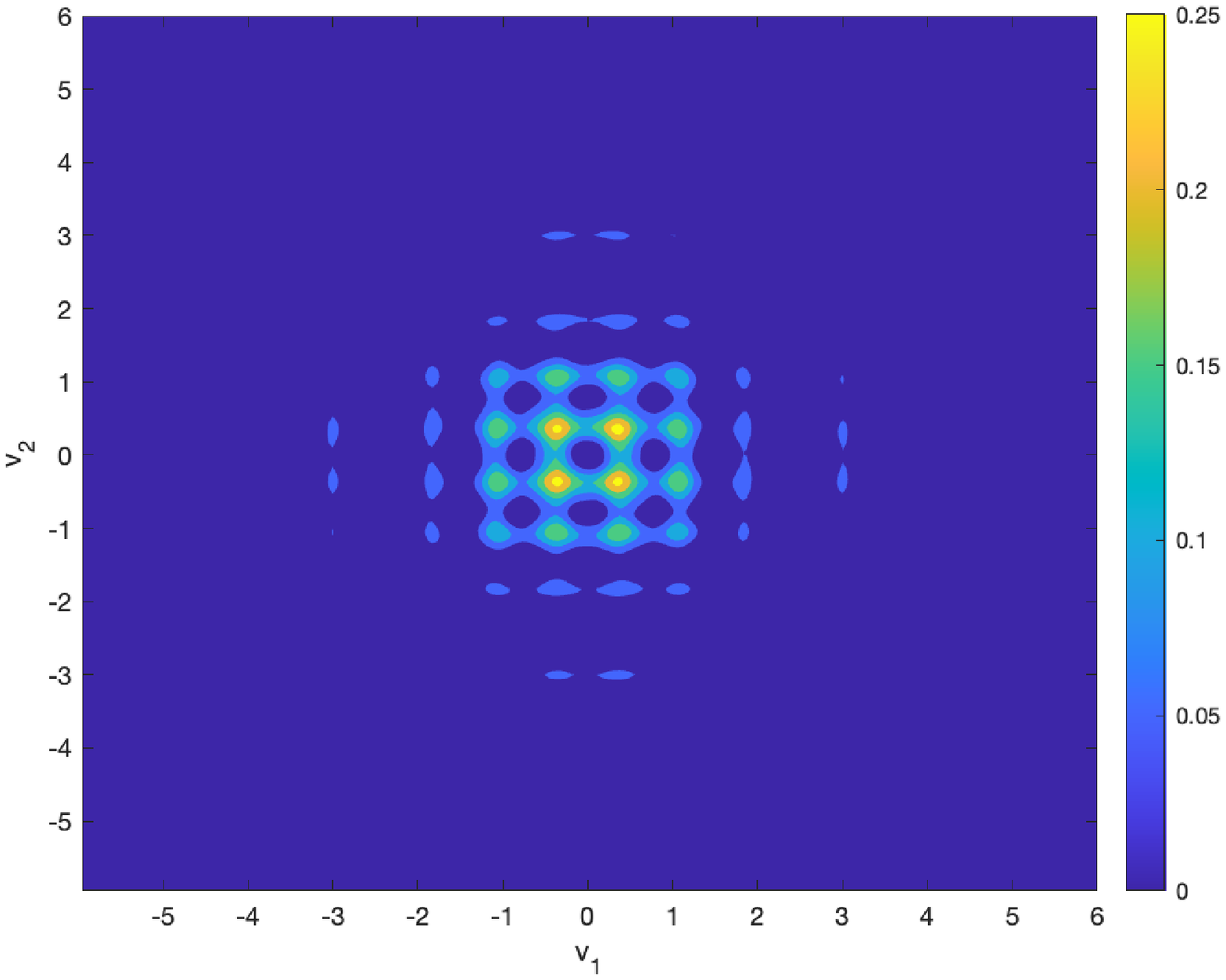}}
	\subfigure[]{\includegraphics[height=60mm]{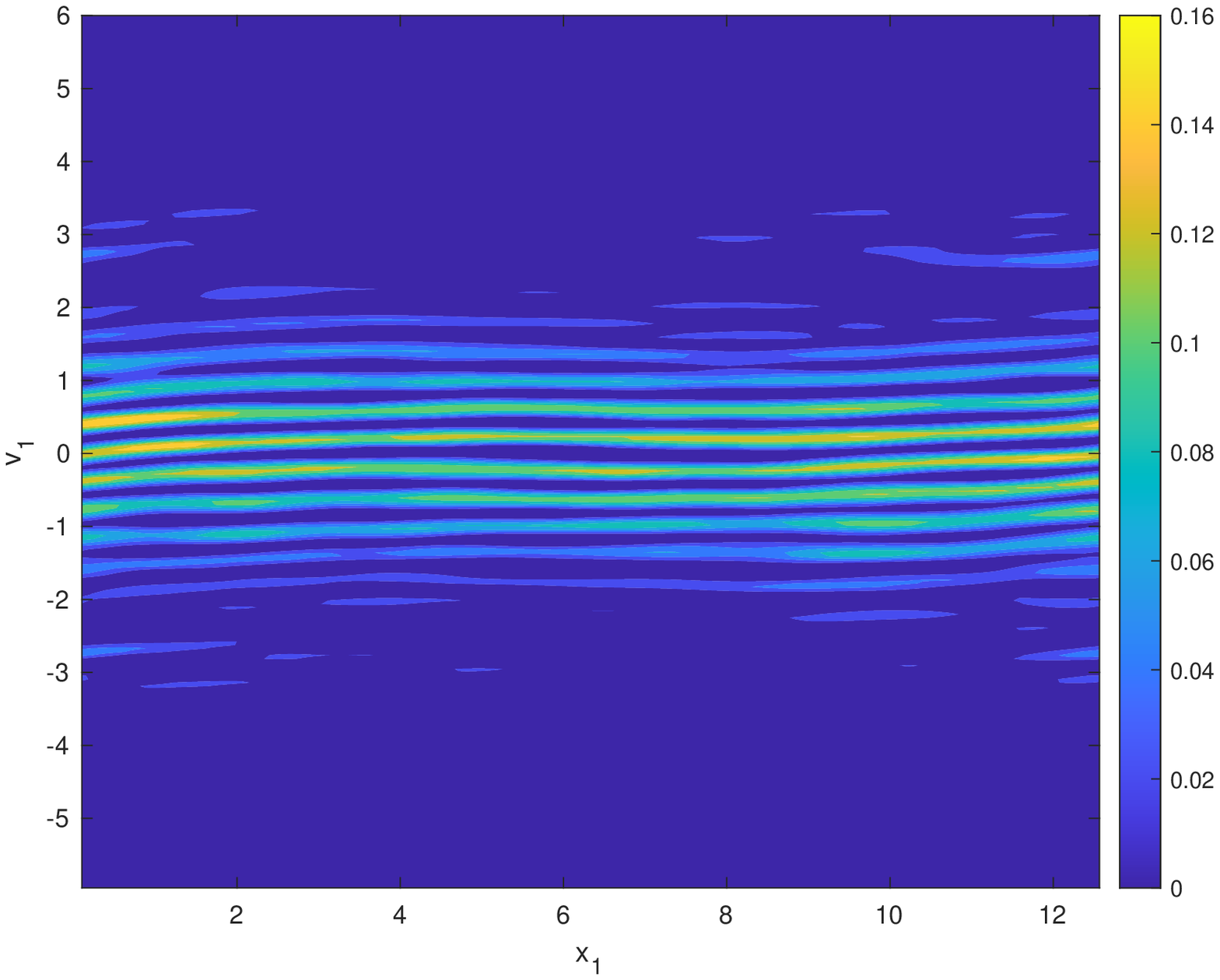}}
	\subfigure[]{\includegraphics[height=60mm]{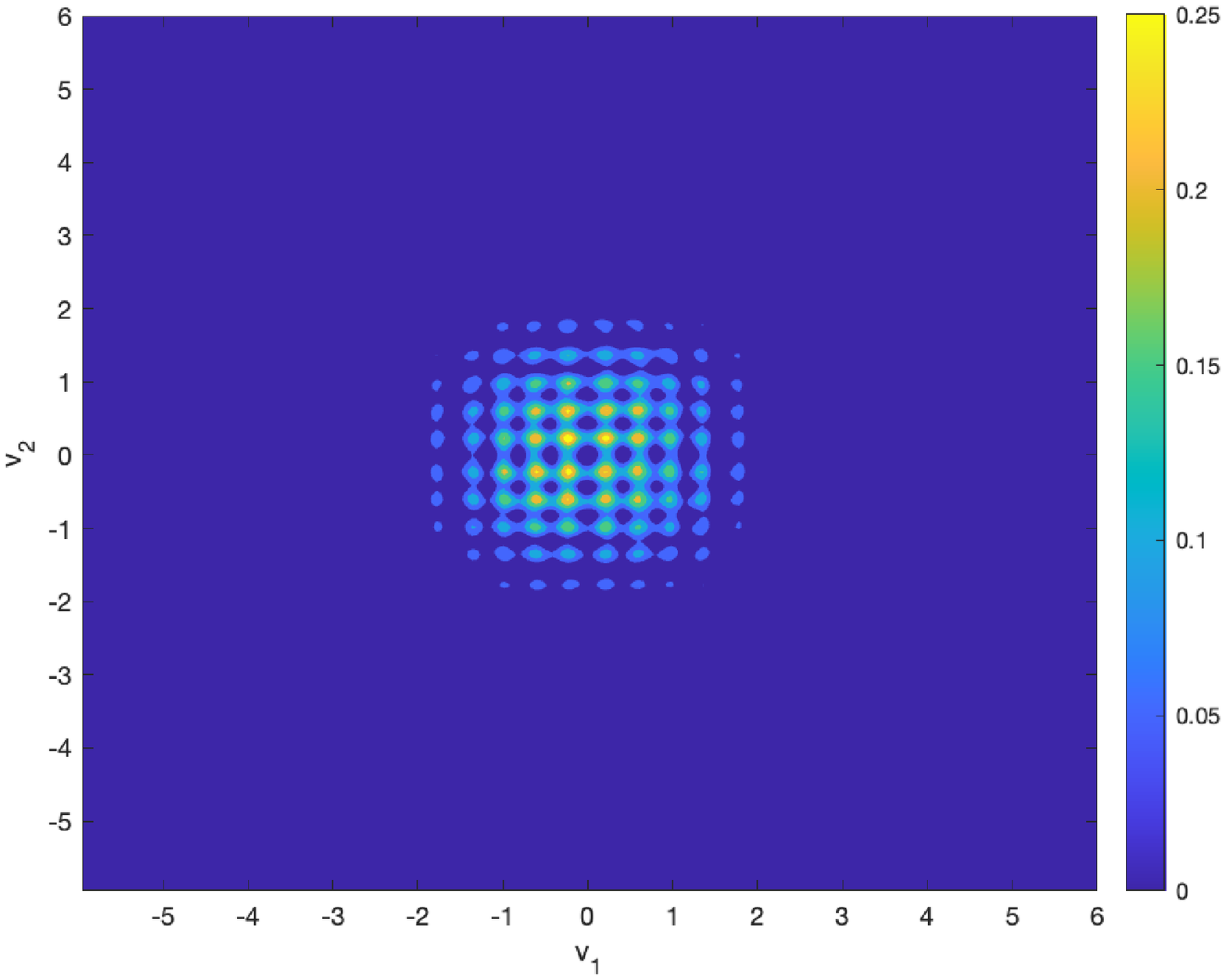}}
	%\subfigure[]{\includegraphics[height=40mm]{weak2d_total_energy.eps}}
	\caption{Example \ref{ex:strong}. Strong Landau damping 2D2V. $d=2$. $\varepsilon=10^{-3}$. $r_{max}=32$. Uniform meshes $N_x^2\times N_v^2=128^2\times256^2$. Approach I. (a) 2D cut at $(x_2,v_2)=(2\pi,0)$ $t=5$. 
(b) 2D cut at $(x_1,x_2)=(2\pi,2\pi)$ $t=5$	. (c) 2D cut at $(x_2,v_2)=(2\pi,0)$ $t=15$. (d) $(x_1,x_2)=(2\pi,2\pi)$ $t=15$. (e) 2D cut at $(x_2,v_2)=(2\pi,0)$ $t=30$. (f) $(x_1,x_2)=(2\pi,2\pi)$ $t=30$.
	\label{fig:strong_2dcuts}}
\end{figure}

\begin{figure}[h!]
	\centering
	%\subfigure[]{\includegraphics[height=40mm]{landau2d_strong_elec_energy.eps}}
	\subfigure[]{\includegraphics[height=60mm]{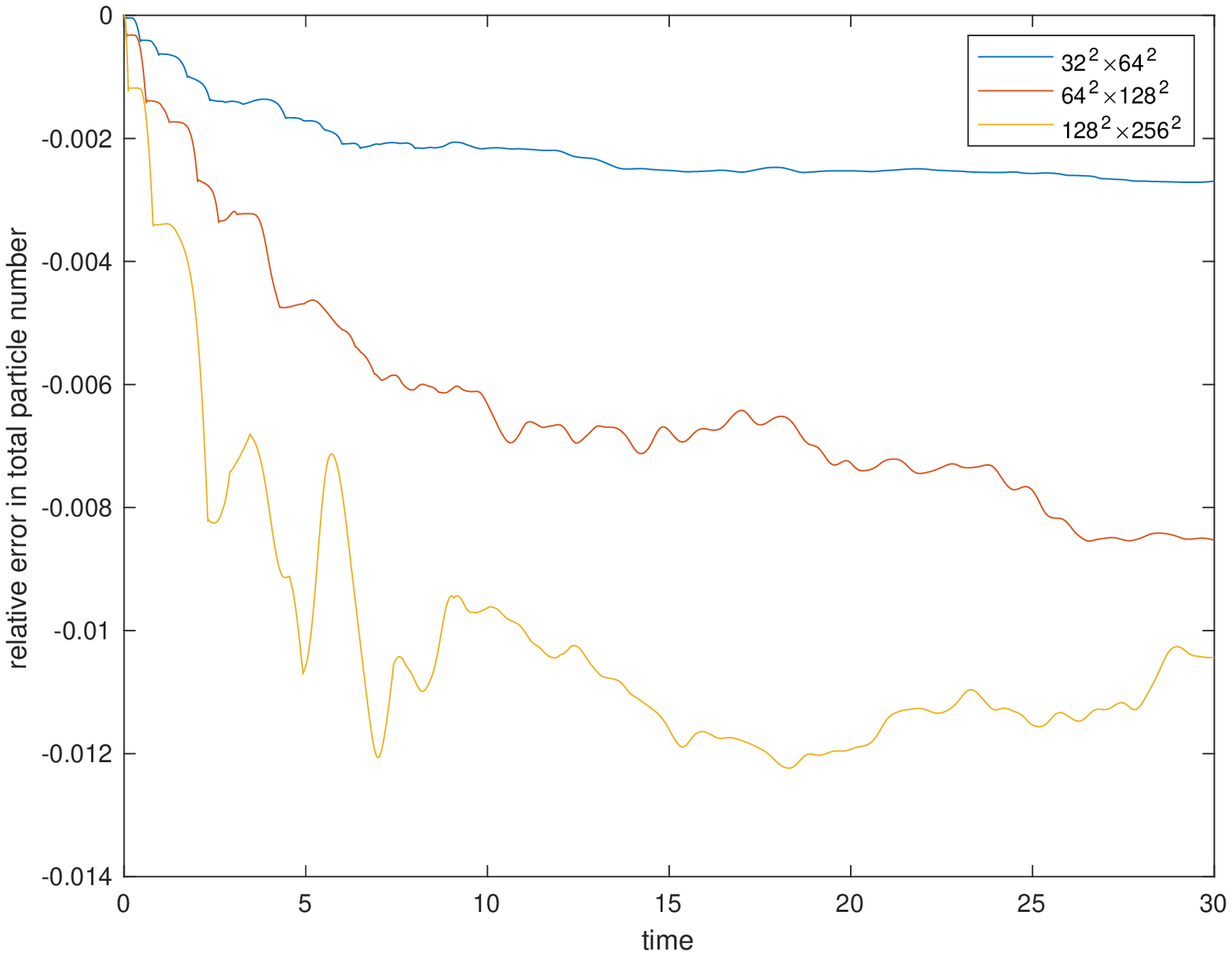}}
	\subfigure[]{\includegraphics[height=60mm]{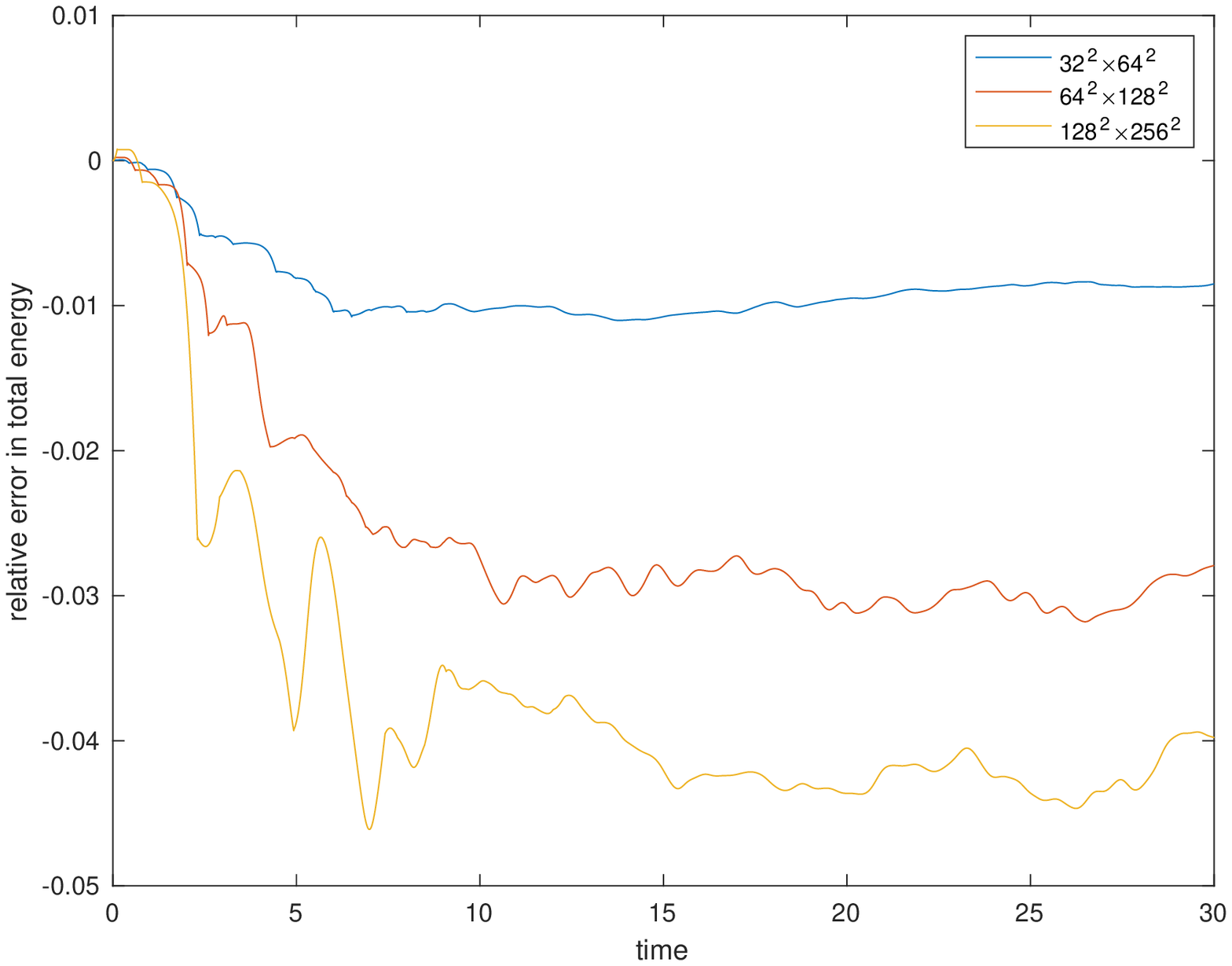}}
	%\subfigure[]{\includegraphics[height=40mm]{weak2d_total_energy.eps}}
	\caption{Example \ref{ex:strong}. Strong Landau damping 2D2V. $d=2$. $\varepsilon=10^{-3}$. $r_{\max}=32$. Approach I. (a) The time evolution of  the relative errors in total particle number. (b) The time evolution of  the relative errors in total energy.\label{fig:strong2d_his}}
\end{figure}

\end{exa}

%% file: conclusion.tex
\section{Conclusion}
\label{sec4}
\setcounter{equation}{0}
\setcounter{figure}{0}
\setcounter{table}{0}

In this paper, we proposed a novel low-rank tensor approach to approximate transport equations in high dimensions with application to Vlasov simulations. In particular, the solution is represented in the low-rank HT format, and the associated basis is  dynamically and adaptively updated by the proposed adding and removing basis procedure. High order spatial and temporal discretizations are employed for accurate capture of complex solution structures. For the transport problems that do not exhibit low-rank structures, we further propose to solve the associated flow maps in a similar low-rank fashion, which may  enjoy the desired low-rank structures. We plan to extend the approach to the Vlasov-Maxwell system and other kinetic models with relaxation/collision terms, such as  the BGK model. We also plan to address the open problems associated with the novel flow map approach, such as extension to general boundary condition.